\documentclass[12pt]{alggeom}
\usepackage[]{amsmath, amsthm, amsfonts, verbatim, amssymb, mathrsfs, mathtools, booktabs, multirow}

\usepackage[all]{xy}

\usepackage{hyperref}
\usepackage{pifont}

\usepackage{enumerate}
\usepackage{graphicx}
\usepackage{epstopdf}
\DeclareMathAlphabet{\mathpzc}{OT1}{pzc}{m}{it}
\usepackage{tikz,tikz-cd, color}
\usepackage[makeroom]{cancel}
\usetikzlibrary{arrows}
\usepackage{mathdots}
\usepackage{float}
\usepackage{caption, subcaption}
\usepackage{rotating}
\usepackage{adjustbox}
\usetikzlibrary{matrix}

\tikzset{
  symbol/.style={
    draw=none,
    every to/.append style={
      edge node={node [sloped, allow upside down, auto=false]{$#1$}}}
  }
}

\usetikzlibrary{decorations.shapes}

\tikzset{decorate sep/.style 2 args=
{decorate,decoration={shape backgrounds,shape=circle,shape size=#1,shape sep=#2}}}

\newcommand{\nospaceperiod}{\makebox[0pt][l]{\,.}}

\newtheorem{theorem}{Theorem}[section]
\newtheorem{lemma}[theorem]{Lemma}
\newtheorem{proposition}[theorem]{Proposition}

\newtheorem{conjecture}[theorem]{Conjecture}
\newtheorem{definition}[theorem]{Definition}

\newtheorem{question}[theorem]{Question}

\theoremstyle{remark}
\newtheorem{remark}[theorem]{Remark}
\newtheorem{example}[theorem]{Example}

\theoremstyle{definition}
\newtheorem{notation}[theorem]{Notation}

\newcommand{\bC}{\mathbb{C}}

\newcommand{\bF}{\mathbb{F}}

\newcommand{\bQ}{\mathbb{Q}}
\newcommand{\bP}{\mathbb{P}}
\newcommand{\bR}{\mathbb{R}}
\newcommand{\bZ}{\mathbb{Z}}
\newcommand{\bfE}{\mathbf{E}}

\newcommand{\cO}{\mathcal{O}}

\newcommand{\cQ}{\mathcal{Q}}
\newcommand{\cS}{\mathcal{S}}

\newcommand{\fc}{\mathfrak{c}}
\newcommand{\fq}{\mathfrak{q}}

\newcommand{\sA}{\mathcal{A}}
\newcommand{\sB}{\mathcal{B}}
\newcommand{\sC}{\mathcal{C}}
\newcommand{\sF}{\mathcal{F}}
\newcommand{\sE}{\mathcal{E}}
\newcommand{\sI}{\mathcal{I}}
\newcommand{\sL}{\mathcal{L}}
\newcommand{\sV}{\mathcal{V}}

\newcommand{\tbP}{\widetilde{\mathbb{P}}}

\newcommand{\tY}{\widetilde{Y}}

\newcommand{\id}{\mathrm{id}}

\newcommand{\Mov}{\mathrm{Mov}}

\newcommand{\cMov}{\overline{\mathrm{Mov}}}
\newcommand{\Bir}{{\rm Bir}}
\newcommand{\Aut}{{\rm Aut}}

\def\Bs{{\bf Bs}}

\newcommand{\Bl}{{\rm Bl}}

\newcommand{\coker}{{\rm coker}\,}
\newcommand{\Eff}{\mathrm{Eff}}

\newcommand{\Exc}{{\rm Exc}}
\newcommand{\Ext}{{\rm Ext}}

\newcommand{\ev}{{\rm ev}}
\newcommand{\Gr}{{\rm  Gr}}

\newcommand{\Hom}{{\rm Hom}}

\newcommand{\sHom}{\mathscr{H}om}

\newcommand{\Nef}{{\rm Nef}}

\newcommand{\GL}{{\rm GL}}

\newcommand{\Pic}{{\rm Pic}\,}
\newcommand{\Proj}{{\bf Proj}}
\newcommand{\rank}{{\rm rank}}

\newcommand{\Sing}{{\rm Sing}}
\newcommand{\Sym}{{\rm Sym}}

\def\<{\langle}
\def\>{\rangle}

\def\ra{\rightarrow}
\def\dra{\dashrightarrow}

\author{Ching-Jui Lai}
\email{cjlai72@mail.ncku.edu.tw}
\address[]{Department of Mathematics, National Cheng Kung University, Tainan 70101, Taiwan
}

\author{Sz-Sheng Wang}
\email{sswangtw@icloud.com}
\address[]{Shing-Tung Yau Center of Southeast University, Southeast University,
Nanjing 211189, China}

\classification{14J32, 14E05}
\keywords{Calabi-Yau threefold, birational geometry, movable cone}
\thanks{The first author is supported by the grant MOST 107-2115-M-006-020 and an internal grant of National Cheng Kung University. The second author is partially supported by the Fundamental Research Funds for the Central Universities 2242020R10048, and he thanks Southeast University, Shing-Tung Yau Center of Southeast University, Tsinghua University, and Yau Mathematical Sciences Center for providing support and a stimulating environment, and also thanks to the math department of National Cheng Kung University for its hospitality. Some of the work on this paper was done while he was visiting the first author at NCKU. We thank for referees for their very useful comments.}

\begin{document}
\title[The movable cone of certain CY3 of Picard number two]{The movable cone of certain Calabi--Yau threefolds of Picard number two}

% Enter the publication year and the ID number of the paper

\begin{abstract}
We describe explicitly the chamber structure of the movable cone for a general smooth complete intersection Calabi--Yau threefold $X$ of Picard number two in certain $\mathbb{P}^{r}$-ruled Fano manifold and hence verify the Morrison--Kawamata cone conjecture for such $X$. Moreover, all birational minimal models of such Calabi--Yau threefolds are found, whose number is finite up to isomorphism.
\end{abstract}

\maketitle

\tableofcontents

%%%%%%%%%%%%%%%%%%%%
%%%%%%%%%%%%%%%%%%%%
\section{Introduction}\label{Introsec}
A smooth projective variety $X$ of dimension $n$ is called a Calabi--Yau $n$-fold if it satisfies $\omega_X:=\wedge^n\Omega_X \cong \cO_X$ and $H^1 (\cO_X) = 0$. Such manifolds are fundamental objects in birational geometry and theoretical physics.

To understand the geometry of a variety, one considers linear systems of different divisors. Modulo numerical equivalence, this leads to the study of convex geometry of various cones of divisors in the N\'eron--Severi space. For example, being the dual of the Mori cone of curves, the nef cone of divisors plays essential roles in the cone theorem \cite[Theorem 3.7]{KM98}. Another crucial example is the movable cone of divisors, which encodes the birational geometry of a given variety \cite{Kawamata88}.

Inspired by mirror symmetry of Calabi--Yau manifolds, Morrison \cite{Morrison93} and Kawamata \cite{Kawamata97} proposed the conjectures which would give a clear picture of relevant cones for Calabi--Yau manifolds. To be more precise, let $N^1 (X)$ be the N\'eron--Severi group, generated by the classes of the divisors on $X$ modulo numerical equivalence. Inside the N\'eron--Severi space $N^1 (X)_{\bR}= N^1 (X) \otimes_{\bZ} \bR$ we have the effective cone $\Eff (X)$, the nef cone $\Nef (X)$, and the movable cone $\cMov (X)$ (that is, the closure of the convex hull of movable divisor classes). Recall that a divisor $D$ is \emph{movable} if the linear system $|m D|$ has no fixed component for some positive integer $m$. As usual, $\Bir (X)$ denotes the group of birational automorphisms of $X$. Notice that every $g \in \Bir (X)$ of the Calabi--Yau manifold $X$ is an isomorphism in codimension $1$ by negativity lemma \cite[Lemma 3.39]{KM98}. Thus, there is an induced homomorphism
$$r:\Bir(X)\ra\GL(N^1(X)),\ g\mapsto g^*.$$
Moreover, if $D$ is movable (resp.~effective), then $g^{\ast} D$ is again movable (resp.~effective).

For our purpose, we state the movable cone conjecture as follows (and a similar statement can be made for the action of $\Aut (X)$ on $\Nef (X) \cap \Eff (X)$):

\begin{conjecture}\label{coneconjM}
Let $X$ be a Calabi--Yau manifold. The action of $\Bir (X)$ on the movable effective cone $\cMov (X) \cap \Eff (X)$ has a rational polyhedral cone\footnote{It is a closed convex cone in $N^1(X)_\bR$ spanned by finitely many equivalence classes of Cartier divisors on $X$.} $\Pi$ as a fundamental domain, in the sense that
\begin{equation*}\label{coneconjM_f}
    \cMov (X) \cap \Eff (X) = \bigcup_{g \in \Bir (X)} g^{\ast} \Pi \tag{$\star$}
\end{equation*}
and the interiors of $\Pi$ and $g^{\ast} \Pi$ are disjoint unless $g^{\ast} = \id$.
\end{conjecture}

In this article, we restrict our attention to the case that $X$ is a Calabi--Yau threefold of Picard number $\rho (X) = 2$.

By the work \cite{Og-CY2} and \cite{LP13}, if $\Bir (X)$ is infinite, or if one of the boundary rays of $\cMov (X)$ is  rational,
% (which implies that $\cMov (X)$ is rational),
then the movable cone conjecture holds on $X$, cf.~\cite[Proposition 4.1 and Theorem 4.5]{LP13}. The hypothesis that $\rho (X) = 2$ is essentially used.

When $\Bir (X)$ is finite, Conjecture \ref{coneconjM} implies that effective movable cone \eqref{coneconjM_f} is closed and therefore it equals $\cMov (X)$. Moreover, if the answer to Question \ref{quest_introM} \eqref{quest_introM1} below is positive, then Conjecture \ref{coneconjM} holds for Calabi--Yau manifolds of Picard number two as discussed before.

\begin{question} \label{quest_introM}
Let $X$ be a Calabi--Yau manifold. Assume that $\Bir (X)$ is finite.
\begin{enumerate}[(1)]
    \item\label{quest_introM1}  Is $\cMov (X)$ always a rational polyhedral cone?

    \item\label{quest_introM2} Is the number of minimal models of $X$ finite up to isomorphism?
\end{enumerate}
\end{question}

A more detailed discussion of movable cone conjecture and \eqref{quest_introM2} in Question \ref{quest_introM} could be found in \cite[Theorem 2.14]{CL14}.

Conjecture \ref{coneconjM} has been verified for several special cases, see \cite{Borcea91, Kawamata97, Fryers01, Og-CY2, LP13, CO2015, BorisovNuer2016} and references therein, but the full cone conjecture remains open.

Our main result here is to construct a class of smooth complete intersection Calabi--Yau (CICY) threefolds and compute explicitly their birational models. We verify that Conjecture \ref{coneconjM} holds for these Calabi--Yau threefolds, where most of them have finite birational automorphism groups. One of the main interesting examples is the following, see Theorem \ref{mainK3} for the notations.

\begin{example} Consider on $\bP^4$ the vector bundle $\sF=\cO(2)^{2} \oplus \cO(1)$. Then $\bP(\sF)$ is Fano with $\cO_{\bP(\sF)}(-K_{\bP(\sF)})\cong\cO_{\bP(\sF)}(3).$ A complete intersection $X_\sF$ in $\bP(\sF)$ defined by three general sections in $H^0(\cO_{\bP(\sF)}(1))$ is a smooth Calabi--Yau threefold. It has only two flops, denoted by  $X_{\sF}^+$ and $X_{\sE}$, where $X_{\sF}^+$ possesses a K3 fibration and $X_{\sE}$ possesses an elliptic fibration. Let $H$ (resp.~$L$) denote the restriction of the pullback of the hyperplane class on $\bP^4$ (resp.~the corresponding divisor class of $\cO_{\bP (\sF)} (1)$) to $X_{\sF}$. The slice of the movable cone $\cMov(X_\sF)$ is a subdivision of a closed interval, which comes from the chamber structure of the cone:
%Let us denote by $H$ the restriction to $X_{\sF}$ of the pullback of the hyperplane class on $\bP^4$, and $L$ by the restriction to $X_{\sF}$ of the corresponding divisor of $\cO_{\bP (\sF)} (1)$. 

\begin{center}
\begin{tikzpicture}
\draw(0,0)--(9,0);

\foreach \x in {0,3,6,9}
  \draw[very thick, teal] (\x,5pt)--(\x,-5pt);

\foreach \y/\ytext in {0/$5 H - L$,3/$H$,6/$L - H$,9.5/$L - 2 H$}
  \draw (\y,0) node[above=1ex] {\ytext};

\foreach \z/\ztext in {-0.5/$\bP^2$,1.5/$X_{\sE}$,4.5/$X_{\sF}$,7.8/$X_{\sF}^+$,9.5/$\bP^1$}
  \draw (\z,0) node[below] {\ztext};
\end{tikzpicture}
\end{center}
We remark that $X_{\sE}$ is a smooth CICY threefold of bidegrees $(2, 1), (2, 1)$ and $(1, 1)$ in $\bP^4 \times \bP^2$.
\end{example}

We say that a Fano manifold $P$ is $\bP^{n}$\emph{-ruled over} $M$ if $P = \bP (\sF)$ for some vector bundle $\sF$ of rank $n + 1$ over a projective manifold $M$. Such $\sF$ is also called a \emph{Fano bundle}, see Definition \ref{ruledFano}.

We will consider smooth Calabi--Yau threefolds contained in certain $\bP^{n}$-ruled Fano manifold with Picard number $2$. The following theorem is the prototype of the result we aim to establish.

\begin{theorem} Let $P = \bP (\sF)$ be a $\bP^n$-ruled Fano manifold over $\bP^4$ of Fano index $n + 1 \geqslant 2$. We assume that $P$ is normalized, that is, $\sF$ is ample and $\cO (K_P) \cong \cO_P (n + 1)$. Then a complete intersection
$$
X_{\sF} = Z_1 \cap \cdots \cap Z_{n + 1} \subseteq P
$$
of general hypersurfaces $Z_i \in |\cO_P (1)|$ is a smooth Calabi--Yau threefold. Moreover, all the birational models of $X_\sF$ are constructed and the movable cone conjecture holds on $X_{\sF}$.
\end{theorem}

Our proof depends on the classification of Fano bundles, where for most cases the vector bundle $\sF$ splits, see Theorem \ref{r>2} and Theorem \ref{r=2}.
As a generalization, we establish the following theorem, see Theorems  \ref{mainthmrk2}, \ref{mainthma11}, \ref{mainK3}, \ref{mainBordiga}, \ref{mainP4quin} and \ref{mainGr41} for the details.

\begin{theorem} \label{mainthmM_intro}
Let $M$ be a smooth Fano fourfold with $\Pic(M)=\bZ[\cO_M(1)]$ and Fano index $r_M \geqslant 2$, i.e., $\cO_M(-K_M)\cong\cO_M(r_M) $. Let $\displaystyle{\sF = \oplus_{i = 1}^{n + 1} \cO_M (a_i)}$ and $\displaystyle{\sE = \oplus_{i = 1}^{n + 1} \cO_M (b_i)}$ be direct sums of line bundles, where $(a_i)_i$ and $(b_i)_i$ are sequences of nonnegative integers. Suppose that the Calabi--Yau condition holds for the pair $(\sF, \sE)$, i.e.,
$\sum_{i = 1}^{n + 1} (a_i + b_i) = r_M$, and $M$ is not del Pezzo of degree $1$. Then for a general section
\begin{equation}\label{mainthmM_intro1}
    s \in H^0 (\bP (\sF), \sE \boxtimes \cO_{\bP (\sF)} (1)),
\end{equation}
the zero scheme $X_{\sF} = Z (s)$ is a smooth Calabi--Yau threefold of Picard number $2$. Moreover, all the birational models of $X_\sF$ are constructed and the movable cone conjecture holds for $X_{\sF}$.
\end{theorem}

Notice that replacing $\sF$ with its tensor product with $\cO_M (c)$ has the effect of replacing the line bundle $\cO_{\bP (\sF)} (1)$ by $\cO_M (-c) \boxtimes \cO_{\bP (\sF)} (1)$, but does not change $\bP (\sF)$. In particular, it does not affect the Calabi--Yau condition. We will give the complete list of such pairs $(\sF, \sE)$ up to a twist by $\cO_M (c)$, see Proposition \ref{cicylist}.

Our theorem unifies several known examples in the literature and provides evidence to the movable cone  conjecture. In general, it is very hard to find explicit birational models of a given Calabi--Yau threefold. In our case, this is overcome by two key ingredients. First of all, any CICY threefold we consider is naturally equipped with a small contraction together with its flop, see Section \ref{detcontrsubsec}. The flop is over a determinantal hypersurface $D$ in a smooth Fano fourfold $M$. For a general $s$ in \eqref{mainthmM_intro1}, $D$ is a nodal hypersurface, that is, it has only ordinary double points (ODPs). This part is established in \cite{SSW18}.

Second of all, by using the geometric construction of Eagon--Northcott complexes, see Proposition \ref{negdiag} and Remark \ref{negdiagRmk}, it gives rise to a special surface $S_{\sF}$ in our Calabi--Yau $X_{\sF}$. This enables us to find all birational models and hence the full movable cone with its chamber structure, except in two cases. The remaining cases are when $(M, \sF) = (\bP^4, \cO(1)^5)$ or $(\Gr (2, 4), \cO (1)^4)$. We will treat these cases in Section \ref{movfansecII}. In contrast to Section \ref{movfansecI}, the birational automorphism groups of smooth Calabi--Yau threefolds associated with these two exceptional cases have infinite order.

Finally, we make two remarks. Firstly, the non-split case $(\Gr(2,4),\cS(2)\oplus\cO(1))$ and the del Pezzo of degree $1$ (cf.~Proposition \ref{fanofdbpf}) involve more complicated computations and will be discussed in a forthcoming paper. Secondly, the construction in this paper applies to higher dimensional $\bP^n$-ruled Fano manifolds. It is interesting to know what kind of higher dimensional Calabi--Yau manifolds appear and investigate their birational geometry.

%cases in Theorem \ref{r>2}

The paper is organized as follows. In Section \ref{prelsec}, we have included some basic facts and results about Fano bundles and extremal contractions from smooth Calabi--Yau threefolds. The list of Fano bundles we consider is given in Proposition \ref{cicylist}. In Section \ref{degenlocisec}, we recall some general results about degeneracy loci, including Bertini-type and Lefschetz-type Theorems. We also provide the construction and results for determinantal contractions. Section \ref{flopsec} contains the geometric construction of Eagon--Northcott complexes. Section \ref{movfansecI} and \ref{movfansecII} are devoted to the proof of the main results, Theorem \ref{mainthmM_intro}. To streamline our exposition, we recall the definition of Chern classes of virtual quotient bundles and collect the computation of the Hodge numbers of our Calabi--Yau threefolds only in the Appendices \ref{chernclsubsec} and \ref{hodgenumsubsec}.

\begin{notation} Throughout this paper we work over the complex field $\bC$. All varieties are reduced and irreducible, and we do not distinguish a vector bundle and its associated locally free sheaf. For a vector bundle  $\sF$, we write $\bP(\sF)=\Proj(\Sym^\bullet\sF)$ for the projective bundle of $1$-dimensional quotients of $\sF$ as in \cite{Har:AG} and $\cO_{\sF}(1) \coloneqq \cO_{\bP(\sF)} (1)$ for the tautological line bundle. For a morphism $\sigma \colon \sE^{\vee} \to \sF$ of vector bundles, we say that a property holds for a \emph{general} $\sigma$ if it holds for each $\sigma$ in a Zariski open subset of $H^0 (\sE \otimes \sF)$. The Grassmannian $\Gr(k,n)$ stands for the variety of $k$-dimensional subspaces in a fixed $n$-dimensional vector space, and $\cS$ and $\cQ$ are the universal sub- and quotient bundles of the Grassmannian. For a Fano manifold $M$, $\cO_M(1)$ is the line bundle corresponding to a fundamental divisor of $M$. If $\cO(a)$ is a line bundle, then $\cO(a)^t$ stands for $\cO(a)^{\oplus t}$. The self-intersection cycle of a Cartier divisor $D$ is also denoted by $D^t.$ There should be no confusion from the context for the use of these two similar notations. 
\end{notation}

%%%%%%%%%%%%%%%%%%%%
%%%%%%%%%%%%%%%%%%%%

\section{Preliminaries}\label{prelsec}
We prepare some preliminary results on Fano bundles and contractions of Calabi--Yau threefolds to be used in later sections.

\subsection{Fano Manifolds}\label{Fanomfdsubsec}

A smooth projective variety $M$ is called \emph{Fano} if its anticanonical divisor $- K_M$ is ample. It is known that the Picard group of a Fano variety is a finitely generated torsion-free $\bZ$-module. Therefore the greatest integer $r_M$ which divides $\cO (- K_M)$ in $\Pic (M)$ is called the \emph{index} of $M$, i.e., $- K_M \sim r_M H_M$ for some $\cO (H_M) \in \Pic (M)$. The corresponding divisor $H_M$ defined up to the linear equivalence is called a \emph{fundamental divisor} of $M$. We denote by $\cO_M(1)\cong\cO_M(H_M)$ the corresponding invertible sheaf.

It is well-known also that the index of $M$ is at most $\dim M + 1$. Furthermore, $r_M  = \dim M + 1$ if and only if $M \cong \bP^n$, and $r_M  = \dim M$ if and only if $M \cong Q^n \subseteq \bP^{n + 1}$ is a smooth quadric \cite{Kobayashi73}. Note that every $4$-dimensional smooth quadric $Q^4$ is isomorphic to the Grassmannian $\Gr (2, 4)$. %cf.~\cite[Example 1.5]{Ot:sp}.

A Fano variety $M$ is \emph{del Pezzo} or \emph{Mukai} if $r_M = \dim M - 1$ or $\dim M - 2$ respectively. For a modern survey on the classification of such varieties, we refer the reader to \cite{bookAGV} and references therein, see also Section \ref{hodgenumsubsec}.

The following lemma will be used in the proof of Theorem \ref{mainthmrk2}.

\begin{lemma}\label{c2Fano}
Let $M$ be a smooth Fano fourfold of index $r_M$. Let $H_M$ be a fundamental divisor on $M$ and $d_M = H_M^4$ be the degree of $M$. Then
$$
\int_M c_2 (T_M) \cdot H_M^2 =
\begin{cases}
2 d_M + 12 & \mbox{if $M$ is del Pezzo}, \\
d_M +24 & \mbox{if $M$ is Mukai}.
\end{cases}
$$
\end{lemma}

\begin{proof}
To shorten notation, we let $r = r_M$, $d = d_M$, and $H = H_M$. By Kodaira vanishing, Riemann--Roch and $- K_M \sim r H$, we find that
$$
h^0 (\cO(H)) = \frac{(r + 1)^2 }{24} H^4 + \frac{r +1}{24} c_2 (T_M) \cdot H^2  + 1.
$$
Now from standard arguments using the Riemann--Roch, Serre duality, and Kodaira vanishing \cite[Corollary 2.1.14]{bookAGV}, we see that
$$
h^0 (\cO(H)) =
\begin{cases}
\frac{1}{2} d (r - 1) + 3 & \mbox{if } r > 2, \\
\frac{1}{2} d + 4 & \mbox{if } r = 2.
\end{cases}
$$
The lemma follows by comparing the two expressions of $h^0 (\cO (H))$.
\end{proof}

%\begin{remark}
%For the case $r_M = 4$ and $5$, we see that $\int_M c_2 (T_M) \cdot H_M^2 = 14$ and $10$ respectively, as we already knew.
%\end{remark}

\subsection{Fano Bundles}\label{Fanobdsubsec}

\begin{definition} \label{ruledFano}
A vector bundle $\sF$ of rank $r \geqslant 2$ on a projective manifold $M$ is called a Fano bundle if the projective bundle $\bP(\sF)$ is a Fano manifold. We will call such $\bP (\sF)$ a $\bP^{r - 1}$-ruled Fano manifold.
\end{definition}

On $\bP(\sF)$, we denote its canonical divisor by $K_{\sF}=K_{\bP(\sF)}$ and the natural projection morphism by $p_\sF:\bP(\sF)\ra M$. We say that $\sF$ is \emph{ample} if $\cO_{\sF} (1)$ is an ample line bundle on $\bP (\sF)$. From the relative Euler sequence \cite[Ex.III.8.4]{Har:AG}, we have
\begin{equation} \label{KF}
    \cO (K_{\sF}) \cong p_{\sF}^{\ast} \left(\cO (K_M) \otimes \det \sF \right) \otimes \cO_{\sF} (- r).
\end{equation}

%\begin{lemma}[\cite{NO07}] \label{fanobd}
%Suppose that $\sF$ is a vector bundle of rank $r$ on a projective manifold $M$. Then $\bP (\sF)$ is a $\bP^{r - 1}$-ruled Fano manifold of index $r$ if and only if $\sF \otimes \sL$ is ample with $c_1(\sF \otimes \sL) = c_1(T_M)$ for some $\sL \in \Pic (M)$. Moreover, in this case $M$ is Fano.
%\end{lemma}

\begin{lemma}[\cite{NO07}] \label{fanobd}
Suppose that $\sF$ is a vector bundle of rank $r$ on a projective manifold $M$. Then $\bP (\sF)$ is a $\bP^{r - 1}$-ruled Fano manifold of index $r$ if and only if there exists $\sL \in \Pic (M)$ such that $\sF \otimes \sL$ is ample and $c_1(\sF \otimes \sL) = c_1(T_M)$. Moreover, in this case $M$ is Fano.
\end{lemma}

\begin{proof}%[Sketch of the proof]
Let $\sF' = \sF \otimes \sL$. By $c_1 (\sF') = c_1 (T_M)$ and \eqref{KF},
$$
c_1(\cO (K_\sF)) = p_{\sF'}^*(- c_1(T_M)+c_1(\sF')) + c_1(\cO_{\sF'}(-r)) = c_1(\cO_{\sF'}(-r)).
$$
Hence the first assertion that $\bP (\sF)$ is Fano follows as ampleness is a numerical condition. Since $\cO_{\sF}(1)$ can not be expressed as a multiple of other line bundles, it follows that the Fano index of $\bP(\sF)$ is $r$.

Conversely, by \cite[Proposition 3.3]{NO07}, there is a ample twist $\sF'$ of $\sF$ such that $\cO (K_M) \otimes \det \sF'$ is trivial, and hence $c_1 (\sF') = c_1 (T_M)$. As $M$ is the base of a smooth morphism from a Fano manifold, $M$ is Fano by \cite[Corollary 2.9]{KMM}.
\end{proof}

We now turn to the case of $\bP^{r - 1}$-ruled Fano manifolds $\bP (\sF)$ of index $r$, normalized so that $\sF$ is ample and $c_1(\sF)=c_1(T_M)$. As a generalization of Mori's proof on Hartshorne's conjecture \cite{M:aTB}, the classification of such pairs $(M,\sF)$ has attracted intense attention. For our construction of Calabi--Yau threefolds, we focus on the cases when $\dim M=4$. A complete list has been established in a series of works \cite{Wis90,Wis91,P90,P91,PSW92,Wis93,Occ01,Occ05,NO07, Kan19}. The list is shorter when $r\geqslant 3$, see \cite[Proposition 7.4]{PSW92} and \cite{Occ05}.

\begin{theorem}[\cite{PSW92}, \cite{Occ05}] \label{r>2} Let $M$ be a projective manifold of dimension $d$ and $\sF$ a Fano bundle of rank $r\geqslant 3$ with $c_1(\sF)=c_1(T_M)$. Then $r\leqslant d+1$. When $d=4$, the pair $(M,\sF)$ is exactly one of the following:
    \begin{enumerate}[(i)]
        \item $M = \bP^4$ and $\sF$ is given by
            \begin{enumerate}
                \item $r=5$: $\cO(1)^{\oplus 5}$;
                \item $r=4$: $T_{\bP^4}$ or $\cO(1)^{\oplus 3}\oplus\cO(2)$;
                \item $r=3$: $\cO(1)\oplus\cO(2)^{\oplus 2}$ or $\cO(1)^{\oplus 2}\oplus\cO(3)$.
            \end{enumerate}
        \item $M =\Gr(2,4)\cong Q^4\subseteq\bP^5$ via Pl\"ucker embedding and $\sF$ is given by
            \begin{enumerate}
                    \item $r=4$: $\cO(1)^{\oplus 4}$;
                    \item $r=3$: $\cO(1)^{\oplus 2}\oplus\cO(2)$ or  $\mathbf{E}(2)\oplus\cO (1)$, where $\mathbf{E}$ is a spinor bundle with $c_1=-1$ and $c_2=(1,0)$ or $(0,1)$.
            \end{enumerate}
        \item $M$ is del Pezzo with $\Pic(M)=\bZ[\cO(1)]$ and $\sF = \cO(1)^{\oplus 3}$.
        %\footnote{Here $M=V_d$ with $d = H_M^4$ is one of the following, cf.~\cite[Theorem 3.3.1]{bookAGV}:
               %\begin{enumerate}[$(1)$]
                % \item $V_1$: $X_6\subseteq\bP(1^4, 2, 3)$
                % \item $V_2$: $X_4\subseteq\bP(1^5,2)$
                % \item $V_3$: $X_3\subseteq\bP^5$
                % \item $V_4$: $X_{2,2}\subseteq\bP^6$
                % \item $V_5$: a 2-dimensional linear section of $\Gr(2,5)\subseteq\bP^9$
              % \end{enumerate}}
        \item $M = \bP^2 \times \bP^2$ and $\sF=\cO(1,1)^{\oplus 3}$.
    \end{enumerate}
\end{theorem}

We remark that a spinor bundle $\bfE$ on $Q^4 \cong \Gr (2, 4)$ is either the universal subbundle $\cS$ or the dual of the universal quotient bundle $\cQ^\vee$.

When $r=2$, the Picard number $\rho(M)$ can be bigger than 2 and the list is much longer. For simplicity, we only list the classification of $(M,\sF)$ over a Fano manifold $M$ which appears in Theorem \ref{r>2} or is Mukai with $\rho (M) = 1$, cf.~\cite[Theorem 1.1 (3), (4) and Theorem 1.3]{NO07} and \cite{Kan19}.

\begin{theorem}[\cite{NO07, Kan19}] \label{r=2} Let $\sF$ be a Fano bundle of rank two on $M$ with $c_1(\sF)=c_1(T_M)$ and $\rho (M) = 1$. When $M$ is $\bP^4$, $\Gr(2,4)$, del Pezzo, or  Mukai, the pair $(M,\sF)$ is exactly one of the following:
\begin{enumerate}[(i)]
    \item $M = \bP^4$ and $\sF$ is $\cO(1)\oplus\cO(4)$ or $\cO(2)\oplus\cO(3)$;
    \item $M =\Gr(2,4)$ and $\sF$ is $\cO(1)\oplus\cO(3)$ or $\cO(2)\oplus\cO(2)$;
    \item $M$ is del Pezzo with $\Pic(M)=\bZ[\cO(1)]$ and $\sF=\cO(1)\oplus\cO(2)$;
    \item $M$ is Mukai with $\Pic(M)=\bZ[\cO(1)]$ and $\sF=\cO(1)\oplus\cO(1)$.
\end{enumerate}
When $M = \bP^2 \times \bP^2$, we have that $\sF$ is $\cO(1,2)\oplus\cO(2,1)$ or $\cO(1,1)\oplus\cO(2,2)$:
    \begin{enumerate}[(a)]
        \item $\cO(1,2)\oplus\cO(2,1)$: $\bP(\sF)=\Bl_{P_1\cup P_2}\bP^5$ where $P_i$'s are two non-meeting planes.

        \item $\cO(1,1)\oplus\cO(2,2)$: $\bP(\sF)$ is the blow up of a cone in $\bP^9$ over the Segre embedding $\bP^2\times \bP^2\subseteq\bP^8$ along its vertex.
    \end{enumerate}
\end{theorem}

\begin{lemma}\label{ggvb} Let $\sF$ be a vector bundle on a variety $M$. Fix $\ell \geqslant 1$. Then $\cO_\sF(\ell)$ is globally generated if and only if $\Sym^{\ell}\sF$ is globally generated.
\end{lemma}
\begin{proof} Denote by $p_{\ell}:\bP(\Sym^{\ell}\sF)\ra M$ and $p=p_{1}:\bP(\sF)\ra M$ the natural projections. There is a Segre embedding
$$\iota_{\ell}:\bP(\sF)\hookrightarrow\bP(\Sym^{\ell}\sF)$$
such that $\iota_{\ell}^*\cO_{\bP(\Sym^{\ell}\sF)}(1)=\cO_{\bP(\sF)}(\ell)$ with the universal quotient $q_{\ell}:p^*\Sym^{\ell}\sF\twoheadrightarrow \cO_{\bP(\sF)}(\ell)$ given by symmetrizing $q=q_1:p^*\sF\twoheadrightarrow\cO_{\bP(\sF)}(1)$.

The lemma follows from the following set theoretic identity $$p_{\ell}(\Bs(\cO_{\bP(\sF)}(\ell)))=\Bs(\Sym^{\ell}\sF),$$
which we now prove: If $\ev_{\ell}:H^0(M,\Sym^{\ell}\sF)\otimes\cO_M\ra\Sym^{\ell}\sF$
is surjective at $m\in M$, then we can pull it back by $p_l$, compose it with the universal quotient and use the fact that
$$
H^0(M,\Sym^{\ell}\sF)\cong H^0(\bP(\sF),\cO({\ell})),
$$
to conclude that $m\notin p_{\ell}(\Bs(\cO(\ell)))$. Conversely, if $m\in\Bs(\Sym^{\ell}\sF)$, then the image of $\ev_\ell(m)$ is contained in a hyperplane $H(m)\subseteq \Sym^{\ell}\sF(m)$. By the construction of the universal quotient, the evaluation map
$$H^0(\bP(\sF),\cO(\ell))\otimes\cO_{\bP (\sF)}\ra\cO(\ell),$$
which factors through $q_{\ell}\circ p_{\ell}^*\circ\ev_{\ell}$, is then zero at any point $x\in p_l^{-1}(m)$.
\end{proof}

\begin{proposition} \label{fanofdbpf}
The tautological bundle $\cO_\sF(1)$ associated to any Fano bundle in Theorems  \ref{r>2} and \ref{r=2} is globally generated, except when $M$ is a del Pezzo fourfold\footnote{A complete list of classification of del Pezzo fourfold can be found in Theorem \ref{ClassdP4} in the Appendix B.} of degree $1$, which is a degree $6$ hypersurface $X_6\subseteq\bP(1^4,2,3)$.
\end{proposition}
%\footnote{A smooth del Pezzo fourfold of degree $d$ is $M=V_d$, where $-K_M=3H_M$ and $d=H_M^3$, and is one of the following, cf.~\cite[Theorem 3.3.1]{AGV}: 
%$V_1=X_6\subseteq\bP(1^4, 2, 3)$; $V_2=X_4\subseteq\bP(1^5,2)$; $V_3=X_3\subseteq\bP^5$; $V_4=X_{2,2}\subseteq\bP^6$; $V_5$ is a 2-dimensional linear section of $\Gr(2,5)$ in its Pl\"ucker embedding; $V_6=\bP^2\times\bP^2$.} 

\begin{proof} From Lemma \ref{ggvb}, it is enough to show that each Fano bundle $\sF$ on $M$ is generated. Moreover, we only have to check this on each direct summand.

From the Euler sequence there is a surjection $\cO(1)^{5}\twoheadrightarrow T_{\bP^4}$ and thus
$$
H^0 (\bP(T_{\bP^4}), \cO (1)) \otimes \cO \ra \cO(1)
$$
is the restriction of the corresponding sequence on $\bP(\cO(1)^{\oplus 5}) \cong \bP^4\times\bP^4$. Hence $T_{\bP^4}$ is generated.

It is shown in \cite[Theorem 2.8.(ii)]{Ot:sp} that if $\bfE=\cS$, then $\bfE(1)\cong\cS(1)\cong \cS^\vee$ is a quotient of the universal trivial bundle and hence generated. The same holds in case  $\bfE(1)=\cQ^\vee(1)\cong\cQ$. Hence $\bfE(2)$, as a twist of $\bfE(1)$ by the very ample line bundle $\cO(1)$, is generated.

For del Pezzo varieties in Theorem \ref{r>2} and \ref{r=2}, it follows from Fujita's results that $|\cO(1)|$ is generated except the degree $1$ case, see  ~\cite[Proposition 3.2.4]{bookAGV}. For the degree $1$ case, $M=X_6\subseteq\bP(1^4,2,3)$ is a hypersurface of degree $6$. If $x_0,\dots, x_3, y, z$ with $\deg(x_0)=\cdots=\deg(x_3)=1$, $\deg(y)=2$, and $\deg(z)=3$ are homogeneous coordinates of $\bP(1^4,2,3)$, then the defining equation of $X_6$ is of the form 
$$f_6(x_0,\dots,x_3,y,z)=z^2+zh_3(x_0,\dots,x_3,y)+h_6(x_0,\dots,x_3,y),$$
where $h_3$ and $h_6$ are homogeneous polynomials of degree $3$ and $6$ respectively. The base locus of $|\cO (1)|$ consists of points of the form  $[0:\cdots:0:*:*]$, which is clearly non-empty.

For Mukai varieties, it follows from \cite[Proposition 1]{Mukai89} and \cite{Me:good}.
\end{proof}

\subsection{Calabi--Yau Condition} \label{CYsubsec}

Suppose that $M$ is a smooth Fano fourfold and $n \geqslant 1$. Let $\sF = \oplus_{i = 1}^{n + 1} \cO (a_i)$ and $\sE = \oplus_{i = 1}^{n + 1} \cO (b_i)$, where $(a_i)_i$ and $(b_i)_i$ are nondecreasing sequences of nonnegative integers. Up to a twist by a line bundle, we are going to find all pairs $(\sF, \sE)$ such that $\sF$ is ample and the Calabi--Yau condition 
\begin{equation} \label{cycond}
    c_1 (\sF - \sE^{\vee}) = c_1 (T_M)
\end{equation}
holds, that is, $\sum_{i = 1}^{n + 1} (a_i + b_i) = r_M$ and $a_i > 0$ for all $i$. Here the virtual bundle $\sF - \sE^{\vee}$ is in the Grothendieck group of vector bundles on $M$. The formulas for Chern classes of such bundles are given in Appendix \ref{chernclsubsec}.

%Using these formulas we can then rewrite \eqref{cycond} as 

%That is, $\displaystyle{\sum_{i = 1}^{n + 1} (a_i + b_i) = r_M}$ and $a_i > 0$ for all $i$.

\begin{proposition} \label{cicylist}
Under the above assumptions, the triples $(M, \sF, \sE)$ are the following (up to a twist with a line bundle):
\begin{enumerate}[(i)]
    \item\label{cicylist1} If $\sE$ is a trivial bundle, then we have
    \begin{enumerate}
        \item $M = \bP^4$ and the sequence $(a_i)_i$ is one of the following:
        \begin{equation*}
              (4, 1), (3, 2), (3, 1, 1), (2, 2, 1), (2, 1, 1, 1), (1, 1, 1, 1, 1).
        \end{equation*}

        \item $M = \Gr (2, 4)$ and $(a_i)_i = (3, 1), (2, 2), (2, 1, 1)$, or $(1, 1, 1, 1)$.

        \item\label{cicylist1c} $M$ is a smooth del Pezzo fourfold and $(a_i)_i = (2, 1)$ or $(1, 1, 1)$.

        \item\label{cicylist1d} $M$ is a smooth Mukai fourfold and $(a_i)_i = (1, 1)$.
    \end{enumerate}

    \item\label{cicylist2} If $\sE$ is not trivial, then $\sE = (\cO (1) \oplus \cO^n)$ and $(M, \sF)$ is given by
    \begin{enumerate}
        \item $n = 2$: $(\bP^4, \cO (2) \oplus \cO(1)^2)$

        \item $n = 1$: $(\bP^4, \cO (3) \oplus \cO(1))$ or $(\Gr (2, 4), \cO (2) \oplus \cO (1))$
    \end{enumerate}
\end{enumerate}
\end{proposition}

\begin{proof}
First, we adopt the convention that if all  $a_i$'s are the same, say equal to $a$, then we replace $(a_i, b_i)$ by $(0,b_i + a)$ for all $i$ and interchange $\sF$ and $\sE$, and similarly for $b_i$'s. Note that the replacements preserve the Calabi--Yau condition.

Notice that $2 \leqslant r_M \leqslant 5$. Indeed, it is known that the Fano index $r_M$ of $M$ is less than or equal to $\dim M + 1 = 5$. By assumption, $r_M$ is greater than or equal to $\sum_i a_i \geqslant 2$.

With our convention and $a_i > 0$ for all $i$, the cases when $r_M=2$ or 3 are easy, which correspond to items \eqref{cicylist1} (c) and (d). We only need to consider $r_M = 5$ or $4$. In the latter cases, we have $M = \bP^4$ or $\Gr (2,4)$ and there are six and four nontrivial partitions of $5$ and $4$ respectively. Then the proposition follows from an easy calculation. Remark that on $\bP^4$ we identify the case $(\cO (2) \oplus \cO (1), \cO (2) \oplus \cO )$ with $(\cO (3) \oplus \cO (1), \cO (1) \oplus \cO)$ by tensoring $\cO (- 1)$ and interchanging $\sF$ and $\sE$.
\end{proof}

\begin{remark}
The list \eqref{cicylist1} in Proposition \ref{cicylist} is a special case of Theorem \ref{r>2} and \ref{r=2}.
\end{remark}

\subsection{Primitive Contractions}\label{pricontrsubsec}

%Let $X$ be a smooth Calabi--Yau threefold. Recall from \cite[p.~566]{Wilson92} that a birational morphism $\pi : X \to Y$ is a \emph{primitive contraction} if $Y$ is normal and $\pi$ cannot be factored in the algebraic category. This is equivalent to the condition that the relative Picard number of $\pi$ is $1$.

We first recall some terminology from \cite[p.~566]{Wilson92}.

\begin{definition}
Let $X$ be a smooth Calabi--Yau threefold. We say that a birational morphism $\pi : X \to Y$ is a \emph{primitive contraction} if $Y$ is normal and the relative Picard number of $\pi$ is $1$.
\end{definition}

This is equivalent to the condition that $\pi$ cannot be factored in the algebraic category.

\begin{definition}
We say that a birational morphism is \emph{small} if it contracts only finitely many curves, and a primitive contraction is
\begin{enumerate}[(i)]
    \item of type I if it is small;
    \item of type II if it contracts an irreducible surface down to a point;
    \item of type III if it contracts an irreducible surface down to a curve.
\end{enumerate}
\end{definition}

We conclude this section with the following three simple results, which will be used in  Section \ref{movfansecI}. For the convenience of the reader, we supply proofs here.

\begin{lemma}\label{ODPsmsurf}
Let $\pi : X \to Y$ be a small resolution of a normal threefold $Y$ and $C$ an irreducible exceptional curve. Suppose that $K_X$ is $\pi$-trivial and there is a smooth surface $S$ in $X$ such that $C\subseteq S$ is a $(- 1)$-curve in $S$. Then the normal bundle of $C$ in $X$ is isomorphic to $\cO_{\bP^1} (-1) \oplus \cO_{\bP^1} (-1)$.
\end{lemma}

\begin{proof}
Consider the normal bundle sequence
\begin{equation} \label{normalbdlC}
    0 \to N_{C / S} \to N_{C / X} \to N_{S / X} |_C \to 0.
\end{equation}
Since $C \cong \bP^1$ is a $(- 1)$-curve in $S$, this implies that the normal bundle $N_{C / S}$ is $\cO_{\bP^1} (- 1)$. By \eqref{normalbdlC}, $K_X \cdot C = 0$, and adjunction formula, we get
\[
    \deg N_{S / X} |_C = S\cdot C=K_S\cdot C - K_X\cdot C= - 1.
\]
Then $N_{S / X} |_C$ is also $\cO_{\bP^1} (- 1)$ and thus the exact sequence \eqref{normalbdlC} splits.
\end{proof}

\begin{proposition}\label{divcontr}
Let $X$ be a smooth Calabi--Yau threefold and $D \subseteq X$ an irreducible smooth surface.
\begin{enumerate}[(i)]
    \item\label{divcontr1} If $D$ is a K3 surface, then $| D|$ is a base point free linear pencil that induces a fibration $X \to \bP^1$ whose a general fiber is a K3 surface.

    \item\label{divcontr3} If $D \cong \bP^2$, then there is a  primitive contraction $X \to Y$ which contracts the divisor $D \subseteq X$ to a $\frac{1}{3}(1,1,1)$-point $p \in Y$.
\end{enumerate}
\end{proposition}

\begin{proof}
In case \eqref{divcontr1}, the exact sequence
\[
0\ra\cO_X\ra\cO_X(D)\ra\cO_D(D)\ra0
\]
and $H^1 (\cO_X) = 0$ imply that $|D|$ is base point free and $h^0(\cO_X (D))=2$ as $\cO_D(D)\cong \cO_D (K_D) \cong\cO_D$. It follows that $X\ra\bP(|D|)$ has connected fibers. Notice that by upper semicontinuity \cite[III Theorem 12.8]{Har:AG} we have $h^1 (\cO_F) \leqslant h^1(\cO_D) = 0$ for a general fiber $F$. Therefore $K_F = (K_X + F)|_F = 0$ and $h^1 (\cO_F) = 0$, that is, a general fiber $F$ is a K3 surface\footnote{In fact, since we are work over $\bC$, all smooth fibers are diffeomorphic to the K3 surface $D$ by Ehresmann's theorem, and hence are K3 surfaces (see \cite[VII Corollary 3.5]{FM94}).}.

%by upper semicontinuity of $h^1(\cO_D)$, we see that a smooth fiber $F$ is a K3 surface.

In case \eqref{divcontr3}, we have the following more general fact (cf.~\cite[Lemma 2.5]{Kapustka209I}): Let $D$ be a del Pezzo surface. If there is an ample divisor $A$ on $X$ such that $\lambda A |_{D} \sim - K_D$ for some $\lambda > 0$, then $L:=D + \lambda A$ is obviously nef and big on $X$.  By the base-point-free theorem, some multiple of it gives a divisorial contraction $\varphi_L:X \to Y$ contracting the divisor $D$. When $\rho(D)=1$, any ample divisor on $X$ works and $\varphi_L$ contract $D$ to a point. If $D=\bP^2,$ then the argument of \cite[(3.3.5)]{Mori82} proves that $\widehat{\cO}_{Y, p} \simeq \bC [\![x, y, z ]\!]^G$, where $G \coloneqq \bZ / 3 \bZ$ acts on $\bC [\![x, y, z ]\!]$ via the weight $(1, 1, 1)$. Notice that in our case, $\cO_D (D) \cong \cO_{\bP^2}(- 3)$ because $X$ is Calabi--Yau.
\end{proof}

%\footnote{Suppose that $(A, \fm, \kappa = A/ \fm)$ is a Noetherian local ring, $f \in \fm$ and $A / (f)$ is a regular local ring of dimension $\dim A - 1$. We claim that $A$ is a regular local ring, i.e., $\dim A = \dim_{\kappa} \fm / \fm^2$. It follows from the short exact sequence of $\kappa$-vector spaces
%\[
%    0 \to (\fm^2 + (f)) / \fm^2 \to \fm / \fm^2 \to \fm / (f) \big/ (\fm^2 + (f)) / (f) \to 0,
%\]
%$\dim A \leqslant \dim_{\kappa} \fm / \fm^2$ and the hypothesis of $A / (f)$.}

\begin{lemma} \label{A1flop}
Let $Y$ be a nodal threefold with isolated ODPs, denoted by $\Sigma$. Suppose further that there exists a smooth surface $T \supseteq \Sigma$. Then:
\begin{enumerate}[(i)]
    \item\label{A1flopI} The blow-up $\pi \colon X \coloneqq \Bl_{T} Y \to Y$ is a small resolution, and the restriction $\pi^{-1} (T) \to T$ is the blow-up of $T$ at the smooth points of $T$ located at ODPs of $Y$.
    
    \item\label{A1flopII} Let $\pi^+:X^+\ra Y$ be the Atiyah flop obtained by taking $p \colon W \coloneqq \Bl_{\Sigma} Y \ra X$ and then a blow down $q \colon W\ra X^+$ along its exceptional divisors $E$ in the other direction. 
    \begin{equation*} 
    \begin{tikzcd}
    & W \ar[ld, "p", swap] \ar[rd, "q"] &\\
    X && X^{+}
    \end{tikzcd}
\end{equation*}
    Then the proper transform $T^+$ of $T$ in $X^{+}$ is isomorphic to $T$ via $\pi^+$.
    %If we flop $X$ to the other small resolution $X^{+}$ of Y, then the proper transform of $S$ in $X^{+}$ is isomorphic to $S$.
    
    \item\label{A1flopIII} If any two of $X$, $X^{+}$ and $Y$ are projective, then so is the third.
\end{enumerate}
\end{lemma}

\begin{proof}
%Since the flop $X^+$ can be constructed by taking $p \colon W \ra X$ and then a blow down $q:W\ra X^+$ along its exceptional divisors in the other direction, 

The proof of \eqref{A1flopI} can be found in \cite[Lemma 7.1]{BorisovNuer2016} and \eqref{A1flopIII} in \cite[Lemma 7.3]{Friedman91}. To prove \eqref{A1flopII}, we let $T_W$ be the proper transforms of $T$ in $W$. Note that $q (T_W)$ is the proper transform $T^+$ of $T$. According to that $T \supseteq \Sigma$ and $W = \Bl_{\Sigma} Y$, it follows that $T_W$ is the blow-up of $T$ along $\Sigma$ and thus $p \colon T_W \xrightarrow{\sim} \pi^{- 1} (T)$ is an isomorphism by \eqref{A1flopI}. 

%Let $C$ be the exceptional curve of $T_W \cong \pi^{-1} (T)$. ....with normal bundle $\rho_1^{\ast} \cO_{\bP^1} (- 1) \otimes \rho_2^{\ast} \cO_{\bP^1} (- 1)$. Here we choose the isomorphism so that $p|_E$ and $q|_E$ are the projections $\rho_1$ and $\rho_2$ respectively.

For simplicity, we now assume that $\Sigma = \{p\}$ and then the exceptional divisor $E \cong \bP^1 \times \bP^1$. Let $C$ be the rational curve $p (E)$. By the construction of $X$, the curve $C$ is the $(- 1)$-curve $\pi^{- 1} (p)$ on $\pi^{- 1} (T) \cong T_W$. Abusing notation slightly, we use the same letter $C$ for the curve $p^{- 1} (C)$ in $T_W$. On the other hand, by the construction of $X^+$, the induced morphism $q \colon T_W\ra T^+$ for proper transforms contracts the $(-1)$-curve $C$ on $T_W$. We claim that $T^+$ is normal. It follows that $T^+$ is smooth (cf.~\cite[p.415, Step 5]{Har:AG}) and hence the morphism $\pi^+ \colon T^+ \to T$ of smooth surfaces is an isomorphism because they are (set-theoretically) a bijection. For the claim, observe that $T^+$ is smooth outside the point $q (C)$. Since the threefold $X^+$ is smooth, the surface $T^+$ is an effective Cartier divisor and hence is Cohen--Macaulay \cite[Proposition 8.23]{Har:AG}. Therefore it follows from Serre's criterion for normality that $T^+$ is normal.
\end{proof}

%%%%%%By the construction of $X^+$, the induced morphism $q \colon T_W\ra T^+$ for proper transforms contracts exactly those $(-1)$-curves on $T_W$. In particular, $T_W={\rm Bl}_\Sigma T\ra T^+$ factors abstractly as $T_W\ra T\xra{\mu} T^+$. On the other hand,  the natural induced morphism  $\pi^+|_{T^+}:T^+\ra T$ is birational with a discrete exceptional set. Since $\pi^+|_{T^+}\circ\mu:T\ra T$ is birational and identity on a big open set, $\pi^+|_{T^+}\circ\mu={\rm id}_T$ is the identity. The same argument implies that $\mu\circ \pi^+|_{T^+}={\rm id}_{T^+}$ and hence that $\pi^+|_{T^+}:T^+\xra{\cong} T$.

%By the construction of $X^+$, the induced morphism $T_W\ra T^+$ for proper transforms of $T$ on $W$ and $X^+$ contracts exactly those $(-1)$-curves on $T_W\cong \pi^{- 1}(T)$. 

Note that the above proof only uses the local structure of the Atiyah flop, hence the lemma also applies to a singular surface as long as the ODPs on $Y$ are smooth points on $T$.

%%%%%%%%%%%%%%%%%%%%
%%%%%%%%%%%%%%%%%%%%

\section{Degeneracy Loci}\label{degenlocisec}

\subsection{Bertini-type and Lefschetz-type Theorems}\label{BertLefsubsec}

Let $M$ be a variety of dimension $d$, and let $\sigma : \sE^{\vee} \to \sF$ be a morphism of vector bundles on $M$ of rank $e$ and $f$ respectively. For each $k \leqslant \min \{e, f\}$ one can associate to $\sigma$ its $k$th \emph{degeneracy locus}
$$
D_k (\sigma) = \{ x \in M \mid \rank (\sigma (x)) \leqslant k\},
$$
with the convention $D_{- 1} (\sigma) = \varnothing$. Its ideal is locally generated by $(k + 1)$-minors of a matrix for $\sigma$. Notice that the $0$th degeneracy locus of $\sigma$ is the zero scheme $Z (\sigma)$ of the corresponding section of $\sE \otimes \sF$. The \emph{expected codimension} of $D_k (\sigma)$ in $M$ is $(e - k) (f - k)$, though the degeneracy locus may be empty or have strictly smaller codimension.

The following Bertini-type theorem is well known and relies on generic smoothness.

\begin{theorem}[\cite{Banica91}] \label{Berti}
Let $\sE$ and $\sF$ be vector bundles of ranks $e$ and $f$ on a smooth variety $M$ such that $\sE \otimes \sF$ is globally generated. If $\sigma : \sE^{\vee} \rightarrow \sF$ is a general morphism, then one of the following holds:
\begin{enumerate}[(i)]
  \item\label{berti1} $D_k(\sigma)$ is empty;
  \item\label{berti2} $D_k(\sigma)$ has expected codimension $(e - k)(f -k)$ and the singular locus of $D_k(\sigma)$ is $D_{k - 1}(\sigma)$.
\end{enumerate}
Here "general" means that there is a Zariski open set in the vector space $H^0 (\sE \otimes \sF)$ such that either \eqref{berti1} or \eqref{berti2} holds for all $\sigma$ belonging to the open set.
\end{theorem}

We make the following notion of the generality of morphisms used in \cite{SSW18}.

\begin{definition} \label{generaldef}
For a given integer $r \geqslant 0$, a morphism $\sigma : \sE^{\vee} \to \sF$ is said to be $r$-general if the subset $D_i(\sigma) \setminus D_{i - 1}(\sigma)$ is smooth of (expected) codimension $(e - i)(f - i)$ in the smooth variety $M$ for all $i = 0, 1, \cdots, r$.
\end{definition}

%\begin{remark}\label{r-gen}
%Applying Theorem \ref{Berti} repeatedly, we see that if $\sE \otimes \sF$ is globally generated, then there is a Zariski open set $U$ in $H^0 (\sE \otimes \sF)$ such that all $\sigma$ belonging to $U$ are $r$-general or $D_r (\sigma) = \varnothing$.
%\end{remark}

%Applying Theorem \ref{Berti} repeatedly, we see that if $\sE \otimes \sF$ is ample and globally generated, then there is a Zariski open set $U$ in $H^0 (\sE \otimes \sF)$ such that for all $\sigma \in U$ the subset $D_i(\sigma) \setminus D_{i - 1}(\sigma)$ is smooth of (expected) codimension $(e - i)(f - i)$ in the smooth variety $M$ for all $0 \leqslant i \leqslant r$

%\begin{corollary}
%Let $\sE$ and $\sF$ be as in Theorem \ref{Berti}. For a given integer $r \geqslant 0$, there is a Zariski open set $U$ in $H^0 (\sE \otimes \sF)$ such that for all $\sigma \in U$ either the subset $D_i(\sigma) \setminus D_{i - 1}(\sigma)$ is smooth of (expected) codimension $(e - i)(f - i)$ in the smooth variety $M$ for all $0 \leqslant i \leqslant r$ or else $D_r (\sigma) = \varnothing$.
%\end{corollary}

%\begin{proof}
%Applying Theorem \ref{Berti} with $k = r$, we find that there is a Zariski open set $U_r$ in $H^0 (\sE \otimes \sF)$ such that for all $\sigma \in U$ either $D_r (\sigma) = \varnothing$ or else the condition \eqref{berti2} in Theorem \ref{Berti} holds. 
%\end{proof}

The following is from \cite[Theorem 7.1.1, 7.2.1, Example 7.1.5]{LazPAGII}.

\begin{theorem}[\cite{LazPAGII}] \label{Lefprop}
Let $\sigma : \sE^{\vee} \to \sF$ be a morphism between vector bundles of rank $e$ and $f$ on a projective variety $M$ of dimension $d$, and assume that the bundle $\sE \otimes \sF$ is ample.
\begin{enumerate}[(i)]
    \item\label{Lefprop1} The $k$th degeneracy locus $D_k (\sigma)$ of $\sigma$ is non-empty (resp., connected) if $d \geqslant (e - k) (f - k)$ (resp., $d > (e - k) (f - k)$).
     \item Assume that $M$ is smooth, and let $X = D_0(\sigma)$. Then
        $$
        H^i (M, X ; \bZ) = 0 \text{ for } i \leqslant d -e f,
        $$
       the restriction map $H^i(M , \bZ) \to H^i(X , \bZ)$ is an isomorphism for $i < d - e f$ and injective when $i = d - e f$. In particular, if $X$ is also smooth, then
       \begin{enumerate}
           \item  the restriction maps $H^q (M, \Omega_M^p) \to H^q (X, \Omega_X^p)$ are isomorphisms for $p + q < d - ef$;

           \item the restriction map $\Pic (M) \to \Pic (X)$ on Picard groups is an isomorphism if $2 < d - ef$.
       \end{enumerate}
\end{enumerate}
\end{theorem}

Note that Theorem \ref{Lefprop} does not require $D_k (\sigma)$ to have the expected codimension.

\begin{remark}\label{r-gen}
Suppose $\sE \otimes \sF$ is ample and globally generated. For a given integer $r \geqslant 0$, by applying Theorem \ref{Berti} repeatedly, we find that there is a Zariski open set $U$ in $H^0 (\sE \otimes \sF)$ such that all $\sigma$ belonging to $U$ are $r$-general. Note that if $\dim M \geqslant (e - k)(f- k)$ then $D_k (\sigma) \neq \varnothing$ by Theorem \ref{Lefprop} \eqref{Lefprop1}.
\end{remark}

\subsection{Determinantal Contractions} \label{detcontrsubsec}

Let $\sF$ be a coherent sheaf and $\sE$ a vector bundle on a variety $M$. Recall that $p_{\sF} : \bP (\sF) \to M$ is the projection. For a morphism $\sigma : \sE^{\vee} \to \sF$ of $\cO_M$-modules, we can view the composite of $p_{\sF}^{\ast} \sE^{\vee} \to p_{\sF}^{\ast} \sF$ and the canonical map $p_{\sF}^{\ast} \sF \to \cO_{\sF} (1)$ as a global section $s_{\sigma}$ of the bundle
\begin{equation} \label{Cayleytrick}
    \sHom (p_{\sF}^{\ast} \sE^{\vee}, \cO_{\sF} (1)) \cong p_{\sF}^{\ast} \sE \otimes \cO_{\sF} (1).
\end{equation}
Write $\sC$ for the cokernel sheaf of $\sigma$ and consider the zero locus $Z (s_{\sigma})$.

\begin{lemma} \label{zeroeqproj}
There is an isomorphism $Z (s_{\sigma}) \cong \bP (\sC)$ as subschemes of $\bP (\sF)$. 
\end{lemma}

%In particular, $\bP (\sC)$ is the subscheme of $\bP (\sF)$ defined as the zero locus of $s_{\sigma}$ naturally given by $\sigma$.

\begin{proof}
If we can prove that the functors of points induced by $Z (s_{\sigma})$ and $\bP (\sC)$ are isomorphic, then the lemma follows from Yoneda's lemma.

Recall that the $M$-scheme $\bP (\sC)$ represents the functor that attaches to every $M$-scheme $f : T \to M$ the set of equivalence classes of quotients $\lambda : f^{\ast} \sC \to \sL$ where $\sL$ is a line bundle on $T$. The surjection $\sF \twoheadrightarrow \sC$ induces a closed embedding $\bP(\sC) \hookrightarrow \bP (\sF)$, which sends a $T$-valued point $[\lambda]$ to the class $[\mu]$ of
$$
\mu : f^{\ast} \sF \to f^{\ast} \sC \xrightarrow{\lambda} \sL.
$$
Observe that a quotient $\mu$ of $f^{\ast} \sF$ factors through $f^{\ast} \sC$ if and only if the composition of $\mu$ with $f^{\ast} \sigma : f^{\ast} \sE^{\vee} \to f^{\ast} \sF$ is zero.

Let $g : T \to \bP (\sF)$ be the morphism of $M$-scheme associated to a $T$-valued point $[\mu]$ of $\mu : f^{\ast} \sF \to \sL$, which satisfies $\sL = g^{\ast} \cO_{\sF} (1)$ and $f = p_{\sF} \circ g$. Then the morphism $g : T \to \bP (\sF)$ factors through $Z (s_{\sigma})$ if and only if the morphism of bundles $g^{\ast} p_{\sF}^{\ast} \sE^{\vee} \to g^{\ast} \cO_{\sF} (1)$ induced by $g^{\ast} s_{\sigma}$ is zero, which is equivalent to $\mu : f^{\ast} \sF \to \sL$ factoring through $f^{\ast} \sC$.
\end{proof}

Now assume that $\sF$ is a vector bundle and $\sigma : \sE^{\vee} \to \sF$ is a morphism of bundles of ranks $e \geqslant f$ on $M$. If we write everything in local coordinates, then we see that $\bP (\sC) \cong Z (s_{\sigma})$ maps onto $D_{f - 1}(\sigma)$, cf.~\cite[Example 14.4.10]{fulton98}:
\begin{equation} \label{geomENcpx}
    \begin{tikzcd}
      Z (s_{\sigma}) \ar[r, hook] \ar[d, two heads] & \bP(\sF) \ar[d, "p_{\sF}"] \\
      D_{f - 1} (\sigma) \ar[r, hook] & M\nospaceperiod
    \end{tikzcd}
\end{equation}
We can compute the canonical bundle of $Z(s_{\sigma})$ from \eqref{KF}, \eqref{Cayleytrick}, and the adjunction formula:
\begin{align} 
    \cO (K_{Z (s_{\sigma})}) & \cong \left( \cO (K_{\sF}) \otimes \det (p_{\sF}^{\ast} \sE \otimes \cO_{\sF} (1)) \right)|_{ Z (s_{\sigma})} \label{KZ} \\
    & \cong (\cO_{\sF} (e - f) \otimes p_{\sF}^{\ast} (\cO (K_M) \otimes \det \sF \otimes \det \sE))|_{Z (s_{\sigma})}. \notag
\end{align}
Note that the expected codimension of $Z (s_{\sigma})$ is $e$, and given $x \in D_{f - 1} (\sigma)$ the fiber of $Z (s_{\sigma})$ over $x$ is $\bP(\coker \sigma (x))$.

\begin{lemma} \label{ENcpxsurface}
If $D_{f - 2} (\sigma) = \varnothing$, then $Z(s_{\sigma}) \to D_{f - 1} (\sigma)$ is an isomorphism.
\end{lemma}

\begin{proof}
Since $Z(s_{\sigma})$ is the projectivization $\bP(\sC)$ of the cokernel sheaf $\sC$ of $\sigma$, it suffices to show that the restriction of $\sC$ to $D_{f - 1} (\sigma)$ is a line bundle. By assumption, for any point $x \in D_{f - 1} (\sigma)$ the linear map $\sigma (x)$ has constant rank $f - 1$, so the cokernel $\sC$ is a vector bundle of rank $1$, which completes the proof.
\end{proof}

%\begin{definition} \label{derterdef}
%We denote by $X_{\sF}$ the zero scheme of $s_{\sigma}$ if $e = f = n + 1$. The restriction of $p_{\sF}$ to $X_{\sF}$ is called the \emph{determinantal contraction} of $X_{\sF}$, denoted by $\pi_{\sF} : X_{\sF} \to D_n (\sigma)$.
%\end{definition}

\begin{definition} \label{derterdef}
If $\sE$ and $\sF$ are vector bundles of the same rank $n + 1$, then we denote by $X_{\sF}$ the zero scheme of the section $s_{\sigma}$. The restriction of $p_{\sF}$ to $X_{\sF}$ is called the \emph{determinantal contraction} of $X_{\sF}$, denoted by $\pi_{\sF} : X_{\sF} \to D_n (\sigma)$.
\end{definition}

We state the main results of  \cite[Proposition 3.6 and Theorem 4.4]{SSW18}, which will be used in Sections \ref{movfansecI} and \ref{movfansecII}.

\begin{theorem}[\cite{SSW18}] \label{nodal3F}
With notation as in Definition \ref{derterdef}, we assume that $M$ is a smooth projective fourfold. If $\sigma$ is $n$-general and $X_{\sF}$ is connected, then $D_n (\sigma)$ is a nodal hypersurface and the determinantal contraction $\pi_{\sF}$ is a small resolution.
\end{theorem}

The number of singularities of the nodal determinantal hypersurface is determined by Chern classes of $\sF$ and $\sE$ (cf.~\cite[Remark 3.3]{SSW18}).

\begin{proposition}\label{numbODP}
For an $n$-general $\sigma$, the number of ODPs of $D_{n} (\sigma)$ is
\[
    \int_M c_2 (\sF -\sE^{\vee})^2 - c_1 (\sF -\sE^{\vee}) \cdot c_3 (\sF -\sE^{\vee}).
\]
\end{proposition}

\begin{proof}
Since $\Sing (D_{n} (\sigma)) = D_{n - 1} (\sigma)$ and the (expected) codimension of $D_{n - 1} (\sigma)$ in the smooth fourfold $M$ is $4$, the result follows from Giambelli--Thom--Porteous formula \cite[Theorem 14.4]{fulton98}.
\end{proof}

%\begin{remark} \label{numbODP}
%For an $n$-general $\sigma$, we have $\Sing (D_{n} (\sigma)) = D_{n - 1} (\sigma)$ and the (expected) codimension of $D_{n - 1} (\sigma)$ in the smooth fourfold $M$ is $4$ (cf.~\cite[Remark 3.3]{SSW18}). Therefore the number of ODPs of $D_{n} (\sigma)$ is
%$$
%|\Sing (D_{n} (\sigma))| = \int_M c_2 (\sF -\sE^{\vee})^2 - c_1 (\sF -\sE^{\vee}) \cdot c_3 (\sF -\sE^{\vee})
%$$
%by Giambelli--Thom--Porteous formula \cite[Theorem 14.4]{fulton98}.
%\end{remark}

To study the birational geometry of $X_{\sF}$, we first compute the intersection numbers on it in terms of Chern classes of $\sE$ and $\sF$ (see Section \ref{chernclsubsec}). The following is from \cite[Proposition 4.5]{SSW18}.

\begin{proposition}[\cite{SSW18}] \label{LHinter}
With the assumptions as in Theorem \ref{nodal3F}, let $H_M$ be a Cartier divisor on $M$, $H_\sF = (\pi_{\sF}^*H_M)|_{X_{\sF}}$ and $L_\sF = c_1 (\cO_{\sF} (1)|_{X_\sF})$.
\begin{enumerate}[(i)]
    \item For $k = 0, 1, 2, 3$,
    $$
    \int_{X_{\sF}} H_{\sF}^k \cdot L_{\sF}^{3-k} = \int_M H_M^k\cdot c_{4 - k}(\sE - \sF^{\vee}).
    $$

    \item Under the Calabi--Yau condition $c_1(\sF - \sE^{\vee}) = c_1(T_M)$, we have
    \begin{align*}
        \int_{X_{\sF}} c_2(T_{X_{\sF}})\cdot H_{\sF} &= \int_M c_2(T_M)\cdot c_1(\sE - \sF^{\vee})\cdot H_M \\
        \int_{X_{\sF}} c_2(T_{X_{\sF}})\cdot L_{\sF} &= \int_M c_2(T_M)\cdot c_2(\sE - \sF^{\vee}) - |\Sing (D_n (\sigma))|.
    \end{align*}
\end{enumerate}
\end{proposition}

%%%%%%%%%%%%%%%%%%%%
%%%%%%%%%%%%%%%%%%%%

\section{Birational maps via Matrix Transpositions} \label{flopsec}

From now on, we let $\sE$ and $\sF$ be vector bundles of rank $n + 1$ on a smooth projective fourfold $M$. Assume that $\sF$ is an ample Fano bundle. Then we see that $M$ is Fano (cf.~Lemma \ref{fanobd}). Denote by $H_M$ a fundamental divisor on $M$, $r_M$ the Fano index of $M$, and $d_M = H_M^4$ the degree of $M$. We further assume that $H_M$ is base point free (cf.~Proposition \ref{fanofdbpf}).

%To set up our construction, we assume that $\sE$ and $\sF$ are globally generated and there is an exact sequence of vector bundles
%\begin{equation} \label{Fminus}
%    0 \to \cO_M (a) \to \sF \to \sF_{-} \to 0
%\end{equation}
%with $\cO_M (r_M + a) \otimes \det \sE^{\vee}$ being ample.

To set up our construction, we assume that $\sE$ and $\sF$ are globally generated. Suppose that there are an integer $a > 0$ and a vector bundle $\sF_{-}$ of rank $n$ such that we have an exact sequence of vector bundles
\begin{equation} \label{Fminus}
    0 \to \cO_M (a) \to \sF \to \sF_{-} \to 0
\end{equation}
with $\cO_M (r_M + a) \otimes \det \sE^{\vee}$ being ample. 

Note that $\sE \otimes \sF$ and $\sE \otimes \sF_{-}$ are also ample and globally generated (cf.~\cite[Proposition 6.1.12 (i), Theorem 6.2.12 (iv)]{LazPAGII}). Thus we can apply Bertini-type theorem to these bundles (see Remark \ref{r-gen}). By Theorem \ref{Berti}, we can pick a general morphism $\sigma : \sE^{\vee} \to \sF$, and it induces a general morphism $\sigma_{-} : \sE^{\vee} \to \sF_{-}$. Indeed, we set $\sL = \cO_M (r_M + a) \otimes \det \sE^{\vee}$. By Griffths vanishing theorem (\cite[7.3.2]{LazPAGII}) and the assumption that $\sL$ is ample, we get
\[
\Ext^1 (\sE^{\vee}, \cO_M (a)) \cong H^1 (M, \cO_M (K_M) \otimes \sE \otimes \det \sE \otimes \sL) = 0,
\]
and thus $\Hom (\sE^{\vee}, \sF) \to \Hom (\sE^{\vee}, \sF_{-})$ is surjective.

\begin{remark}
In Sections \ref{movfansecI} and \ref{movfansecII}, vector bundles $\sE$ and $\sF$ are direct sum of line bundles $\cO_M (a_i)$. We will take $a = \max \{a_i\}$, and the above assumptions are easily achieved.
\end{remark}

\begin{notation}
Let $H_\sF = p_{\sF}^*H_M$ and $L_\sF = c_1 (\cO_{\sF} (1))$ on $\bP(\sF)$, and similarly for $H_\sE$ and $L_{\sE}$ on $\bP(\sE)$. Fix a bundle $\sV$ on $M$. By abuse of notation, we write $\sV \boxtimes \cO_{\sF} (1)$ for $(p_{\sF}^{\ast} \sV \otimes \cO_{\sF} (1))$ and use the same notations $L_{\sF}$ and $H_{\sF}$ for their restrictions to $X_{\sF}$, and similarly for bundles $\sE$ and $\sF_{-}$.
\end{notation}
%By abuse of notation, we write $\cO_M (c) \boxtimes \cO_{\sF} (1)$ for $(p_{\sF}^{\ast} \cO_M (c) \otimes \cO_{\sF} (1))$, $L_{\sF}$ for the restriction $L_{\sF} |_{X_{\sF}}$ on $X_{\sF}$, and similarly for $H_{\sF}$, $L_{\sE}$, and $H_{\sE}$.

Recall from Definition \ref{derterdef} that $\sigma$ induces the zero scheme $X_{\sF}$ in $\bP (\sF)$. The zero scheme induced by $\sigma_{-}$ is defined similarly:

\begin{definition}
We denote $S_{\sF} \subseteq \bP (\sF_{-})$ by the zero locus of the global section of $\sE \boxtimes \cO_{\sF_{-}} (1)$ induced by $\sigma_{-}$.
\end{definition}

%We denote by $X_{\sF}$ and $S_{\sF}$ the zero schemes induced by $\sigma$ and $\sigma_{-}$ respectively (see Section \ref{detcontrsubsec}).

Using the existence of \eqref{Fminus}, we can construct the basic diagram \eqref{basic_diag} in the following proposition, which will play an important role in Section \ref{movfansecI}.

\begin{proposition} \label{negdiag}
Under the above assumptions, we have
\begin{enumerate}[(i)]
    \item $X_{\sF}$ is a smooth (irreducible) threefold with
    \begin{equation}\label{KX}
        \cO (K_{X_{\sF}}) \cong (p_{\sF}^{\ast} (\cO (K_M) \otimes \det \sF \otimes \det \sE))|_{X_{\sF}},
    \end{equation}
    the Picard number $\rho(X_\sF)=\rho (\bP (\sF))$ and $H^i (\cO_{X_\sF}) = 0$ for $i = 1, 2$.

    \item $S_{\sF}$ is a smooth (irreducible) surface and belongs to the linear system $|L_{\sF} - a H_{\sF}|$ on $X_{\sF}$.

    \item There is a commutative diagram
    \begin{equation} \label{basic_diag}
        \begin{tikzcd}[row sep=.7em]
          & \bP (\sF_{-}) \ar[rd, symbol=\supseteq]  &\\
          \bP (\sF) \ar[ru, symbol=\supseteq] \ar[rd, symbol=\supseteq] \ar[dd, "p_{\sF}"] & & S_{\sF} \ar[dd, "\rotatebox{270}{\(\sim\)}"]\\
          & X_{\sF}  \ar[ru, symbol=\supseteq] \ar[d, "\pi_{\sF}"] &\\
          M \ar[r, symbol=\supseteq] & D_n(\sigma) \ar[r, symbol=\supseteq] & D_{n - 1}(\sigma_{-}),
        \end{tikzcd}
    \end{equation}
    where the natural contraction $S_{\sF} \to D_{n - 1} (\sigma_{-})$ is an isomorphism.
\end{enumerate}
\end{proposition}

Note that $D_{n - 1} (\sigma_{-})$ contains $D_{n - 1} (\sigma)$, the singular locus of $D_{n} (\sigma)$.

\begin{proof}
Recall that $\sE \otimes \sF$ and $\sE \otimes \sF_{-}$ are ample. By Theorem \ref{Lefprop}, we have isomorphisms $\Pic (\bP (\sF)) \xrightarrow{\sim} \Pic (X_{\sF})$, $H^i (\bP (\sF) , \bZ) \xrightarrow{\sim} H^i (X_{\sF} , \bZ)$ for $i<3$, and $H^j (\bP (\sF_{-}) , \bZ) \xrightarrow{\sim} H^j (S_{\sF} , \bZ)$ for $j < 2$. Hence $X_{\sF}$ and $S_{\sF}$ are connected. According to Theorem \ref{nodal3F}, it follows that $X_{\sF}$ is smooth.
Since $\sF$ and $\sE$ are of the same rank, the formula \eqref{KX} is given by \eqref{KZ}. From the assumption that $\sF$ is Fano, it implies that $H^i (\cO_{X_\sF}) = H^i (\cO_{\bP (\sF)}) = 0$ for $i = 1, 2$.

Let $\sC$ be the cokernel sheaf of $\sigma$, and similarly for $\sC_{-}$. By \eqref{Fminus} and diagram chasing, we get the exact sequence
\begin{equation} \label{Cminus}
    \cO_M (a) \to \sC \to \sC_{-} \to 0.
\end{equation}
The commutative diagram \eqref{basic_diag} follows from the isomorphism $X_{\sF} \cong \bP (\sC)$ and $S_{\sF} \cong \bP (\sC_{-})$ by applying Lemma \ref{zeroeqproj}. According to \eqref{Cminus} and that the tautological line bundle of $\bP (\sC)$ is the restriction $\cO_{\sF} (1) |_{\bP (\sC)}$, it follows that $S_{\sF}$ is defined by a global section of the line bundle $(\cO_M (- a) \boxtimes \cO_{\sF} (1))|_{X_{\sF}}$ and hence $S_{\sF} \in |L_{\sF} - a H_{\sF}|$.

Notice that $D_{n - 2} (\sigma_{-}) = \varnothing$ for a general $\sigma_{-}$ because the expected codimension of $D_{n - 2} (\sigma_{-})$ in the smooth fourfold $M$ is $6$. By Lemma \ref{ENcpxsurface} and Theorem \ref{Berti}, the epimorphism $S_{\sF} \to D_{n - 1} (\sigma_{-})$ is an  isomorphism and $\Sing (D_{n - 1} (\sigma_{-})) = \varnothing$.
\end{proof}

\begin{remark} \label{negdiagRmk}
The commutative diagram \eqref{basic_diag} is the geometric picture that arises if we compare the Eagon--Northcott complex induced by $\sigma$ with that induced by $\sigma_{-}$. See p.321 and $(\mathrm{EN}_0)$ in \cite[Appendix B.2]{LazPAGI}.
\end{remark}

There is the other determinantal contraction $\pi_{\sE} : X_{\sE} \to D_{n} (\sigma^{\vee})$ via the dual morphism $\sigma^{\vee} : \sF^{\vee} \to \sE$, and observe that $D_n (\sigma) = D_n (\sigma^{\vee})$. We can certainly assume that $D_n (\sigma)$ is singular, since otherwise $\pi_{\sF}$ and $\pi_{\sE}$ are isomorphisms. Therefore the determinantal contraction $\pi_{\sE}$ gives rise to a diagram
\begin{equation} \label{flopdiag}
    \begin{tikzcd}
    X_\sF \ar[dr,swap,"\pi_{\sF}"] \ar[rr,dashed,"\chi"] & & X_{\sE} \ar[dl,"\pi_{\sE}"] \\
    & D_n(\sigma) &  .\\
    \end{tikzcd}
\end{equation}

\begin{proposition}
Let $\chi = \pi_{\sE}^{- 1} \circ \pi_{\sF}$. Then the rational map $\chi$ is not an isomorphism.
\end{proposition}

\begin{proof}
Let $T$ denote the smooth surface $D_{n - 1} (\sigma_{-}) = D_{n - 1} (\sigma_{-}^{\vee})$, and let $\sI_T$ be its ideal sheaf in $M$. Apply the Eagon--Northcott complex (see ($\mathrm{EN}_0$) in \cite[p.322]{LazPAGI}) to the morphism $\sigma_{-}^{\vee}$, we get
\begin{equation*}
    0 \to \sF_{-}^{\vee} \xrightarrow{\sigma_{-}^{\vee}} \sE \to \sI_T \otimes \det \sF_{-}^{\vee} \otimes \det \sE \to 0.
\end{equation*}
Since $M$ and $T$ are smooth, the Rees algebra $\bigoplus_{k \geqslant 0} \sI_T^k$ is isomorphic to the symmetric algebra $\Sym^\bullet \sI_T$. Thus the projectivization of the cokernel of $\sigma_{-}^{\vee}$ is isomorphic to  $\Bl_{T} M \cong \bP (\Sym^\bullet \sI_T)$. 

On the other hand, the projectivization of the cokernel of $\sigma^{\vee}$ is contained in that of $\sigma_{-}^{\vee}$ and is isomorphic to $X_{\sE}$ by Lemma \ref{zeroeqproj}. Since $D_n (\sigma^{\vee})$ is a nodal hypersurface and $T$ contains the singular locus of $D_n (\sigma^{\vee})$, we find that $X_{\sE} \cong \Bl_T D_{n} (\sigma^{\vee})$ and thus $\pi_{\sE}^{- 1} (T)$ is isomorphic to the blow-up of $T$ at its smooth points located at ODPs of $D_n (\sigma^{\vee})$ by Lemma \ref{A1flop}. Then proposition follows from that $\chi_{\ast}S_{\sF} = \pi_{\sE}^{- 1} (T) \to T$ is not an isomorphism while $S_{\sF}$ is isomorphic to $T$ by Proposition \ref{negdiag}
\end{proof}

\begin{remark}
The morphism $\sigma$ is defined locally by a matrix of elements in a coordinate ring of an affine open set. The transpose of the matrix is then the corresponding matrix of $\sigma^{\vee}$. Hence the birational map $\chi$ is locally induced by the matrix transposition and $D_n (\sigma) = D_n (\sigma^{\vee})$.
\end{remark}

The remainder of this section will be devoted to compute $\chi_{\ast} L_{\sF}$ under certain assumptions, which is extremely useful in Section \ref{movfansecI}. To simplify the notations, we let $H \coloneqq H_\sF$ and $L \coloneqq L_\sF$. We recall that $\chi$ is an isomorphism in codimension one, and clearly the proper transform $\chi_{\ast} H$ is $H_{\sE}$.

\begin{lemma}\label{ab}
Under the assumptions as in Proposition \ref{negdiag}, if we write $\chi_*L = \alpha L_\sE + \beta H_\sE$ in $\Pic(X_{\sE})_\bQ$, then $\alpha \beta < 0$ and
$$
\begin{cases}
L\cdot H^2 = \alpha L_\sE\cdot H_\sE^2 + \beta H_\sE^3 \\
L^2\cdot H = \alpha^2 L_\sE^2\cdot H_\sE + 2 \alpha \beta L_\sE\cdot H_\sE^2 + \beta^2 H_\sE^3
\end{cases}.
$$
\end{lemma}

\begin{proof}
By assumption, the bundle $\sF$ is globally generated, and so is $L$ by Lemma \ref{ggvb}. Since $H$ and $L$ are base point free on $X_\sF$, we may assume $H^2$ and $L^2$ are represented by 1-cycles avoiding the indeterminacy loci of $\chi$.  By the geometric interpretation of intersection numbers, the lemma follows from
$$
L\cdot H^2=\chi_*L\cdot H_\sE^2\ {\rm and}\ L^2\cdot H=(\chi_*L)^2\cdot H_\sE.
$$

Recall that $L$ is ample. Observe that $\alpha \neq 0$ or otherwise $\chi_*L\equiv_\bQ bH_\sE$ and hence $L\sim_\bQ bH$ can not be ample. On the other hand, $\chi_*L$ cannot be ample or otherwise $\chi$ is an isomorphism by  \cite[Lemma 1.5]{Kawamata97}. Hence $\beta \neq 0$ and the only possibility is $\alpha \beta < 0$, as $\chi_{\ast} L$ is big but not ample.
\end{proof}

\begin{proposition} \label{deterwall}
Under the assumptions as in Lemma \ref{ab}, we assume furthermore that $\sE$ and $\sF$ are direct sums of line bundles $\cO_M (a_i)$ satisfying the Calabi--Yau condition $c_1 (\sF - \sE^{\vee}) = c_1 (T_M)$. If
\begin{equation} \label{deterwallcond}
    \int_M c_2 (\sF - \sE^{\vee})\cdot H_M^2 > 0 > \int_M \left(c_1 (\sE - \sF^{\vee})^2 - 2 c_2 (\sE - \sF^{\vee}) \right)\cdot H_M^2,
\end{equation}
then
$$
\chi_*L = - L_\sE + r_M H_\sE,
$$
where $r_M$ is the Fano index of $M$.
\end{proposition}

\begin{proof}
We begin by proving the equality
\begin{equation} \label{betaeq}
    H_{\sE}^3 \beta^2 - 2 (L \cdot H^2) \beta + L^2\cdot H - L_{\sE}^2\cdot H_{\sE} = 0.
\end{equation}
To deduce \eqref{betaeq} from Lemma \ref{ab}, we write $\alpha = (L_{\sE}\cdot H_{\sE}^2)^{- 1} (L\cdot H^2 - H_{\sE}^3 \beta)$. Substituting this into the equation of $L^2\cdot H$ in Lemma \ref{ab} and denoting the constant $(L_{\sE}^2\cdot H_{\sE}) (L_{\sE}\cdot H_{\sE}^2)^{- 2}$ by $C$, we get
\begin{equation} \label{betaeq1}
    \left[H_{\sE}^3 - C (H_{\sE}^3)^2\right] \beta^2 - 2 \left[L\cdot H^2 - C (H_{\sE}^3) (L\cdot H^2) \right] \beta + L^2\cdot H - C (L\cdot H^2)^2 = 0.
\end{equation}

To apply Proposition \ref{LHinter}, we need some recurrence relations of Chern classes of virtual quotient bundles. Write
\[
    c_1 \coloneqq c_1 (\sF - \sE^{\vee}) = c_1 (\sE - \sF^{\vee}).
\]
By the total Chern classes $c (\sF - \sE^{\vee}) \cdot c(\sE^{\vee} - \sF) = 1$ and $c_k (\sF - \sE^{\vee}) = (- 1)^k c_k (\sF^{\vee} - \sE)$, we have the following recurrence relations
\begin{align}
    c_2 (\sF - \sE^{\vee}) + c_2 (\sE - \sF^{\vee}) &= c_1^2, \label{reuchern1} \\
    c_3 (\sF - \sE^{\vee}) - c_3 (\sE - \sF^{\vee}) &= c_1 \cdot ( c_2 (\sF - \sE^{\vee}) - c_2 (\sE - \sF^{\vee})). \label{reuchern2}
\end{align}
From the assumption that $\sE$ and $\sF$ are direct sums of line bundles $\cO_M (a_i)$, the Chern class $c_{4 - k} (\sE - \sF^{\vee})$ is a multiple of the class $H_M^{4 - k}$. For example, the class $c_{1}$ is the multiple $r_M H_M$ by the Calabi--Yau condition. We denote the constant by $\fc_{4 - k} (\sE - \sF^{\vee})$, and similarly for $\fc_{4 - k} (\sF - \sE^{\vee})$.

Multiplying \eqref{betaeq} by $(L_{\sE}\cdot H_{\sE}^2)^{2}$ gives
\begin{align}
      \left[(L_{\sE}\cdot H_{\sE}^2)^2 - (L_{\sE}^2\cdot H_{\sE}) (H_{\sE}^3)\right] & \left[H_{\sE}^3 \beta^2 - 2 (L\cdot H^2) \beta\right] \notag \\
      &+ (L_{\sE} \cdot H_{\sE}^2)^2 (L^2 \cdot H) - (L^2_{\sE} \cdot H_{\sE}) (L \cdot H^2)^2 = 0. \label{betaeq2}
\end{align}
By Proposition \ref{numbODP} and \ref{LHinter}, we see that
\begin{align*}
    (L_{\sE}\cdot H_{\sE}^2)^2 - (L_{\sE}^2\cdot H_{\sE}) (H_{\sE}^3) &= (\fc_2 (\sF - \sE^{\vee})^2 - \fc_3 (\sF - \sE^{\vee}) \fc_1) d_M^2
\end{align*}
and it equals $|D_{n - 1} (\sigma)| d_M$. Also, we rewrite the term with no $\beta$ in \eqref{betaeq2} as 
\begin{equation*}
\begin{aligned}
    & (L_{\sE} \cdot H_{\sE}^2)^2 (L^2 \cdot H) - (L^2_{\sE} \cdot H_{\sE}) (L \cdot H^2)^2 && \\
    &= \left( \fc_2 (\sF - \sE^{\vee})^2 \fc_3 (\sE - \sF^{\vee}) -\fc_3(\sF - \sE^{\vee}) \fc_2 (\sE - \sF^{\vee})^2 \right) d_M^3 && \\
    &= \big\{ \fc_2 (\sF - \sE^{\vee})^2 [\fc_3 (\sE - \sF^{\vee}) - \fc_3 (\sF - \sE^{\vee})] && \\
    &\phantom{==} + \fc_3 (\sF - \sE^{\vee}) [\fc_2 (\sF - \sE^{\vee})^2 - \fc_2 (\sE - \sF^{\vee})^2] \big\} d_M^3 && \\
    &= \{ [\fc_3 (\sE - \sF^{\vee}) - \fc_3 (\sF - \sE^{\vee})] [\fc_2 (\sF - \sE^{\vee})^2 - \fc_3 (\sF - \sE^{\vee}) \fc_1] \} d_M^3&& \text{by \eqref{reuchern1}, \eqref{reuchern2}} \\
    &= (L^2 \cdot H - L_{\sE}^2 \cdot H_{\sE}) |D_{n - 1} (\sigma)| d_M
\end{aligned}
\end{equation*}
Therefore equation \eqref{betaeq2} becomes 
\[
    (|D_{n - 1} (\sigma)| d_M) (H_{\sE}^3 \beta^2 - 2 (L \cdot H^2) \beta + L^2\cdot H - L_{\sE}^2\cdot H_{\sE}) = 0,
\]
and the equation \eqref{betaeq} follows from $|D_{n - 1} (\sigma)| d_M \neq0$.

We are going to compute the discriminant of the quadratic equation \eqref{betaeq}. From Proposition \ref{LHinter}, we find that
\begin{equation*}
\begin{aligned}
    & 4 (L\cdot H^2)^2 - 4 (H_{\sE}^3) (L^2\cdot H - L_{\sE}^2\cdot H_{\sE}) && \\
    &= 4 d_M^2 \left(\fc_2 (\sE -\sF^{\vee})^2 - \fc_1 (\fc_{3} (\sE - \sF^{\vee}) - \fc_{3} (\sF - \sE^{\vee}))\right) && \\
    &= 4 d_M^2 \left(\fc_2 (\sE -\sF^{\vee})^2 + \fc_1^2 (\fc_2 (\sF - \sE^{\vee}) - \fc_2 (\sE - \sF^{\vee}))\right) && \text{by \eqref{reuchern2}}\\
    &= 4 d_M^2 \left(\fc_2 (\sE -\sF^{\vee})^2 -2 \fc_1^2 \fc_2 (\sE - \sF^{\vee}) + \fc_1^4 \right) && \text{by \eqref{reuchern1}} \\
    &= \left(2 d_M \left(\fc_1^2 - \fc_2 (\sE -\sF^{\vee})\right)\right)^2  && \\
    &= \left(2 d_M \fc_2 (\sF -\sE^{\vee})\right)^2 && \text{by \eqref{reuchern1}},
\end{aligned}
\end{equation*}
and $H_{\sE}^3 = \fc_1 d_M = r_M d_M$ by the Calabi--Yau condition.

By the above equalities, assumption \eqref{deterwallcond}, and quadratic formula, the solutions of quadratic equation \eqref{betaeq} are $\beta_{+} = (r_M)^{- 1} (2 \fc_2 (\sE -\sF^{\vee}) - r_M^2)$ and $\beta_{-} = r_M$ where $\beta_{\pm}$ are positive numbers. According to
\begin{equation*}
\begin{aligned}
    L\cdot H^2 - H_{\sE}^3 \beta_{\pm} &= \fc_2 (\sE - \sF^{\vee}) d_M - r_M d_M \beta_{\pm} && \\
    &= \pm (\fc_1^2 - \fc_2 (\sE - \sF^{\vee})) d_M && \\
    &= \pm \fc_2 (\sF - \sE^{\vee}) d_M = \pm L_{\sE}\cdot  H_{\sE}^2, &&
\end{aligned}
\end{equation*}
it follows that
$$
\alpha_{\pm} = (L_{\sE}\cdot H_{\sE}^2)^{- 1} (L\cdot H^2 - H_{\sE}^3 \beta_{\pm}) = \pm 1.
$$
Hence the only possibility is $(\alpha_{-}, \beta_{-}) = (-1, r_M)$, by $\alpha_{+} \beta_{+} > 0$ and Lemma \ref{ab}.
\end{proof}

%%%%%%%%%%%%%%%%%%%%
%%%%%%%%%%%%%%%%%%%%

\section{Birational Models and Movable Cones I} \label{movfansecI}

Throughout this section, we will use the same notation as in Section \ref{flopsec}. We assume that $M$ is a smooth Fano fourfold with $\rho (M) = 1$ and is not a del Pezzo fourfold of degree $1$. Hence a fundamental divisor $\cO (1)$ of $M$ is globally generated, cf.~Proposition \ref{fanofdbpf}.

We shall apply Proposition \ref{negdiag} to a general morphism $\sigma : \sE^{\vee} \to \sF$, which gives rise to a smooth Calabi--Yau threefold $X_{\sF}$ with Picard number $2$ and a smooth surface $S_{\sF}$ in $X_{\sF}$.  For simplicity of notation, we continue to write $S$, $L$ and $H$ for $S_{\sF}$, $L_{\sF}|_{X_{\sF}}$ and $H_{\sF}|_{X_{\sF}}$ respectively.

\subsection{Rank Two Cases} \label{rk2subsec}

Suppose that $\sE$ and $\sF$ are vector bundles of rank two, and $\sF_{-} \cong \cO (b)$. Let $G = L - b H$ on $X_{\sF}$.

\begin{lemma} \label{contrrk2}
Assume that $G$ is base point free and big and $\rho (X_{\sF}) = 2$. There is a morphism $\varphi_{G} : X_{\sF} \to Y_{\sF}$ which contracts the exceptional divisor $S$ to a point, where $Y_{\sF}$ is a normal variety.
\end{lemma}

\begin{proof}
Let $\varphi_{G} : X_{\sF} \to Y_{\sF}$ be the Stein factorization of the morphism given by $|G|$. By assumption, $\varphi_{G}$ is birational. Observe that  $\sF_{-}$ is a line bundle and thus the natural projection $q : \bP(\sF_{-}) \xrightarrow{\sim} M$ is an isomorphism. From $\sF_{-} \cong \cO (b)$ and the formula \eqref{KF}, we see that
$$
q^{\ast} K_{M} \sim K_{\sF_{-}} \sim (q^{\ast} K_{M} + b H_{\sF_{-}}) - L_{\sF_{-}},
$$
where $K_{\sF_{-}} = K_{\bP(\sF_{-})}$ and $L_{\sF_{-}} = \cO_{\sF_{-}}(1).$
Hence $L_{\sF_{-}} \sim b H_{\sF_{-}}$, $G |_{S} = (L_{\sF} - b H_{\sF})|_{\bP (\sF_{-})} |_{S} \sim 0$, and the birational morphism $\varphi_G$ contracts $S$ to a point.
\end{proof}

Our result in the case of rank two is following:

\begin{theorem} \label{mainthmrk2}
Let $\sF = \cO (a) \oplus \cO (b)$ and $\sE = \cO (c) \oplus \cO (d)$ with $a \geqslant b > 0$ and $c \geqslant d \geqslant 0$.  Assume in addition that $(\sF, \sE)$ satisfies the Calabi--Yau condition \eqref{cycond}. Then, for a general morphism $\sigma : \sE^{\vee} \to \sF$, $X_{\sF}$ is a smooth Calabi--Yau threefold with Picard number $2$,
$$
\Nef(X_\sF)=\bR_{\geqslant 0}[L-bH]+\bR_{\geqslant 0}[H]
$$
and the determinantal contraction $\pi_{\sF}$ is induced by $|H|$.

The movable cone $\cMov(X_{\sF})$ is the convex cone generated by the divisors $L -b H$ and $r_M H - L$ and covered by the nef cones of $X_{\sF}$ and $X_{\sE}$. There are no more minimal models of $X_{\sF}$. Furthermore,
\begin{enumerate}[(i)]
    \item\label{mainthmrk2I} if $\sE$ is a trivial bundle and $a > b$, then $|L - b H|$ induces a primitive contraction $X_{\sF} \to Y_{\sF}$ of type II and the flop $X_{\sE} \to D_{1} (\sigma)$ of $\pi_{\sF}$ admits a K3 fibration induced by $|r_M H - L|$;

    \item\label{mainthmrk2II} if $a = b$, then $M$ is $\Gr (2, 4)$ or Mukai, and $X_{\sF}$, $X_{\sE}$ admit K3 fibrations induced by $|L - b H|$, $|r_M H - L|$ respectively;

    \item\label{mainthmrk2III} if $\sE$ is not trivial, then $M= \bP^4$ or $\Gr (2, 4)$ and $X_{\sF}$, $X_{\sE}$ admit primitive contractions of type II induced by $|L - b H|$, $|r_M H - L|$ respectively.
 \end{enumerate}
\end{theorem}

The following picture is $\cMov (X_{\sF})$ in $N^1 (X_{\sF})_{\bR}$. We depict $X_{\sF}$ and $X_{\sE}$ inside their nef cones. Note that $Y_{\sE} = \bP^1$ if $\sE$ is trivial.

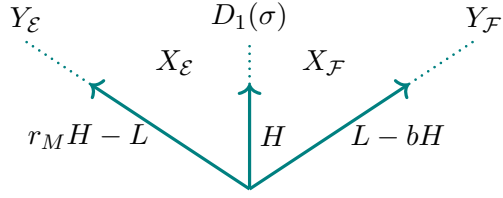
\begin{figure}[H]
    \centering
    \begin{tikzpicture}[xscale=0.7,yscale=0.7]
    %\draw[->, very thick, red] (0,0)--(-1.2,2);
    \draw[->, very thick, teal] (0,0)--(0,2);
    \draw[->, very thick, teal] (0,0)--(3,2);
    \draw[->, very thick, teal] (0,0)--(-3,2) ;

    %\node [left] at (-0.7,1) {$\textcolor{red}{L}$};
    \node [right] at (0,1) {$H$};
    \node [right] at (1.7,1) {$L - b H$};
    \node [left] at (-1.7,1) {$r_M H - L$};

    \node [above] at (0,2.8) {$D_1(\sigma)$};
    \node [above] at (4.4,2.8) {$Y_\sF$};
    \node [above] at (-4.2,2.8) {$Y_\sE$};
    \draw[decorate sep={.2mm}{1mm},fill, teal] (0,1)--(0,2.8);
    \draw[decorate sep={.2mm}{1mm},fill, teal] (3,2)--(4.2,2.8);
    \draw[decorate sep={.2mm}{1mm},fill, teal] (-3,2)--(-4.2,2.8);

    \node at (-1.4,2.4) {$X_{\sE}$};
    \node at (1.4,2.4) {$X_{\sF}$};
    \end{tikzpicture}
    %\caption{The movable cone of $X_{\sF}$}
    \label{fig:movrk2}
\end{figure}

\begin{proof}
By Proposition \ref{fanofdbpf}, the line bundle $\cO (1)$ of $M$ is globally generated and so are $\sF$ and $\sE$. Since $\bP (\sF (- b)) \cong \bP (\sF)$ and $a \geqslant b$, we see that $\cO_{\sF (-b)} (1)$ is globally generated by Proposition \ref{ggvb} and thus $G = L - b H$ is base point free. 

Choosing $\sF_{-} = \cO (b)$, there is a short exact sequence \eqref{Fminus}. It is easy to check that the assumptions of Proposition \ref{negdiag} are satisfied. Therefore, by the Calabi--Yau condition, $X_{\sF}$ is a smooth Calabi--Yau threefold with Picard number $\rho (\bP (\sF)) = 2$ and contains the smooth surface $S \in |L - a H|$ induced by $\sigma_{-} : \sE^{\vee} \to \sF_{-}$.

Assume that $\sE$ is trivial. Note that $a + b = r_M$ by the Calabi--Yau condition. From $\bP (\sE) = M \times \bP^1$, we have a diagram
\begin{center}
      \begin{tikzcd}
      X_\sF\arrow[dr,swap,"\pi_{\sF}"]\arrow[rr,dashed,"\chi"]& &X_{\sE}\arrow[dl,"\pi_{\sE}"]\arrow[dr] & \\
      & D_1(\sigma) &    &\bP^1
      \end{tikzcd}
\end{center}
where $X_{\sE} \to \bP^1$ is the restriction of $X_{\sE}$ to the second projection of $\bP (\sE)$.

If we can prove that the intersection number of $c_2 (T_{X_{\sE}})$ with the general fiber $F$ of $X_{\sE} \to \bP^1$ is $24$, then $F$ is a K3 surface (see \cite[Lemma 3.3]{Oguiso93}). Notice that $F \in |L_{\sE}|$ on $X_{\sE}$. By Proposition \ref{numbODP}, \ref{LHinter} and Lemma \ref{c2Fano}, we find that
\begin{align}
    \int_{X_{\sE}} c_2(T_{X_{\sE}}) \cdot L_{\sE} &= \int_M c_2 (T_M) \cdot c_2 (\sF - \sE^{\vee}) - c_2(\sF)^2 \notag\\
%    &= \int_M c_2 (T_M) . c_2 (\sF) - c_2(\sF)^2 \\
    &= (ab) \int_M c_2 (T_M) \cdot H_M^2 - (ab)^2 d_M \label{rk2K324}\\
%    &= a (r_M - a)(c_2 . H_M^2 - a (r_M - a) d_M) \\
    &= a (r_M - a) \left(d_M (a - 1) (a - r_M + 1) + \frac{24}{r_M - 1} \right) = 24, \notag
\end{align}
for $(r_M, a) = (2, 1)$ or $(3, 2)$. Hence $X_{\sE} \to \bP^1$ is a K3 fibration, and similarly for $(r_M, a) = (4, 2), (4, 3), (5, 3)$ or $(5, 4)$, where $d_M=2$ if $r_M=4$ (resp. $d_M=1$ if $r_M=5$).

In case $a = b$, we see that $M$ is $\Gr (2, 4)$ or Mukai and $a = 2$ or $1$ by Proposition \ref{cicylist}. Then the restriction $X_{\sF} \to \bP^1$ of $X_{\sF}$ to
\begin{equation}\label{rk2K3a=b}
    \bP (\sF) \cong \bP (\sF (- b)) = M \times \bP^1 \to \bP^1
\end{equation}
%$$
%\bP (\sF) \cong \bP (\sF (- b)) = M \times \bP^1 \to \bP^1
%$$
is also a K3 fibration. This can be proved in the same way as shown before.

On the other hand, by Proposition \ref{LHinter}, we have
\begin{equation}\label{rk2Gbig}
    G^3 = \sum_{k = 0}^3 (- b)^k \binom{3}{k} \int_M H_M^k \cdot s_{4 - k} (\sF^{\vee}) = a^2 (a - b)^2 d_M.
\end{equation}
%$$
%G^3 = \sum_{k = 0}^3 (- b)^k \binom{3}{k} \int_M H_M^k \cdot s_{4 - k} (\sF^{\vee}) = a^2 (a - b)^2 d_M.
%$$
In case $a > b$ (i.e., $G^3 > 0$), $G$ is big and there is a primitive contraction $\varphi_{G} \colon X_{\sF} \to Y_{\sF}$ of type II with exceptional set $S$ by Lemma \ref{contrrk2}.

To apply Lemma \ref{deterwall}, we need to verify the inequality \eqref{deterwallcond}. In case $\sE$ is trivial, it follows from the fact that $\sF$ is ample and
\begin{align*}
    \int_M \left(c_1 (\sE - \sF^{\vee})^2 - 2 c_2 (\sE - \sF^{\vee}) \right) \cdot H_M^2 &= \int_M \left(s_1 (\sF^{\vee})^2 - 2 s_2 (\sF^{\vee}) \right) \cdot H_M^2 \\
    &= - \left((a + b)^2 - 2 a b \right) d_M < 0.
\end{align*}
Hence the matrix of $\chi_*:N^1(X_\sF)\ra N^1(X_\sE)$ with respect to $\{L,H\}$ and $\{L_\sE,H_\sE\}$ is given by
$$
[\chi_*] =
\begin{bmatrix}
-1 & 0 \\
r_M & 1
\end{bmatrix}
=
[(\chi^{-1})_*],
$$
where the last equality is straight forward now. Therefore \eqref{mainthmrk2I} is established by \eqref{rk2K324}, \eqref{rk2Gbig} and the above geometric  argument, and similar for \eqref{mainthmrk2II} with the geometry in \eqref{rk2K3a=b}.

We now turn to the case $\sE \ncong \cO^{\oplus 2}$, that is, $r_M = 4$ or $5$ and
$$
(a, b, c, d) = (r_M - 2, 1, 1, 0).
$$
The inequality \eqref{deterwallcond} follows from a direct computation\footnote{Notice that $c_i (\sE) = 0$ for $i >1$ in this case.}:
\begin{itemize}
    \item $\int_{M} c_2 (\sF - \sE^{\vee}) \cdot H_M^2 = 2 (r_M - 1) d_M$,

    \item $\int_M \left(c_1 (\sE - \sF^{\vee})^2 - 2 c_2 (\sE - \sF^{\vee}) \right) \cdot H_M^2 = - (r_M - 2)^2 d_M$.
\end{itemize}
We can see that the base point free divisor $G \sim L - H$ on $X_{\sF}$ is big, which follows from
\begin{align}\label{G3}
    G^3 &= \sum_{k = 0}^3 (-1)^k \binom{3}{k} \int_M H_M^k \cdot \left(s_{4 - k} (\sF^{\vee}) + s_{3 - k} (\sF^{\vee}) c_1 (\sE) \right) \\
    &= (r_M - 2)^2 (r_M - 3)^2 d_M + (r_M - 2) (r_M - 3)^2 d_M \notag \\
    &= (r_M - 1) (r_M - 2) (r_M - 3)^2 d_M > 0 \notag
\end{align}
for $r_M = 4$ or $5$. By Lemma \ref{contrrk2}, there is a primitive contraction $\varphi_{G} \colon X_{\sF} \to Y_{\sF}$ of type II with exceptional set $S$.

Replacing the pair $(\sF, \sE)$ by $(\sE(1), \sF(-1))$, we can use the same argument as shown before to show that the linear system $|L_{\sE}|$ also induces a primitive contraction $X_{\sE} \to Y_{\sE}$ of type II. Hecne \eqref{mainthmrk2III} is established and the proof is completed.
\end{proof}

%\begin{remark}\label{rk2surfRmk}
%For a primitive contraction $X_{\sF} \to Y_{\sF}$ of type II, it is known that the exceptional set $S \subseteq X_{\sF}$ is a del Pezzo surface. Thus $S$ is determined\footnote{We also have $K_S \cdot H|_S = - (a - b)^{- 1} K_S^2$ in this case.} by
%$$
%K_S^2 = (a - b)^2 \left((a + b)(a + b + c) - a (b + r_M)\right) d_M,
%$$
%where $\sF = \cO (a) \oplus \cO (b)$ and $\sE = \cO (c) \oplus \cO$ as shown in the proof of Theorem \ref{mainthmrk2}: By adjuction and Proposition \ref{negdiag}, we have $K_S \sim  - (a - b) H|_{S}$. As in the proof of Lemma \ref{contrrk2}, we see that $L|_S \sim b H|_S$. The formula of $K_S^2$ follows from $K_S^2 = (a - b)^2 \left(H^2 \cdot (L - aH)\right)_{X_{\sF}}$ and Proposition \ref{LHinter}.
%\end{remark}

\begin{remark}\label{rk2surfRmk}
For a primitive contraction $X_{\sF} \to Y_{\sF}$ of type II, it is known that the exceptional set $S \subseteq X_{\sF}$ is a del Pezzo surface. We have a formula for the self-intersection of $K_S$:
$$
K_S^2 = (a - b)^2 \left((a + b)(a + b + c) - a (b + r_M)\right) d_M,
$$
where $\sF = \cO (a) \oplus \cO (b)$ and $\sE = \cO (c) \oplus \cO$. Indeed, we first observe that $L|_S \sim b H|_S$ as in the proof of Lemma \ref{contrrk2}. By adjunction and Proposition \ref{negdiag}, we get $K_S \sim  - (a - b) H|_{S}$. Then the formula follows from $K_S^2 = (a - b)^2 \left(H^2 \cdot (L - aH)\right)_{X_{\sF}}$ and Proposition \ref{LHinter}. On the other hand, we see that $K_S \cdot H|_S = - (a - b)^{- 1} K_S^2$. Note the Hirzebruch surface $\bF_1$ and $\bP^1 \times \bP^1$ are distinguished by $K_S \cdot H|_S$ being $- 5$ or $- 4$ (or by the Fano index $a - b$ of the surface being $1$ or $2$). Therefore the del Pezzo surface $S$ is determined by $K_S^2$ in our case.
\end{remark}

%Thus $S$ is determined\footnote{We also have $K_S \cdot H|_S = - (a - b)^{- 1} K_S^2$ in this case.} by

%the Fano index $r_S=a-b$ in $-K_S=(a-b)H|_S$ or $d_M=H|_S^2$ to distinguish $S=\bF_1$ ($r_S=1$) or $\bP^1\times\bP^1 (r_S=2).$

%$$
%K_S^2 = (a - b)^2 \left((a + b)(a + b + c) - a (b + r_M)\right) d_M,
%$$
%where $\sF = \cO (a) \oplus \cO (b)$ and $\sE = \cO (c) \oplus \cO$ as shown in the proof of Theorem \ref{mainthmrk2}: 
%By adjuction and Proposition \ref{negdiag}, we have $K_S \sim  - (a - b) H|_{S}$. As in the proof of Lemma \ref{contrrk2}, we see that $L|_S \sim b H|_S$. The formula of $K_S^2$ follows from $K_S^2 = (a - b)^2 \left(H^2 \cdot (L - aH)\right)_{X_{\sF}}$ and Proposition \ref{LHinter}.

\begin{remark}
In case $M = \bP^4$, $\sF = \cO (3) \oplus \cO (1)$ and $\sE = \cO (1) \oplus \cO$, we see that $S_{\sE}$ is isomorphic to a cubic surface in $\bP^3$ and the Calabi--Yau $X_{\sE}$ was studied in \cite[Theorem 5.5]{Kapustka209I}. On the other hand, $S = S_{\sF} \cong \bP^1 \times \bP^1$ and $X_{\sF}$ was studied in \cite[Section 2.3]{Kapustka09II}. In that paper, $X_{\sF}$ is the case of $\deg (X') = 5$ in Table $5$ and our divisor $L$ is $2 H^* + D$. Theorem \ref{mainthmrk2} tells us that these two Calabi--Yau threefolds $X_{\sF}$ and $X_{\sE}$ are connected by the flop $\chi$.

In case $M = \bP^4$ and $\sF = \cO (3) \oplus \cO (2)$, $X_{\sF}$ was studied in \cite[Theorem 5.3]{Kapustka209I}, and we note that the flop $\chi$ connects $X_{\sF}$ and the complete intersection $X_{\sE}$ of two hypersurfaces of bidegrees $(3, 1)$ and $(2, 1)$ in $\bP (\sE) = \bP^4 \times \bP^1$.
\end{remark}

\begin{remark}
In case that $M$ is del Pezzo and $\sF = \cO (2) \oplus \cO (1)$, we see that the del Pezzo surface $S$ has degree $2 \leqslant K_S^2 = d_M \leqslant 5$ by Remark \ref{rk2surfRmk}. The Calabi--Yau $X_{\sF}$ in the cases $d_M = 2, 4$ and $5$ were studied in \cite{Kapustka09II} (see the cases of Number $19, 6$, and $11$ in Table $1$ of that paper) and the remaining case in \cite[Remrak 5.9]{Kapustka209I}.
\end{remark}

\subsection{Rank Three Cases} \label{rk3subsec}

Assume that $\sF$ is of the form
\begin{equation}\label{Fbda11}
    \cO (a) \oplus \cO (1) \oplus \cO (1)
\end{equation}
and $\sE$, $\sF$ satisfy the Calabi--Yau condition. From Section \ref{CYsubsec}, there are four examples of such bundles. For $a = 1$, $M$ is a smooth del Pezzo fourfold and $\sE = \cO^3$. For $a = 2$, we have $(M, \sE) = (\bP^4, \cO (1) \oplus \cO^2)$ or $(\Gr (2, 4), \cO^3)$. And for $a = 3$, $(M, \sE) = (\bP^4, \cO^3)$. We see that $\sF_{-} = \cO (1)^2$ and $c_1 (\sE) = 0$ or $1$ in these cases.

Set $G = L - H$ on $X_{\sF}$. As in the proof of Theorem \ref{mainthmrk2}, we see that $G$ is base point free and hence nef.

We start with the cases $a = 2, 3$. Applying Proposition \ref{LHinter}, we get Table \ref{table_rk3}. In these cases, the top self-intersection number of $G$ is positive, and thus the nef divisor $G$ is big. %(see, e.g., \cite[Proposition 2.61]{KM98}). %\cite[Theorem 2.2.16]{LazPAGI}
\begin{table}[H]
\begin{center}
\begin{tabular}{ccccccc}
\toprule
M & $\sE$ & $L^3$ & $L^2 \cdot H$ & $L \cdot H^2$ & $H^3$ & $G^3$ \\ \midrule
$\bP^4$ & $\cO (1) \oplus \cO^2$ & 83 & 37 & 15 & 5 & 12 \\
 & $\cO^3$ & 179 & 58 & 18 & 5 & 54 \\
$\Gr (2, 4)$ & $\cO^3$ & 114 & 52 & 22 & 8 & 16 \\
\bottomrule
\end{tabular}
\end{center}
    \caption{The intersection numbers on $X_{\sF}$.}
    \label{table_rk3}
\end{table}

Let $\varphi_{G} : X_{\sF} \to Y_{\sF}$ be the Stein factorization of the morphism given by $|G|$, which is birational. Let $q_S :S \to \bP^1$ be the restriction to $S$ of the second projection $q : \bP (\sF_{-}) \cong M \times \bP^1 \to \bP^1$. Note that the corresponding divisor of $q^{\ast} \cO_{\bP^1} (1)$ is the divisor $L_{\sF_{-}} - H_{\sF_{-}}$ on $\bP (\sF_{-})$.

\begin{lemma} \label{rk3IIILem}
The birational morphism $\varphi_G : X_{\sF} \to Y_{\sF}$ determined by $|G|$ is primitive of type III. Moreover, $S$ is the exceptional divisor and $\varphi_G |_S = q_S$.
\end{lemma}

\begin{proof}
We first show that $(q_S)_{\ast} \cO_S \cong \cO_{\bP^1}$ and hence $q_S : S \to \bP^1$ has connected fibers. Consider the case $\sE = \cO^3$. According to the definition of $s_{\sigma_{-}}$, it follows that $S$ is the complete intersection of three smooth hypersurfaces $D_1$, $D_2$ and $D_3$ in $\bP (\sF_{-})$, where $D_i \in |L_{\sF_{-}}|$ for all $i$. Since $q : D_1 \to \bP^1$ is surjective, we see that $\cO_{\bP^1} \to q_{\ast} \cO_{D_1}$ is injective. Hence $\cO_{\bP^1} \xrightarrow{\sim} q_{\ast} \cO_{D_1}$ follows from the commutative diagram
\begin{equation*}
    \begin{tikzcd}
      \cO_{\bP^1} \ar[r, equal] \ar[d, "\rotatebox{270}{\(\sim\)}"] & \cO_{\bP^1} \ar[d, rightarrowtail] & \\
      q_{\ast} \cO_{\bP (\sF_{-})} \ar[r] & q_{\ast} \cO_{D_1} \ar[r] & R^1 q_{\ast} \cO (- D_1) = 0
    \end{tikzcd}
\end{equation*}
where the lower right corner is the relative Kodaira vanishing theorem. Consequently, $D_1 \to \bP^1$ has connected fibers. The same computation applies inductively to $D_2$, $D_3$ and hence the claim follows, and similarly for the case\footnote{In this case, $S = D_1 \cap D_2 \cap D_3$ where $D_1 \in |L_{\sF_{-}} - H_{\sF_{-}}|$ and $D_2, D_3 \in |L_{\sF_{-}}|$. We remark that $\cO (D_1) \cong q^{\ast} \cO_{\bP^1} (1)$.} $\sE = \cO^2 \oplus \cO (1)$.

By definition, $q_S$ is defined by $|G |_S|$, i.e., $G |_S$ is the divisor corresponding to $q_S^{\ast} \cO_{\bP^1} (1)$. Since $G - S \sim (a - 1) H$ is nef and big (for $a = 2$ or $3$), we get 
\[
    H^1 (X_{\sF}, \cO_{X_{\sF}}(G - S)) = 0
\] 
by Kawamata--Viehweg vanishing. Then $|G|_S$, the trace of $|G|$ on $S$, is the complete linear system $|G |_S|$. Hence every (connected) fiber of $q_S$ is contracted by $\varphi_G$ and $\varphi_G |_S = q_S$ (cf.~\cite[Proposition 1.14]{Debarre01}).

By Proposition \ref{negdiag}, the Picard number of $X_{\sF}$ is $\rho (\bP (\sF)) = 2$. Then the relative Picard number of $\varphi_G$ is $1$ and thus $\Exc (\varphi_G) = S$ \cite[Proposition 2.5]{KM98}.
\end{proof}

Our result in this case \eqref{Fbda11} is the following:

\begin{theorem}\label{mainthma11}
Let $\sF = \cO (a) \oplus \cO (1)^2$ and $\sE = \cO (c) \oplus \cO (d) \oplus \cO (e)$ with $a > 0$ and $c \geqslant d \geqslant e \geqslant 0$.  Assume in addition that $(\sF, \sE)$ satisfies the Calabi--Yau condition \eqref{cycond}. Then for a general morphism $\sigma : \sE^{\vee} \to \sF$, $X_{\sF}$ is a smooth Calabi--Yau threefold with Picard number $2$,
$$
\Nef(X_\sF)=\bR_{\geqslant 0}[L-H]+\bR_{\geqslant 0}[H]
$$
and the determinantal contraction $\pi_{\sF}$ is induced by $|H|$.

The movable cone $\cMov(X_{\sF})$ is the convex cone generated by the divisors $L - H$ and $r_M H - L$ and covered by the nef cones of $X_{\sF}$ and $X_{\sE}$ such that there are no more minimal models of $X_{\sF}$. Furthermore,
\begin{enumerate}[(i)]
    \item\label{mainthma11(i)} if $a > 1$, then $|L -  H|$ induces a primitive contraction $X_{\sF} \to Y_{\sF}$ of type III and the flop $X_{\sE} \to D_{1} (\sigma)$ of $\pi_{\sF}$ admits an elliptic fibration induced by $|r_M H - L|$ unless $(M, \sE) = (\bP^4, \cO (1) \oplus \cO^2)$, for which it has a primitive contraction $X_{\sE} \to Y_{\sE}$ of type III;

    \item\label{mainthma11(ii)} if $a = 1$, then $M$ is del Pezzo, and $X_{\sF}$, $X_{\sE}$ admit elliptic fibrations over $\bP^2$ induced by $|L - H|$, $|r_M H - L|$ respectively.
\end{enumerate}
\end{theorem}

The following picture is $\cMov (X_{\sF})$ in $N^1 (X_{\sF})_{\bR}$. We depict $X_{\sF}$ and $X_{\sE}$ inside their nef cones.

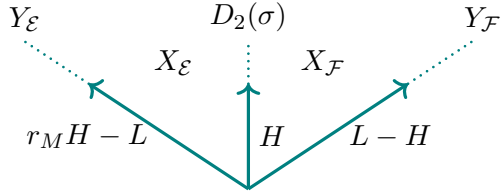
\begin{center}
\begin{tikzpicture}[xscale=0.7,yscale=0.7]
%\draw[->, very thick, red] (0,0)--(-1.2,2);
\draw[->, very thick, teal] (0,0)--(0,2);
\draw[->, very thick, teal] (0,0)--(3,2);
\draw[->, very thick, teal] (0,0)--(-3,2) ;

%\node [left] at (-0.7,1) {$\textcolor{red}{L}$};
\node [right] at (0,1) {$H$};
\node [right] at (1.7,1) {$L - H$};
\node [left] at (-1.7,1) {$r_M H - L$};

\node [above] at (0,2.8) {$D_2(\sigma)$};
\node [above] at (4.4,2.8) {$Y_\sF$};
\node [above] at (-4.2,2.8) {$Y_{\sE}$};

\draw[decorate sep={.2mm}{1mm},fill, teal] (0,1)--(0,2.8);
\draw[decorate sep={.2mm}{1mm},fill, teal] (3,2)--(4.2,2.8);
\draw[decorate sep={.2mm}{1mm},fill, teal] (-3,2)--(-4.2,2.8);

\node at (-1.4,2.4) {$X_{\sE}$};
\node at (1.4,2.4) {$X_{\sF}$};
\end{tikzpicture}
\end{center}

\begin{proof}
As in the proof of Theorem \ref{mainthmrk2}, we can verify that the inequality \eqref{deterwallcond} in Proposition \ref{deterwall} holds. For example, in the case $(M, \sE) = (\bP^4, \cO (1) \oplus \cO^2)$,
$$
c_2 (\sF - \sE^{\vee}) = c_2 (\sF) + c_1(\sF)c_1(\sE) + c_1 (\sE)^2 = 10 H_{M}^2.
$$
From Lemma \ref{deterwall}, we see that the image of $L_{\sE}$ in $N^1(X_{\sF})$ is $r_M H - L$.

For $a = 1$, $M$ is del Pezzo and $\sE$ is trivial by Proposition \ref{cicylist}. Then $|L_{\sE}|$ induces an elliptic fibration on $X_{\sE}$ which is the restriction of the natural projection $\bP (\sE) = M \times \bP^2 \to \bP^2$ to $X_{\sE}$. Similarly, $\bP (\sF(-1)) = M \times \bP^2$ and $|L - H|$ induces an elliptic fibration on $X_{\sF}$.

For $a > 1$, \eqref{mainthma11(i)} follows from Lemma \ref{rk3IIILem}. Notice that if $(M, \sE) = (\bP^4, \cO (1) \oplus \cO^2)$, then $\sE \cong \sF (1)$.
\end{proof}

Note that we could also characterize the exceptional surface $S$. For abbreviation, we let $\tbP^2 (r)$ stand for the blow-up of $\bP^2$ in the points $x_1, \cdots, x_r$, which can be infinitely near. 

%Let $\tbP^2 (x_1, \cdots, x_r)$ denote the blow-up of $\bP^2$ in the points $x_1, \cdots, x_r$, which can be infinitely near.

\begin{proposition} \label{StypeIII}
Let $S$ be the smooth surface as in Lemma \ref{rk3IIILem}. Then
\begin{equation*}
    S \cong
    \begin{cases}
    \tbP^2 (9 - d_M) & \mbox{if } a = 2,\\
    \tbP^2 (1) & \mbox{if } a = 3.
    \end{cases}
\end{equation*}
\end{proposition}
%\begin{equation*}
%    S \cong
%    \begin{cases}
%    \tbP^2 (x_1, \cdots, x_{9 - d_M}) & \mbox{if } a = 2,\\
%    \tbP^2 (x_1) & \mbox{if } a = 3.
%    \end{cases}
%\end{equation*}

\begin{proof}
As we have seen in the proof of Lemma \ref{rk3IIILem}, $q_S : S \to \bP^1$ has connected fibers. Let $\ell$ be the fiber class of $q_S$. Recall that $K_S \sim S|_S$ and $S \sim L - a H$. Since $\ell^2=0$ on $S$ and 
\begin{align*}
   \left(K_S \cdot \ell\right)_S &= \left((L - a H)^2 \cdot (L - H)\right)_{X_{\sF}} \\
   &= (1 + c_1 (\sE)) (1 - a) d_{M} = - 2,
\end{align*}
it implies that a general fiber of $q_S$ is a smooth rational curve and hence $S$ is rational. On the other hand, we have 
\begin{align*}
    K_S^2 &= \left((L - a H)^3\right)_{X_{\sF}} \\
    &= \left( (3 a - 5) (a - 1) + 2 c_1 (\sE) (a - 1)(a - 2) \right) d_{M}\\
    &=
    \begin{cases}
    d_M & \mbox{if } a = 2, \\
    8 & \mbox{if } a = 3.
    \end{cases}
\end{align*}
To prove the proposition, it remains to show that if $a=3$, then $S\cong\bF_1.$

By running a relative minimal model program of $q_S:S\rightarrow\bP^1$ over $\bP^1$, there is an $n\in\bZ_{\geqslant 0}$ and a birational morphism $S\rightarrow \bF_n$ over $\bP^1$ consisting of $m$ finitely many smooth blow-downs with $K_S^2=K_{\bF_n}^2-m=8-m.$
Hence for $a=3$, we get $m=0$ and $K_S=K_{\bF_n}=-2C_n-(n+2)l$, where $C_n$ satisfying $C_n^2=-n$ is the unique negative section over $\bP^1.$ As $$(H|_S \cdot \ell)_S = \left(H\cdot(L-H)\cdot(L-3H)\right)_{X_{\sF}}=1.$$
and 
\begin{equation*}
    \left(K_S \cdot H |_S\right)_S = \left((L - 3 H)^2 \cdot H \right)_{X_{\sF}} = -5,
\end{equation*}
we have 
$$C_n\cdot H|_S=\frac{n-3}{-2}\in\bZ_{>0},$$
and $n=1$ is the unique possibility. Hence $S=\tbP^2 (1)\cong D_{n-1}(\sigma_-) \hookrightarrow \bP^4$, via the very ample linear system $|H|_S|$ (see the diagram \eqref{basic_diag} in Proposition \ref{negdiag}), is a rational scroll of degree $(H|_S)^2=3$ and $H|_S\sim C_1+2\ell$, i.e., $|H|_S|$ is the linear system of quadrics on $\bP^2$ passing through a fixed point $x_1$.
\end{proof}

We now deal with the remaining case $M = \bP^4$ and $\sF = \cO (2)^{2} \oplus \cO (1)$. In this case, we see that $\sE = \cO^{3}$, $\sF_{-} = \cO (1) \oplus \cO (2)$ and
$$
q \colon \bP (\sF_{-}) \cong \bP (\cO \oplus \cO (1)) = \Bl_{o} \bP^5 \to \bP^5
$$
is the blow-up of $\bP^5$ at a point $o$ (cf.~\cite[Example V.2.11.4]{Har:AG}). Applying Proposition \ref{numbODP} and \ref{LHinter}, we get Table \ref{table_221}.
\begin{table}[H]
  \centering
  \begin{tabular}{ccccccc}
    \toprule
    $L^3$ & $L^2 \cdot H$ & $L \cdot H^2$ & $H^3$ & $L \cdot c_2(T_Z)$ & $H \cdot c_2(T_Z)$ & $\#$ of ODPs \\
    \midrule
     129 & 49 & 17 & 5 & 126 & 50 & 44 \\
    \bottomrule
  \end{tabular}
  \caption{The intersection numbers on $X_{\sF}$.}
  \label{table_221}
\end{table}

Set $S = S_{\sF}$. Let us denote by $q_S : S \to \bP^5$ the restriction of $q$ to $S$ and by $S_0$ its image. Observe that the pullback divisor of the hyperplane class $q^{\ast} H_{\bP^5}$ is the divisor $L_{\sF_{-}} - H_{\sF_{-}}$ on $\bP (\sF_{-})$.

\begin{lemma}\label{K3lem}
$S_0 \subseteq \bP^5$ is a K3 surface of degree $8$. Moreover, $q_S$ is the blow-up of $S_0$ at a point $o$.
\end{lemma}

\begin{proof}
Let $E \subseteq \bP (\sF_{-})$ be the exceptional divisor of $q$. From relative Euler sequence, $\det \sF_{-} \cong \cO (3)$ and $\bP (\sF_{-}) \cong \Bl_{o} \bP^5$, we see that
$$
- 2 H_{\sF_{-}} - 2 L_{\sF_{-}} \sim K_{\sF_{-}} \sim - 6 (L_{\sF_{-}} - H_{\sF_{-}}) + 4 E
$$
and thus $E \sim L_{\sF_{-}} - 2 H_{\sF_{-}}$.

Let $H_S$ and $C$ denote the restrictions of $H_{\sF_{-}}$ and $E$ to $S$ respectively. Notice that $H_S$ is very ample, because it is the pullback of the hyperplane class $H_{\bP^4}$ on $\bP^4$ via $S \xrightarrow{\sim} D_{1} (\sigma_{-})$ (see the commutative diagram \eqref{basic_diag}). Since $(C \cdot H_S)_S = \left( (L - 2 H)^2 \cdot H \right)_{X_{\sF}} = 1$, $C$ is a straight line. Moreover, it is a $(-1)$-curve on $S$ because $\left(C^2\right)_S= \left( (L - 2 H)^3 \right)_{X_{\sF}} = - 1$.

We first observe that the irregularity of $S$ is zero. Indeed, $K_S \sim C$ by $K_{X_{\sF}} \sim 0$, $S \sim L- 2 H$ and the adjunction formula. From $\cO_{C} (C) \cong \cO_{\bP^1} (- 1)$ and
\begin{equation*}
    0 \to \cO_S \to \cO_S (C) \to \cO_{C} (C) \to 0,
\end{equation*}
it implies that $p_g (S) = h^0 (\cO_S) = 1$. By the double point formula \cite[p.~434]{Har:AG} and $(H_S^2)_{S} = (H^2 \cdot (L - 2 H))_{X_{\sF}} = 7$, we see that $\chi (\cO_S) = 2$ and hence $h^1 (\cO_S) = 0$.

%\footnote{This indeed almost follows from the proof of Castelnuovo's criterion in \cite{Har:AG}.}

We claim that $S_0$ is smooth and $q_S$ is the blow-up of $S_0$ at $o$. Observe that
$$
H_S + C \sim (L_{\sF_{-}} - H_{\sF_{-}}) |_S \sim q_S^{\ast} H_{\bP^5}.
$$
From $(H_S + C) \cdot C = 0$, we see that $q_S$ must map $C$ to a point $o$ and $H^0 (\cO_S (H_S + C)) \to H^0 (\cO_C (H_S + C)) \cong \bC$ is surjective. On the other hand, since $H_S$ is very ample, the linear system $|H_S + C|$ separates points and tangent vectors away from $C$, and also separates points of $C$ from points not on $C$, so $q_S : S \setminus C \xrightarrow{\sim} S_0 \setminus \{o\}$.

If we prove that $H^1 (\cO_S (H_S - C)) = 0$, then the claim follows from the step $7$ in the proof of \cite[Theorem V.5.5]{Har:AG}. Consider the exact sequences
\begin{equation} \label{122exactseq}
    0 \to \cO_S ( H_S + \ell C) \to \cO_S (H_S + (\ell + 1) C) \to \cO_{C} (H_S + (\ell + 1) C) \to 0
\end{equation}
for $\ell = -1, 0$. By the long exact sequence in cohomology for \eqref{122exactseq} with $\ell = 0$ and Kodaira vanishing, we see that $H^1 (\cO_S (H_S)) = 0$. Note that
$$
H^0 (\cO_S (H_S)) \to H^0 (\cO_C (H_S))
$$
is surjective. In fact, we already know $\cO_C (H_S) \cong \cO_{\bP^1} (1)$. Given any $D$ belonging to the very ample linear system $|H_S|$ that is either tangent to the straight line $C$ or contains two points of $C$, then $C \subseteq D$. Therefore $h^0 (\cO_S (H_S - C)) = h^0 (\cO_S (H_S)) - 2$.

%This is obvious since $H_S . C = 1$ and $H_S$ is very ample.\footnote{Given any $D \in |H_S|$ that is either tangent to the straight line $C$ or contains $2$ of $C$'s points, then $C \subseteq D$. Hence $h^0 (\cO_S (H_S - C)) + 2 = h^0 (\cO_S (H_S))$.}

From above facts and the long exact sequence  \eqref{122exactseq} with $\ell = - 1$, we see that
$H^1 (\cO_S (H_S - C)) = 0$. Consequently, $S_0$ is a smooth surface of degree $(H_S + C)^2 = 8$ in $\bP^5$ with $K_{S_0} \sim 0$ and $h^1 (\cO_{S_0}) = h^1 (\cO_S) = 0$. The proof is complete.
\end{proof}

\begin{remark}
We know that $S \cong D_2 (\sigma_{-})$ is a smooth surface in $\bP^4$ of degree $7$. The structure of $S$ was studied by Okonek \cite[Theorem 6]{Okonek84}. The proof given above is to verify that the adjunction map defined by $|K_S + H_S|$ is just the natural projection $q_S$.
\end{remark}

Set $G = L - H$ on $X_{\sF}$. As in the proof of Theorem \ref{mainthmrk2}, we see that $G$ is base point free and $G^3 = 28 > 0$. Let $\varphi_{G} : X_{\sF} \to Y_{\sF}$ be the Stein factorization of the morphism given by $|G|$, which is birational.

\begin{lemma}\label{K3prilem}
The birational morphism $\varphi_G : X_{\sF} \to Y_{\sF}$ is a small contraction. Moreover, $Y_{\sF}$ is a nodal Calabi--Yau threefold with one ODP and $\varphi_G |_S = q_S$.
\end{lemma}

\begin{proof}
We first observe that the locus $\Exc (\varphi_G)$ is contained in $S$. Indeed, if $C'$ is a curve contracted by $\varphi_G$, then $(L - H) \cdot C' = 0$ on $X_{\sF}$. Therefore,
\[
S \cdot C' = (L - 2 H) \cdot C' = - L \cdot C' < 0
\]
which implies that $C'$ is contained in $S$.

Since $G - S \sim H$ is nef and big, we get $H^1 (\cO_{X_{\sF}}(G - S)) = 0$ by Kawamata--Viehweg vanishing. Then we see as in the proof of Lemma \ref{rk3IIILem} that $\varphi_G |_S = q_S$. By Lemma \ref{K3lem}, $q_S$ is the blow-up of the K3 surface $S_0$ at a point $o$. Combining these with $\Exc(\varphi_G) \subseteq S$, the exceptional sets of $\varphi_G$ and $q_S$ are the same, which consists of one rational curve $C$. Thus $\varphi_G$ is small.

It remains to show that $Y_{\sF}$ has only one ODP, i.e., $N_{C / X_{\sF}} \cong \cO_{\bP^1} (- 1)^2$. Since $C$ is a $(- 1)$-curve in the smooth surface $S$, this last claim follows from Lemma \ref{ODPsmsurf} and the proof is completed.
\end{proof}

\begin{theorem}\label{mainK3}
Let $(M, \sE, \sF) = (\bP^4, \cO^{3}, \cO (2)^{2} \oplus \cO (1))$. Then for a general morphism $\sigma : \sE^{\vee} \to \sF$, $X_{\sF}$ is a smooth Calabi--Yau threefold with Picard number $2$ and
$$
\Nef(X_\sF)=\bR_{\geqslant 0}[L-H]+\bR_{\geqslant 0}[H],
$$
such that
\begin{enumerate}[(i)]
    \item\label{mainK31} the determinantal contraction $\pi_{\sF}$ is induced by $|H|$;

    \item\label{mainK32} $|L -  H|$ induces a primitive contraction $X_{\sF} \to Y_{\sF}$ of type I and $Y_{\sF}$ is a Calabi--Yau threefold with exactly one ODP singular point;

    \item\label{mainK33} the flop $X_{\sF}^{+} \to Y_{\sF}$ of $X_{\sF} \to Y_{\sF}$ admits a K3 fibration induced by $|L - 2 H|$;

    \item\label{mainK34} $X_{\sE}$ admits an elliptic fibration over $\bP^2$ induced by $|5 H - L|$.
\end{enumerate}
Moreover, the movable cone $\cMov(X_{\sF})$ is the convex cone generated by the divisors $L - 2 H$ and $5 H - L$ and covered by nef cones of $X_{\sF}$, $X_{\sF}^{+}$, and $X_{\sE}$ such that there are no more minimal models of $X_{\sF}$.
\end{theorem}

The following picture is $\cMov (X_{\sF})$ in $N^1 (X_{\sF})_{\bR}$. We depict $X_{\sF}^+$, $X_{\sF}$, and $X_{\sE}$ inside their nef cones.

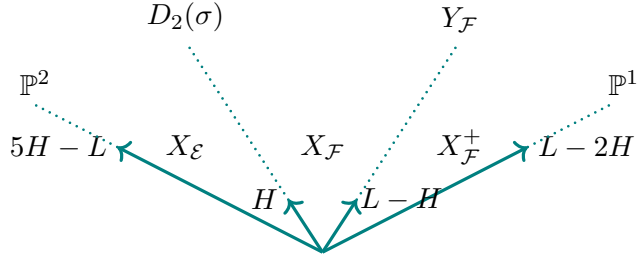
\begin{center}
\begin{tikzpicture}[xscale=0.9,yscale=0.7]
\draw[->, very thick, teal] (0,0)--(.5,1);
\draw[->, very thick, teal] (0,0)--(-.5,1);
\draw[->, very thick, teal] (0,0)--(3,2);
\draw[->, very thick, teal] (0,0)--(-3,2) ;

\node [right] at (.4,1) {$L - H$};
\node [left] at (-.5,1) {$H$};
\node [right] at (3,2) {$L - 2 H$};
\node [left] at (-3,2) {$5 H - L$};

\node [above] at (2,4) {$Y_{\sF}$};
\node [above] at (-2,4) {$D_2(\sigma)$};
\node [above] at (4.4,2.8) {$\bP^1$};
\node [above] at (-4.2,2.8) {$\bP^2$};

\draw[decorate sep={.2mm}{1mm},fill, teal] (.5,1)--(2,4);
\draw[decorate sep={.2mm}{1mm},fill, teal] (-.5,1)--(-2,4);
\draw[decorate sep={.2mm}{1mm},fill, teal] (3,2)--(4.2,2.8);
\draw[decorate sep={.2mm}{1mm},fill, teal] (-3,2)--(-4.2,2.8);

\node at (0,2) {$X_{\sF}$};
\node at (2,2) {$X_{\sF}^{+}$};
\node at (-2,2) {$X_{\sE}$};
\end{tikzpicture}
\end{center}

\begin{proof} The statement \eqref{mainK31} is obvious. By Proposition \ref{negdiag}, $X_{\sF}$ is a smooth Calabi--Yau threefold with Picard number $\rho (\bP (\sF)) = 2$ and the relative Picard number $\rho (X_{\sF} / Y_{\sF})$ is $1$. Then $X_{\sF} \to Y_{\sF}$ is primitive and \eqref{mainK32} follows from Lemma \ref{K3prilem}.

By Lemma \ref{K3lem}, the surface $S$ in $X_{\sF}$ is the smooth blow-up at one point $o\in S_0$. According to Lemma \ref{K3prilem} and \ref{A1flop}, it follows that the Atiyah flop $X_{\sF}^+$ contains the minimal model of $S$, which is isomorphic to the K3 surface $S_0$. By abuse of notation, we continue to write $S_0$ for the K3 surface in $X_{\sF}^+$. From Proposition \ref{divcontr}, we see that the linear system $|S_0|$ determines a fibration $X_{\sF}^+ \to \bP^1$ with $S_0$ as a fiber. Then \eqref{mainK33} follows from the fact that $S \sim L - 2 H$.

The restriction of $\bP(\sE) = \bP^4 \times \bP^2 \to \bP^2$ to $X_{\sE}$ gives rise to a elliptic fibration on $X_{\sE}$ over $\bP^2$, which is induced by $|L_{\sE}|$. From the fact that $\sF$ is ample and
\begin{align*}
    \int_M \left(c_1 (\sE - \sF^{\vee})^2 - 2 c_2 (\sE - \sF^{\vee}) \right) \cdot H_M^2 &= \int_M \left(s_1 (\sF)^2 - 2 s_2 (\sF) \right) \cdot H_M^2 \\
    &= - 9 < 0,
\end{align*}
we see that \eqref{mainK34} holds
by Proposition \ref{deterwall}.
\end{proof}

\begin{remark}
$X_{\sE}$ is a smooth complete intersection of smooth hypersurfaces of bidegrees $(2, 1), (2, 1)$ and $(1, 1)$ in $\bP^4 \times \bP^2$.
\end{remark}

\subsection{A Rank Four Case}

We consider $M = \bP^4$ and $\sF = \cO (2) \oplus \cO (1)^3$. In this case, we see that $\sE = \cO^4$ and $\sF_{-} = \cO (1)^3$. Applying Proposition \ref{LHinter}, we have Table \ref{table_rk4} (cf.~\cite[Lemma 3.2]{CynkRams2012}).
\begin{table}[H]
  \centering
  \begin{tabular}{ccccccc}
    \toprule
    $L^3$ & $L^2 \cdot H$ & $L \cdot H^2$ & $H^3$ & $L \cdot c_2(T_Z)$ & $H \cdot c_2(T_Z)$ & $\#$ of ODPs \\
    \midrule
     99 & 42 & 16 & 5 & 114 & 50 & 46 \\
    \bottomrule
  \end{tabular}
  \caption{The intersection numbers on $X_{\sF}$.}
  \label{table_rk4}
\end{table}

Let $S = S_{\sF}$ and $q_S :S \to \bP^2$ be the restriction to $S$ of the second projection $q : \bP (\sF_{-}) \cong M \times \bP^2 \to \bP^2$. Note that $q^{\ast} \cO_{\bP^2} (1)= \cO (L_{\sF_{-}} - H_{\sF_{-}})$ on $\bP (\sF_{-})$.

\begin{lemma} The surface $S$ is a Bordiga surface, i.e., $S \cong \tbP^2 (10)$. Moreover, $q_S:S\ra\bP^2$ is the blow-up of $\bP^2$ in ten distinct points if $X_\sF$ is chosen in general.
\end{lemma}

\begin{proof}
From the definition of $s_{\sigma_{-}}$, we find that $S$ is the complete intersection of four smooth hypersurfaces $D_i$ in $\bP (\sF_{-})$, where $D_i \in |L_{\sF_{-}}|$ for $i = 1, \cdots, 4$. Using the same argument as in the proof of Lemma \ref{rk3IIILem}, we can show that $q_S : S \to \bP^2$ has connected fibers.

According to the adjunction formula, $K_{X_{\sF}} \sim 0$ and $S \sim L - 2 H$, it follows that $K_S \sim (L - 2 H)|_S$. Therefore
$$
K_S^2 = \left( (L - 2 H)^3 \right)_{X_{\sF}} = - 1 = K_{\bP^2}^2 - 10,
$$
and $q_S:S\ra\bP^2$ consists of ten smooth blow-ups, which might contain infinitely near points.

To finish the proof, simply observe that by construction the subscheme $q_S^{-1}(x)$ of dimension at most one is cut out by linear equations in $\bP^4$ for each $x\in\bP^2$. Hence if $X_\sF$ is general enough, then  $E_x\coloneqq q_S^{-1}(x)\cong\bP^1$ and $K_S \cdot E_x=-1$. In particular, $S$ is the blow-up of ten distinct points on $\bP^2$.
\end{proof}

Set $G = L - H$ on $X_{\sF}$. From $\sF = \cO (2) \oplus \cO (1)^3$, we see that $G$ is base point free and big ($G^3 = 16 > 0$). Then $|G|$ determines a birational morphism $\varphi_G : X_{\sF} \to Y_{\sF}$, where $Y_{\sF}$ is a normal variety.

\begin{lemma} \label{r4Fcontr} The birational morphism
$\varphi_G : X_{\sF} \to Y_{\sF}$ is a small contraction onto a Calabi--Yau threefold $Y_{\sF}$ with $10$ ODPs and $\varphi_G |_S = q_S$.
\end{lemma}

\begin{proof}
Let $C\subseteq X_\sF$ be an integral curve contracted by $\varphi_G$. If $C\nsubseteq  S$, then $(L-H) \cdot C = S \cdot C + H \cdot C  \geqslant 0$ and equality holds only if $S \cdot C = H \cdot C = 0$. But then $\bR_{\geqslant} [C]$ is the extremal ray of $X_\sF\ra D_3 (\sigma)$, which is absurd as $(L-H) \cdot C =L\cdot C> 0$ in this case. If $C \subseteq S$, then $(L-H) \cdot C = q_S^{\ast} H_{\bP^2} \cdot C = 0$ only when $C$ is $q_S$-exceptional.

Recall that $q_S:S\ra\bP^2$ is the blow-up of ten distinct points on $\bP^2$. Then by Lemma \ref{ODPsmsurf}, each irreducible exceptional curve $E_i\cong\bP^1$ has normal bundle $N_{E_i/X_\sF}\cong\cO_{\bP^1}(-1)^2$ and is contracted to an ODP on $Y_\sF$.

\end{proof}

\begin{remark} \label{4QintP7}
The nodal Calabi--Yau threefold $Y_{\sF}$ is an intersection of four quadrics in $\bP^7$. Indeed, by Riemann--Roch and Kawamata--Viehweg vanishing theorem, we have
$$
h^0(\cO_{X_{\sF}} (k G)) = \frac{k^3}{6} G^3 + \frac{k}{12} G \cdot c_2(T_{X_{\sF}}) = \frac{8 k^3 + 16 k}{3}.
$$
From $h^0(\cO_{X_{\sF}} (G)) = 8$, we see that the linear system $|G|$ defines a morphism $X_{\sF} \to \bP^7$. One can show that the image of this morphism is projectively normal. According to $h^0(\cO_{X_{\sF}} (2 G)) = 32$ and $h^0(\cO_{\bP^7}(2)) = 36$, we find four quadrics $Q_0, Q_1, Q_2, Q_3$ containing $Y_{\sF}$. Then $Y_{\sF}$ must be equal to the intersection of these quadrics, which is a threefold of degree $16$.
\end{remark}

Consider the natural map $\phi : X_{\sE} \to \bP^3$, which is the restriction of the second projection  $\bP (\sE) = \bP^4 \times \bP^3 \to \bP^3$ to $X_{\sE}$.

\begin{lemma}[\cite{CynkRams2012}] \label{Bcover_contra}
Let $X_{\sE} \to Y_{\sE} \to \bP^3$ be the Stein factorization of $\phi$. Then the morphism $\phi$ is generically $2:1$ and $\widehat{\phi} : X_{\sE} \to Y_{\sE}$ is a small contraction if $X_{\sE}$ is chosen in general.
\end{lemma}

\begin{proof}
Applying Proposition \ref{LHinter}, we have\footnote{This also can be computed by the fact that $X_\sE$ is a complete intersection of hypersurfaces of degree $(2,1)$ and $3\times (1,1)$ in $\bP^4 \times \bP^3$.} that $L_\sE^3 = 2$, $L_\sE^2 \cdot H_\sE = 7$, $L_\sE \cdot H_\sE^2 = 9$, and $H_\sE^3$ = 5 (cf.~\cite[Lemma 4.1]{CynkRams2012}). Then, by $L_{\sE}^3 = 2$, the surjective morphism $X_{\sE} \to \bP^3$ is generically $2:1$.

To see that $X_{\sE} \to Y_{\sE}$ is small, we recall the description  of the discriminant locus of $\phi$ in \cite{Michalek12, CynkRams2012}. By Remark \ref{4QintP7}, $Y_{\sF} = Q_0 \cap O_1 \cap Q_2 \cap Q_3$, where $Q_i$ is a quadric in $\bP^7$. We define a degree $8$ surface in $\bP^3$:
\[
S_8 \coloneqq \left\{ y \in \bP^3 \,\middle\vert\, \det \left( \sum_{i = 0}^3 y_i \fq_i\right) = 0 \right\},
\]
%\[
%S_8 \coloneqq \{ y \in \bP^3 \mid \det \left( \textstyle\sum_{i = 0}^3 y_i \fq_i\right) = 0\},
%\]
where $Q_i$ is given by the symmetric $8 \times 8$ matrix $\fq_i$. Note that each $Q_i$ contains a fixed plane $\bP^2$ by Lemma \ref{r4Fcontr}. Since $\sigma : \sE^{\vee} \to \sF$ is chosen in general, $S_8$ has only $94$ isolated singular points by \cite[Theorem 2.7]{Michalek12}. From \cite[Theorem 4.6]{CynkRams2012}, we see that the double cover $X_{\sE} \to \bP^3$ is branched along the surface $S_8$ and $X_{\sE} \to Y_{\sE}$ is a small resolution. Note that the set of one-dimensional fibers of $\widehat{\phi}$ coincides with $\widehat{\phi}^{-1} (\Sing (Y_{\sE}))$.
\end{proof}

\begin{remark}\label{Bcover_ODPs}
When the determinantal octic $S_8$ has only $94$ isolated singular points, the $94$ singular points of $Y_{\sE}$ are all ODPs (cf.~\cite[Corollary 5.7]{CynkRams2012}).

Indeed, by \cite[Corollary 2.12]{Michalek12} (or Section \ref{hodgenumsubsec}), we have the Euler number $\chi_{top} (X_{\sE}) = -108$. Let $\widetilde{Y}$ be a double cover of $\bP^3$ branched over a smooth octic surface $\widetilde{S}_8$. Then the Euler number
$\chi_{top}(\widetilde{Y}) = 2 \chi_{top} (\bP^3) - \chi_{top}(\widetilde{S}_8) = -296$. Therefore we get
$$
\chi_{top} (X_{\sE}) - \chi_{top}(\widetilde{Y}) = 188= 2 |\Sing (Y_{\sE})|,
$$
and $Y_{\sE}$ is a nodal Calabi--Yau threefold by Proposition \ref{compEuler}.
\end{remark}

\begin{lemma} \label{invol_matrix}
Let $\iota:X_\sE\dra X_\sE$ be the involution over $\bP^3$. Then with respect to $\{L_\sE,H_\sE\}$, the matrix representation of $\iota$ on $N^1(X_\sE)$ is $$
[\iota_*] =
\begin{bmatrix*}[c]
1 & 7 \\
0 & -1
\end{bmatrix*}
=
[(\iota^{-1})_*].
$$
\end{lemma}

\begin{proof}
Note that $\iota_*L_\sE=L_\sE$ and write $\iota_*H_\sE=xL_\sE+yH_\sE$. Since $\iota\in\Bir(X_\sE)$ is small, we have
$$\begin{cases}
L_\sE \cdot H_\sE^2=L_\sE \cdot (xL_\sE+yH_\sE)^2\\
L_\sE^2 \cdot H_\sE=xL_\sE^3+yL_\sE^2 \cdot H_\sE\\
\end{cases}\ {\rm or}\
\begin{cases}
9=2x^2+14xy+9y^2\\
7=2x+7y\\
\end{cases}.$$
Since $(x,y)=(0,1)$ is impossible, the only solution is $(x,y)=(7,-1)$.
\end{proof}

\begin{theorem}\label{mainBordiga}
Let $(M, \sE, \sF) = (\bP^4, \cO^{4}, \cO (2) \oplus \cO (1)^3)$. Then for a general morphism $\sigma : \sE^{\vee} \to \sF$, the scheme $X_{\sF}$ is a smooth Calabi--Yau threefold with Picard number $2$ with
$$
\Nef(X_\sF)=\bR_{\geqslant 0}[L-H]+\bR_{\geqslant 0}[H],
$$
such that
\begin{enumerate}[(i)]
    \item\label{mainB1} the determinantal contraction $\pi_{\sF}$ is induced by $|H|$;

    \item\label{mainB2} $|L -  H|$ induces a primitive contraction $X_{\sF} \to Y_{\sF}$ of type I and $Y_{\sF}$ is a Calabi--Yau threefold with $10$ ODPs;

    \item\label{mainB3} $|5 H - L|$ induces a primitive contraction $X_{\sE} \to Y_{\sE}$ of type I, $Y_{\sF}$ is a Calabi--Yau threefold with $94$ ODPs, and the double cover $X_{\sE} \to \bP^3$ factors through $Y_{\sE}$;

    \item\label{mainB4} for the flop $X_{\sF}^{+} \to Y_{\sF}$ of $X_{\sF} \to Y_{\sF}$, $X_{\sF}^{+}$ admits a primitive contraction $X_{\sF}^{+} \to Z_{\sF}$ of type II induced by $|4 L - 5 H|$.
\end{enumerate}
Moreover, the movable cone $\cMov(X_{\sF})$ is the convex cone generated by the divisors $4L - 5 H$ and $490 H - 101 L$ which is covered by the nef cones of $X_{\sF}$, $X_{\sF}^{+}$ and $X_{\sE}$, and there are no more minimal models of $X_{\sF}$.
\end{theorem}

The movable cone decomposition of $X_\sF$ is given by the following diagram, where the left hand side is given by the mirror of the right hand side:
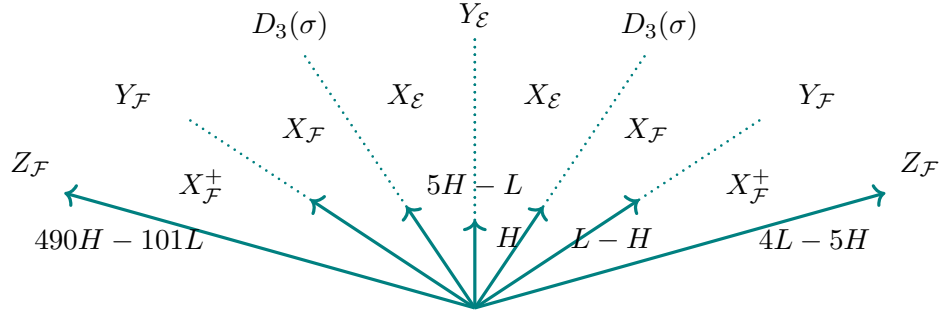
\begin{center}
\begin{tikzpicture}[xscale=0.9,yscale=0.9]
%\draw [help lines] (-5,0) grid (5,5);

\draw[->, very thick, teal] (0,0)--(0,1.3);
\draw[->, very thick, teal] (0,0)--(1,1.5);
\draw[->, very thick, teal] (0,0)--(2.4,1.6);
\draw[->, very thick, teal] (0,0)--(6,1.7);
\draw[->, very thick, teal] (0,0)--(-1,1.5);
\draw[->, very thick, teal] (0,0)--(-2.4,1.6);
\draw[->, very thick, teal] (0,0)--(-6,1.7);

\node [above] at (.5,.75) {$H$};
\node [above] at (2,0.7) {$L - H$};
\node [right] at (4,1) {$4L-5H$};
\node [above] at (0,1.5) {$5H-L$};
\node [left] at (-3.8,1) {$490 H - 101 L$};

%\node [above] at (0,4) {$\bP^3$};
\node [above] at (0,4) {$Y_{\sE}$};
\node [above] at (-2.7,3.8) {$D_3(\sigma)$};
\node [above] at (-5,2.8) {$Y_{\sF}$};
\node [above] at (-4,1.4) {$X_\sF^+$};
\node [above] at (-6.5,1.8) {$Z_\sF$};%{\rm(II)}
\node [above] at (-2.5,2.3) {$X_\sF$};
\node [above] at (-1,2.8) {$X_\sE$};
\node [above] at (2.7,3.8) {$D_3(\sigma)$};
\node [above] at (5,2.8) {$Y_{\sF}$};
\node [above] at (4,1.4) {$X_\sF^+$};
\node [above] at (6.5,1.8) {$Z_\sF$};%{\rm(II)}
\node [above] at (2.5,2.3) {$X_\sF$};
\node [above] at (1,2.8) {$X_\sE$};

\draw[decorate sep={.2mm}{1mm},fill, teal] (0,1.3)--(0,4);
\draw[decorate sep={.2mm}{1mm},fill, teal] (-1,1.5)--(-2.5,3.75);
\draw[decorate sep={.2mm}{1mm},fill, teal] (-2.4,1.6)--(-4.2,2.8);
\draw[decorate sep={.2mm}{1mm},fill, teal] (0,1.3)--(0,4);
\draw[decorate sep={.2mm}{1mm},fill, teal] (1,1.5)--(2.5,3.75);
\draw[decorate sep={.2mm}{1mm},fill, teal] (2.4,1.6)--(4.2,2.8);
%\node at (0,2) {$X_{\sF}$};
%\node at (2,2) {$X_{\sE}$};
%\node at (-2,2) {$X_{\sF}^{+}$};
\end{tikzpicture}
\end{center}

\begin{proof} The statement \eqref{mainB1} is obvious. By Proposition \ref{negdiag}, $X_{\sF}$ is a smooth Calabi--Yau threefold with Picard number $\rho (X_\sF) = 2$ and the relative Picard number $\rho (X_{\sF} / Y_{\sF})$ is $1$. Then $X_{\sF} \to Y_{\sF}$ is primitive and \eqref{mainB2} follows from Lemma \ref{r4Fcontr}.

To determine the supporting divisor of $X_\sF\dra Y_{\sE}$, we can verify that the inequality \eqref{deterwallcond} holds as in the proof of Theorem \ref{mainthmrk2}. By Proposition \ref{deterwall}, we have $\chi_*L=-L_\sE+5H_\sE$ under the map $\chi:X_\sF\dra X_\sE$. In particular, a supporting divisor of $X_\sF\dra Y_{\sE}$ is given by $5H-L$ and \eqref{mainB3} follows from Lemma \ref{Bcover_contra} and Remark \ref{Bcover_ODPs}.

Suppose that $q_S:S\ra\bP^2$ is the blow-up of 10 distinct points on $\bP^2$ and $X_\sF\dra X_\sF^+$ is the flop of $X_\sF\ra Y_\sF$, then as $(L-2H) \cdot C=S \cdot C=K_S  \cdot C=-1$, it is easy to see that the proper transform $(L-2H)^+$ on $X_\sF^+$ is relatively  ample over $Y_\sF$. In particular, $X_\sF\dra X_\sF^+$ is defined by $(L-2H)+\lambda(L-H)$ for $\lambda\gg0$.

The threefold $X_\sF^+$ contains a surface $S^+\cong\bP^2$ and hence there is an extremal contraction $X_\sF^+\ra Z_\sF$ contracting $S^+$ to a $\frac{1}{3}(1,1,1)$ point in $Z_\sF$ (see Proposition \ref{divcontr}). Note that the natural projection $S\ra \bP^2$ factors through $S^+$ and the induced contraction map $f:S\ra S^+$ does not extend to $X_\sF\dra X_\sF^+.$ To find the supporting divisor of $X_\sF\dra Z_\sF$, we need a movable $\bQ$-divisor $A$ on $X_\sF$ to be negative over $Y_\sF$ so that $A^+$ is semiample and $A^+|_{S^+}\equiv 0$.  Say $$A=x(L-H)+y(L-2H)\equiv x(L-H)+yS$$ for some $x,y\in\bQ$. Since $(L-H)|_S\equiv f^*\cO_{\bP^2}(1)$ and $S^+|_{S^+}=K_{S^+}=\cO_{\bP^2}(-3)$, the condition $A^+|_{S^+}\equiv0$ implies that $x=3y$. If $A \coloneqq 4L-5H\sim3(L-H)+S$, then $A$ is mobile as $|L-H|$ is base point free and big. Note that $L-H$ is a pull-back of an ample and base point free divisor on $Y_\sF$ and hence so is $(L-H)^+$.  In particular, $A^+\sim 3(L-H)^++S^+$ is base point free from the exact sequence,
$$H^0(X_{\sF^+},A^+)\rightarrow H^0(S^+,A^+|_{S^+})\rightarrow H^1(X_{\sF^+},3(L-H)^+)=0,$$
where the last equality is the Kawamata--Viehweg vanishing as $K_{X_{\sF^+}}=0$. Hence the linear system $|4L-5H|$ does defines the map $X_\sF\dra Z_\sF$ and fulfills the description in 
\eqref{mainB4}.

Finally, we compute the boundaries of $\cMov (X)$. Since with respect to ordered bases $\{L_\sF,H_\sF\}$ and $\{L_\sE,H_\sE\}$, we have the matrix
\[
[\chi_*] =
\begin{bmatrix}
-1 & 0 \\
5 & 1
\end{bmatrix}
=
[(\chi^{-1})_*],
\]
the composition map $\psi\coloneqq\chi^{-1}\circ\iota\circ\chi:X_\sF\dra X_\sF$ has the matrix representation with respect to the ordered basis $\{L_\sF,H_\sF\}$ as
\[
[\psi_*] =
\begin{bmatrix*}[r]
-34 & -7 \\
165 & 34
\end{bmatrix*}
=
[(\psi^{-1})_*]
\]
Hence $X_\sF\dra X_\sF\ra D_3(\sigma)$ is defined by $(\psi^{-1})_*H= 34H - 7L$. Similarly, $X_\sF\dra X_\sF\dra Y_\sF$ is defined by $(\psi^{-1})_*(L-H)= 131H - 27L$, and $X_\sF\dra X_\sF\dra Z_\sF$ is defined by $(\psi^{-1})_*(4L-5H)= 490H - 101L$.
\end{proof}

%%%%%%%%%%%%%%%%%%%%
%%%%%%%%%%%%%%%%%%%%

\section{Birational Models and Movable Cones II} \label{movfansecII}

In this section, we will treat the remaining cases $(\bP^4,\cO (1)^5)$, $(\bP^4,T_{\bP^4})$, and $(\Gr (2, 4), \cO (1)^4)$. We will see that both boundary rays of the movable cone $\cMov (X_{\sF})$ in these cases are irrational.

As before, we use the same notation as in Section \ref{flopsec} and apply Proposition \ref{negdiag} to construct a smooth Calabi--Yau threefold $X_{\sF}$ with Picard number $2$. We continue to write $L$ and $H$ for $L_{\sF}|_{X_{\sF}}$ and $H_{\sF}|_{X_{\sF}}$ respectively.

\subsection{\texorpdfstring{$M = \bP^4$}{M=P4}} \label{P4IIsubsec}

We first remark that the case $(\sE, \sF) = (\cO^{4}, T_{\bP^4})$ can be regarded as a special case of $(\sE, \sF) = (\cO^{5}, \cO (1)^5)$. Indeed, from the Euler sequence there is a natural embedding \begin{equation} \label{lindual}
    \bP(T_M)\hookrightarrow\bP^4\times\bP^4\cong\bP^4\times(\bP^4)^\vee,
\end{equation}
where $(\bP^4)^\vee$ is the dual projective space. Here we can view $\bP(T_M)$ as the incidence variety and will  only consider the case $(\sE, \sF) = (\cO^{5}, \cO (1)^5)$.

\begin{theorem}[\cite{Fryers01}] \label{mainP4quin}
Let $(M, \sE, \sF) = (\bP^4, \cO^{5}, \cO (1)^5)$. Then for a general morphism $\sigma : \sE^{\vee} \to \sF$, $X_{\sF}$ is a smooth Calabi--Yau threefold with Picard number $2$ with
$$
\Nef(X_\sF)=\bR_{\geqslant 0}[L-H]+\bR_{\geqslant 0}[H],
$$
such that
\begin{enumerate}[(i)]
    \item\label{mainP4quin1} the determinantal contraction $\pi_{\sF}$ is induced by $|H|$ and $D_4 (\sigma)$ is a Calabi--Yau threefold with $50$ ODPs;

    \item\label{mainP4quin2} $|L -  H|$ (resp.~$|5 H - L|$) induces a primitive contraction $X_{\sF} \to Y_{\sF}$ (resp.~$X_{\sE} \to Y_{\sE}$) of type I and $Y_{\sF}$ (resp.~$Y_{\sE}$) is a Calabi--Yau threefold with $50$ ODPs;

    \item\label{mainP4quin3} the flop $X_{\sF}^{+}$ of $X_{\sF} \to Y_{\sF}$ admits a primitive contraction of type I induced by $|4 L -5 H|$ and $X^{+}_{\sF}$ is isomorphic to the flop of $X_{\sE}$, which we denote by $X^+$.
\end{enumerate}
Moreover, the movable cone of $X_{\sF}$ is given by
\begin{equation}\label{P45_movcone}
    \cMov(X_{\sF}) = \bR_{\geqslant 0} [- L + (3 + \sqrt{3}) H] + \bR_{\geqslant 0} [L + (- 3 + \sqrt{3}) H]
\end{equation}
which is covered by the nef cones of $X_{\sF}$, $X_{\sE}$ and $X^+$, and there are no more minimal models of $X_{\sF}$.
\end{theorem}

The picture of $\Mov(X_\sF)$ is the following. The rays accumulate to the boundary rays of slopes $- 3 - \sqrt{3}$ and $- 3 + \sqrt{3}$.

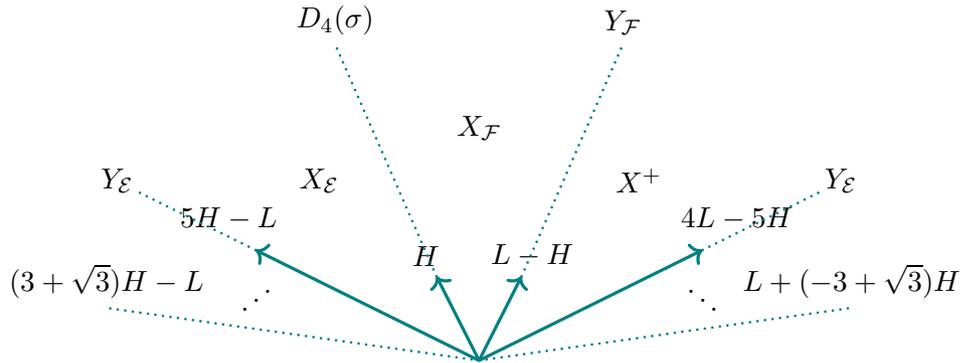
\begin{center}
\begin{tikzpicture}[xscale=0.7,yscale=0.7]
%\draw [help lines] (-5,0) grid (5,5);

\draw[->, very thick, teal] (0,0)--(0.8,1.6);
\draw[->, very thick, teal] (0,0)--(-0.8,1.6);
\draw[->, very thick, teal] (0,0)--(4.2,2.1);
\draw[->, very thick, teal] (0,0)--(-4.2,2.1);
%\draw[->, very thick, teal] (0,0)--(5.6,0.8);
%\draw[->, very thick, teal] (0,0)--(-5.6,0.8);

\node [above] at (1,1.6) {$L-H$};
\node [above] at (-1,1.6) {$H$};
\node [left] at (-3.6,2.7) {$5 H - L$};
\node [right] at (3.6,2.7) {$4 L - 5 H$};

\node [above] at (0,4) {$X_\sF$};
\node [above] at (3,3) {$X^+$};
\node [above] at (-3,3) {$X_{\sE}$};
\node [above] at (2.7,6) {$Y_{\sF}$};
\node [above] at (-2.7,6) {$D_4 (\sigma)$};
\node [above] at (6.8,3) {$Y_{\sE}$};
\node [above] at (-6.8,3) {$Y_{\sE}$};
\node [above] at (4.2,0.7) {$\ddots$};
\node [above] at (-4.2,0.7) {$\iddots$};
\node [above] at (-7,1) {$(3 + \sqrt{3})H - L$};
\node [above] at (7,1) {$L + (- 3 + \sqrt{3})H$};

\draw[decorate sep={.2mm}{1mm},fill, teal] (-0.8,1.6)--(-2.7,6);
\draw[decorate sep={.2mm}{1mm},fill, teal] (0.8,1.6)--(2.7,6);
\draw[decorate sep={.2mm}{1mm},fill, teal] (4.2,2.1)--(6.4,3.2);
\draw[decorate sep={.2mm}{1mm},fill, teal] (-4.2,2.1)--(-6.4,3.2);
%\draw[decorate sep={.2mm}{1mm},fill, teal] (5.6,0.8)--(7,1);
%\draw[decorate sep={.2mm}{1mm},fill, teal] (-5.6,0.8)--(-7,1);
\draw[decorate sep={.2mm}{1mm},fill, teal] (0,0)--(7,1);
\draw[decorate sep={.2mm}{1mm},fill, teal] (0,0)--(-7,1);
\end{tikzpicture}
\end{center}

\begin{proof}
Since the result is known by \cite[Lemma 1]{Fryers01} and \cite{Borcea91}, we only give a rough sketch in our notation. For a general morphism $\sigma : \sE^{\vee} \to \sF$, we have the commutative diagram \eqref{flopdiag} and the birational morphism $\chi : X_{\sF} \dashrightarrow X_{\sE}$. As in the proof of Theorem \ref{mainthmrk2}, we can verify that the inequality \eqref{deterwallcond} holds and thus $\chi_{\ast} L = - L_\sE + 5 H_\sE$.

To construct the flops, we observe that the morphism $\sigma$ corresponds to a $5 \times 5$ matrix $M (z)$ of linear forms
$$
M_{i j} (z) = \sum_{k} a_{i j k} z_k,
$$
and $D_4 (\sigma) = \{z \in \bP^4 \mid \det M (z) = 0\}$. Since $\bP (\sF) \cong \bP^4 \times \bP^4$, we can view $X_{\sF}$ as the variety
$$
\{(z, z') \in \bP^4 \times \bP^4 \mid M (z) \cdot [z']^t = 0\}.
$$
On the other hand, $X_{\sE}$ is induced by the dual morphism $\sigma^{\vee}$, which is defined by the transpose of $M(z)$ in $\bP (\sE) = \bP^4 \times \bP^4 \ni (z, z'')$.

We construct the other matrices of linear forms
\begin{center}
    $M'_{ij} (z') = \sum_k a_{i k j} z_k'$ and $M''_{ij} (z'') = \sum_k a_{k i j} z_k''$
\end{center}
such that
\begin{equation}\label{P4_FE}
    M (z) \cdot [z']^t = M' (z') \cdot [z]^t \text{ and }
    M (z)^t \cdot [z'']^t = M'' (z'') \cdot [z]^t.
\end{equation}
Hence the second projection gives rise to a small contraction  $X_{\sF} \to Y_{\sF}$ (resp.~$X_{\sE} \to Y_{\sE}$) where $Y_{\sF}$ (resp.~$Y_{\sE}$) is the zero locus of $\det M' (z')$ (resp.~$\det M'' (z'')$) and the supporting divisor is $L - H$ (resp.~$L_{\sE}$). By Proposition \ref{numbODP}, the number of ODPs in $D_4 (\sigma)$ is
$$
\int_{\bP^4} c_2 (\cO (1)^5)^2 - c_1 (\cO (1)^5) c_3 (\cO (1)^5) = 50,
$$
and similarly for $Y_{\sF}$ and $Y_{\sE}$.

According to \eqref{P4_FE}, the flop $X_{\sF}^+$ (resp.~$X_{\sE}^+$) of $X_{\sF} \to Y_{\sF}$ (resp.~$X_{\sE} \to Y_{\sE}$) is defined by $M' (z')^t \cdot [z'']^t = 0$ (resp.~$M'' (z'')^t \cdot [z']^t = 0$). More precisely, the matrix $M' (z')^t$ defines a morphism $\cO^5 \to \cO (1)^5$ over $\bP^4 \owns z'$ such that $X_{\sF}^+$ is zero locus of the section defined by this morphism (cf.~\eqref{Cayleytrick}), and similarly for $M'' (z'')^t$. In particular, $X_{\sF}^+$ and $X_{\sE}^+$ are isomorphic, denoted by $X^+$, because
$M' (z')^t \cdot [z'']^t = M'' (z'')^t \cdot [z']^t$. Thus we have the following diagram:
\begin{center}
      \begin{tikzcd}
      X^+ \ar[rd] \ar[rr, dashed, "\iota_{\sF}"] & & X_{\sF}  \ar[ld] \ar[rd] \ar[rr, dashed,"\chi"] & & X_{\sE}  \ar[rd] \ar[ld] & & X^+ \ar[ll, dashed, "\iota_{\sE}", swap] \ar[ld] \\
       & Y_{\sF} & & D_4(\sigma) & &  Y_{\sE} &.
      \end{tikzcd}
\end{center}

Set $\varrho = \iota_{\sF}\circ \iota_{\sE}^{- 1}\circ \chi$. Applying Proposition \ref{deterwall} to morphisms induced by $M (z)$, $M' (z')^t$, and $M'' (z'')^t$, we infer that the matrix representation with respect to the ordered basis $\{L, H\}$ is given by
$$
[\varrho^{\ast}] =
\begin{bmatrix*}[r]
-19 & -15 \\
90 & 71
\end{bmatrix*}
$$
and $L + (- 3 + \sqrt{3}) H$ (resp.~$- L + (3 + \sqrt{3}) H$) is an eigenvector of $\varrho^{\ast}$, corresponding to the eigenvalue $26 - 15 \sqrt{3}$ (resp. $26 + 15 \sqrt{3}$) of $\varrho^{\ast}$. In particular, the birational map $\varrho : X_{\sF} \dashrightarrow X_{\sF}$ is of infinite order.
\end{proof}

\begin{remark}
By linear duality (cf.~\eqref{lindual}), \cite[Section 3.2]{HT18} has also calculated the action of the birational map $\varrho$ and the movable cone (if we put $X_1 = X_{\sF}$, $H_1 = L - H$ and $H_2 = H$ to adapt our notation with the one in \cite{HT18}).
\end{remark}

\subsection{\texorpdfstring{$M = \Gr (2, 4)$}{M=Gr(2,4)}}

We are going to treat the case $\sE=\cO^{4}$ and $\sF=\cO(1)^{4}$. Consider the natural projection $X_{\sE} \to \bP^3$, which is the restriction of $\bP (\sE) = \Gr (2, 4) \times \bP^3 \to \bP^3$ to $X_{\sE}$.

\begin{lemma} \label{Gr4_contra}
Let $X_{\sE} \to Y_{\sE} \to \bP^3$ be the Stein factorization. Then the natural projection $X_{\sE} \to \bP^3$ is generically $2:1$ and $X_{\sE} \to Y_{\sE}$ is a small contraction if $X_{\sE}$ is chosen in general, and similarly for $X_{\sF}$.
\end{lemma}

\begin{proof}
Set $\Gr =  \Gr (2, 4)$. Note that $\bP(\sE) = \Gr \times \bP^3$ is defined by a global section of $\cO_{\bP^5} (2) \boxtimes \cO_{\bP^3}$ under the Pl\"ucker embedding of $\Gr$. We can view 
\[
    \cO_{\sE} (1) = (\cO_{\bP^5}\boxtimes \cO_{\bP^3} (1) )|_{\bP(\sE)}    
\]
and $X_{\sE}$ is defined by four general global sections of $\cO_{\Gr} (1) \boxtimes \cO_{\bP^3} (1)$. Hence $X_{\sE}$ is a complete intersection of type $(2,0),4\times(1,1)$ in $\bP^5\times\bP^3$.

We are going to show that $\pi_{\sE} : X_{\sE} \to Y_{\sE}$ is small. Observe that the fiber of $X_{\sE} \to \bP^3$ over a point $P \in \bP^3$ is determined by the system
\[
L_1(P) = \cdots = L_4(P) = Q = 0,
\]
where $L_i=\sum_{j=0}^5l_{ij}z_j$ are of type $(1,1)$ with coefficients being linear forms $l_{ij}$ on $\bP^3$ and $Q\in H^0(\cO_{\bP^5}(2))$. Note that $\Gr$ is the zero locus of $Q$.

Let $V \subseteq \bP^5 \times \bP^3$ be the complete intersection fourfold defined by $L_1, \cdots, L_4$, i.e., it is defined by a general section of $H^0 (\cO_{\bP^5} (1) \boxtimes \cO_{\bP^3} (1)^4)$. The section corresponds to a general morphism 
\[
    \tau \colon \cO_{\bP^3}^{4} (- 1) \to \cO_{\bP^3}^{ 6}
\]
defined by the matrix $\left[l_{ij} \right]^t$. Let $q_V$ be the restriction of the projection $\bP ^5 \times \bP^3 \to \bP^3$ to $V$. For each $P \in D_k (\tau) \setminus D_{k - 1} (\tau)$, the fiber $q_V^{- 1} (P)$ is $\bP (\coker \tau (P)) \cong \bP^{5 - k}$ for $0 \leqslant k \leqslant 4$. Note that $D_3 (\tau) $ and $D_2 (\tau)$ have the expected codimension $(4-3) \times (6-3) = 3$ and $(4-2) \times (6-2) = 8$ respectively. If $X_{\sE}$ is chosen in general, then $D_2 (\tau) = \varnothing$ and $D_3(\tau)$ consists of (smooth) finitely many points. %whose number can be determined by Giambelli--Thom--Porteous formula \cite[Theorem 14.4(a)]{fulton98},
%\[
%[D_3(\tau)] = c_3 (\cO_{\bP^3}(1)^{ 6}) \cap [\bP^3] \in A_3(\bP^3).
%\]
%We find that the number $|D_3 (\tau)|$ is $20$. 

Now we have that a fiber of $q_V \colon V \to \bP^3$ is $\bP^2$ (resp.~$\bP^1$) if $P \in D_3 (\tau)$ (resp.~$P \in \bP^3 \setminus D_3 (\tau)$). By the fact that $X_{\sE} = V \cap (Q=0)$, the contracting locus of the double cover $X_{\sE} \to \bP^3$ has dimension at most one. Since $X_{\sE} \xrightarrow{\pi_{\sE}} Y_{\sE} \to \bP^3$ is the Stein factorization of $X_{\sE} \to \bP^3$, we conclude that $\pi_{\sE}$ is small. %and the number of singularities of $Y_{\sE}$ is $20$.

After tensoring $\sF$ with $\cO (- 1)$, the same conclusion holds for $X_{\sF}$.
\end{proof}

There are two involutions over $\bP^3$ induced from the natural projections to $\bP^3$ on $X_\sF$ and $X_\sE$. We denote the involutions by $\iota_\sF:X_\sF\dra X_\sF$ and $\iota_\sE:X_\sE\dra X_\sE$, which fit into the following diagram.

$$
\SelectTips{eu}{12}
\xymatrix@!0{
 & X_{\sF} \ar[ld] \ar@{-->}[rr]^{\iota_{\sF}} & & X_{\sF} \ar@{-->}[rr]^{\chi} \ar[ld] \ar[ldd] \ar[rd] & & X_{\sE} \ar@{-->}[rr]^{\iota_{\sE}} \ar[ld] \ar[rd] \ar[rdd] & & X_{\sE} \ar[rd] & \\
Y_{\sF} \ar[rr] \ar[rrd]_{2 : 1} & & Y_{\sF} \ar[d] & & D_3 (\sigma) & & Y_{\sE} \ar[rr] \ar[d] & & Y_{\sE} \ar[lld]^{2 : 1} \\
&& \bP^3 &&&& \bP^3 &&
}
$$

Our aim is to compute the proper transforms of divisors under these involutions. Applying Proposition \ref{LHinter}, we get Table \ref{table_Gr41}.
%Applying Proposition \ref{LHinter}, we summarize the intersection numbers on $X_{\sF}$ in the following table.
\begin{table}[H]
  \centering
  \begin{tabular}{ccccccc}
    \toprule
    $L^3$ & $L^2 \cdot H$ & $L \cdot H^2$ & $H^3$ & $L \cdot c_2(T_Z)$ & $H \cdot c_2(T_Z)$ & $\#$ of ODPs \\
    \midrule
     70 & 40  & 20 & 8 & 100 & 56 & 40 \\
    \bottomrule
  \end{tabular}
  \caption{The intersection numbers on $X_{\sF}$.}
  \label{table_Gr41}
\end{table}

\begin{lemma} \label{Gr4invol}
For the involutions $\iota_\sF:X_\sF\dra X_\sF$ and $\iota_\sE:X_\sE\dra X_\sE$, the matrix representation with respect to $\{L_\bullet,H_\bullet\}$ ($\bullet = \sF$ or $\sE$) is given by
$$[(\iota_\sF)_*] =
\begin{bmatrix}
9 & 8 \\
-10 & -9
\end{bmatrix}
=
[(\iota^{-1}_\sF)_*]
\ {\rm and}\ [(\iota_\sE)_*] =
\begin{bmatrix}
1 & 8 \\
0 & -1
\end{bmatrix}
=
[(\iota^{-1}_\sE)_*].$$
\end{lemma}

\begin{proof} Let $\iota=\iota_\sF$. Note that $\iota_*(L-H)=L-H$ and write $\iota_*H=xL+yH$. Since $\iota\in\Bir(X_\sF)$ is small, we have
$$\begin{cases}
(L-H) \cdot H^2=(L-H) \cdot (xL+yH)^2\\
(L-H)^2 \cdot H=(L-H)^2 \cdot (xL+yH)\\
\end{cases}\ {\rm or}\
\begin{cases}
12=30x^2+40xy+12y^2\\
8=10x+8y\\
\end{cases}$$
Since $(x,y)=(0,1)$ is impossible, the only solution is $(x,y)=(8,-9)$ and the rest is clear.

Note that $L_\sE^3 = 2$, $L_\sE \cdot H_\sE^2 = 12$, $L_\sE^2 \cdot H_\sE = H_\sE^3 = 8$ and $(\iota_{\sE})_{\ast} L_{\sE} = L_{\sE}$. The proof for $\iota_\sE$ is the same as above, and is left to the reader.
\end{proof}

\begin{theorem}\label{mainGr41}
Let $(M, \sE, \sF) = (\Gr (2, 4), \cO^{4}, \cO (1)^4)$. Then for a general morphism $\sigma : \sE^{\vee} \to \sF$, $X_{\sF}$ is a smooth Calabi--Yau threefold of  Picard number $2$ with
$$
\Nef(X_\sF)=\bR_{\geqslant 0}[L-H]+\bR_{\geqslant 0}[H],
$$
such that
\begin{enumerate}[(i)]
    \item\label{mainGr411} the determinantal contraction $\pi_{\sF}$ is induced by $|H|$;

    \item\label{mainGr412} $|L -  H|$ induces a primitive contraction $X_{\sF} \to Y_{\sF}$ of type I, and the double cover $X_{\sF} \to \bP^3$ factors through $Y_{\sF}$;

    \item\label{mainGr413} $|4 H - L|$ induces a primitive contraction $X_{\sE} \to Y_{\sE}$ of type I, and the double cover $X_{\sE} \to \bP^3$ factors through $Y_{\sE}$.
\end{enumerate}
Moreover, the movable cone of $X_{\sF}$ is given by
\begin{equation}\label{Gr4_movcone}
    \cMov(X_{\sF}) = \bR_{\geqslant 0} [- 4 L + (1 0 + \sqrt{3 0})H] + \bR_{\geqslant 0} [4 L + (- 1 0 + \sqrt{3 0}) H]
\end{equation}
which is covered by nef cones of $X_{\sF}$ and $X_{\sE}$, and there are no more minimal models of $X_{\sF}$.
\end{theorem}

The picture of $\Mov(X_\sF)$ is the following. The rays accumulate to the boundary rays of slopes $(- 10  - \sqrt{30}) / 4$ and $(- 10 + \sqrt{30}) / 4$.

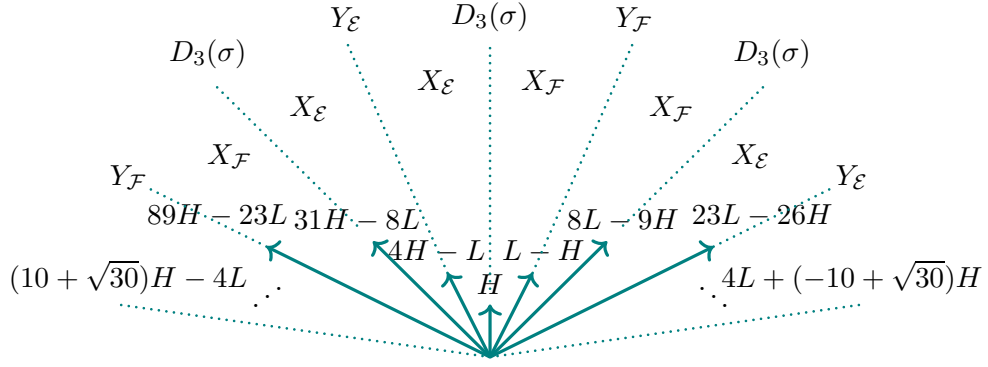
\begin{center}
\begin{tikzpicture}[xscale=0.7,yscale=0.7]
%\draw [help lines] (-5,0) grid (5,5);

\draw[->, very thick, teal] (0,0)--(0,1);
\draw[->, very thick, teal] (0,0)--(0.8,1.6);
\draw[->, very thick, teal] (0,0)--(-0.8,1.6);
\draw[->, very thick, teal] (0,0)--(2.2,2.2);
\draw[->, very thick, teal] (0,0)--(-2.2,2.2);
\draw[->, very thick, teal] (0,0)--(4.2,2.1);
\draw[->, very thick, teal] (0,0)--(-4.2,2.1);
%\draw[->, very thick, teal] (0,0)--(5.6,0.8);
%\draw[->, very thick, teal] (0,0)--(-5.6,0.8);

\node [above] at (0,1) {$H$};
\node [above] at (1,1.6) {$L-H$};
\node [above] at (-1,1.6) {$4H-L$};
\node [above] at (2.5,2.2) {$8L -9H$};
\node [above] at (-2.5,2.2) {$31H- 8 L$};
\node [left] at (-3.6,2.7) {$89H-23L$};
\node [right] at (3.6,2.7) {$23L-26H$};

\node [above] at (0,6) {$D_3(\sigma)$};
\node [above] at (1,4.8) {$X_\sF$};
\node [above] at (-1,4.8) {$X_\sE$};
\node [above] at (3.4,4.3) {$X_\sF$};
\node [above] at (-3.4,4.3) {$X_\sE$};
\node [above] at (2.7,6) {$Y_{\sF}$};
\node [above] at (-2.7,6) {$Y_{\sE}$};
\node [above] at (-4.9,3.4) {$X_\sF$};
\node [above] at (4.9,3.4) {$X_\sE$};
\node [above] at (6.8,3) {$Y_{\sE}$};
\node [above] at (-6.8,3) {$Y_{\sF}$};
\node [above] at (5.3,5.3) {$D_3(\sigma)$};
\node [above] at (-5.3,5.3) {$D_3(\sigma)$};
\node [above] at (4.2,0.7) {$\ddots$};
\node [above] at (-4.2,0.7) {$\iddots$};
\node [above] at (-6.8,1) {$(10+\sqrt{3 0})H-4L$};
\node [above] at (6.8,1) {$4L+(-10+\sqrt{3 0})H$};
%\node [above] at (-7,1) {$\rho_{H/L}\approx0.1.12$};
%\node [above] at (7,1) {$\rho_{L/H}\approx0.258$};

\draw[decorate sep={.2mm}{1mm},fill, teal] (0,1.3)--(0,6);
\draw[decorate sep={.2mm}{1mm},fill, teal] (-0.8,1.6)--(-2.7,6);
\draw[decorate sep={.2mm}{1mm},fill, teal] (0.8,1.6)--(2.7,6);
\draw[decorate sep={.2mm}{1mm},fill, teal] (2.5,2.5)--(5.2,5.2);
\draw[decorate sep={.2mm}{1mm},fill, teal] (-2.5,2.5)--(-5.2,5.2);
\draw[decorate sep={.2mm}{1mm},fill, teal] (4.2,2.1)--(6.4,3.2);
\draw[decorate sep={.2mm}{1mm},fill, teal] (-4.2,2.1)--(-6.4,3.2);
\draw[decorate sep={.2mm}{1mm},fill, teal] (0,0)--(7,1);
\draw[decorate sep={.2mm}{1mm},fill, teal] (0,0)--(-7,1);
%\draw[decorate sep={.2mm}{1mm},fill, teal] (5.6,0.8)--(7,1);
%\draw[decorate sep={.2mm}{1mm},fill, teal] (-5.6,0.8)--(-7,1);
\end{tikzpicture}
\end{center}

\begin{proof} The statement \eqref{mainGr411} is obvious. Statements  \eqref{mainGr412} and \eqref{mainGr413} follow from Lemma \ref{Gr4_contra}.

Under the map $\chi:X_\sF\dra X_\sE$, we have $\chi_*L=-L_\sE+4H_\sE$. Indeed, we can verify that the inequality \eqref{deterwallcond} holds as in the proof of Theorem \ref{mainthmrk2}. In particular, a supporting divisor of $X_\sF\dra Y_{\sE}$ is given by $4H-L$ and with respect to ordered bases $\{L_\sF,H_\sF\}$ and $\{L_\sE,H_\sE\}$, we have the matrix
$$
[\chi_*] =
\begin{bmatrix*}[c]
-1 & 0 \\
4 & 1
\end{bmatrix*}
=
[(\chi^{-1})_*].
$$

Now, we are going to find the boundary of the movable cone $\Mov(X_\sF)$. Set $\theta \coloneqq \chi^{-1}\circ \iota_{\sE}\circ \chi \colon X_{\sF} \dashrightarrow X_{\sF}$. From Lemma \ref{Gr4invol} and above, we see that the set
$$
\Nef(X_{\sF}) \bigcup \chi^{\ast} \Nef (X_{\sE}) \bigcup (\iota_{\sE} \chi)^{\ast} \Nef (X_{\sE}) \bigcup \theta^{\ast} \Nef(X_{\sF})
$$
is given by the rational polyhedral cone
$$
\Pi \coloneqq \bR_{\geqslant 0 } [8 9 H - 2 3 L] + \bR_{\geqslant 0 } [L - H].
$$

Write $\varrho \coloneqq \iota_{\sF}\circ \theta$. From concrete calculations in $2 \times 2$ matrices, we get that, with respect to the ordered basis $\{ L_{\sF}, H_{\sF} \}$,
$$
[\varrho^{\ast}]
=
[(\chi^{- 1})_*][(\iota_\sE^{- 1})_*][\chi_*][(\iota_\sF^{- 1})_*]
=
\begin{bmatrix*}[r]
-199 &-176  \\
770 & 681
\end{bmatrix*}
$$
and $- 4 L + (1 0 + \sqrt{3 0})H$ (resp.~$4 L + (- 1 0 + \sqrt{3 0})H$) is an eigenvector of $\varrho^{\ast}$, corresponding to the eigenvalue $241 + 44 \sqrt{30} > 1$ (resp.~$241 - 44 \sqrt{30} = 1 / (241 + 44 \sqrt{30})$) of $\varrho^{\ast}$. In particular, $\varrho^{\ast}$ is of infinite order.

From the actions of $(\varrho^{\pm})^{\ast}$, we see that $(\varrho^{\pm})^{\ast} \Pi$ and $\Pi$ have non-overlapping interior and $(\varrho^{\pm})^{\ast} \Pi \cup \Pi$ is a cone. Let us denote by $C$ the union
\[
    C = \bigcup_{n \in \bZ} (\varrho^{n})^{\ast} \Pi,
\]
which is a cone and $\overline{C} \subseteq \cMov (X_{\sF{}})$. Since $\Pi$ is contained in the r.h.s.~of \eqref{Gr4_movcone} which is spanned by eigenvectors of $\varrho^{\ast}$, we find that the closure $\overline{C}$ coincides with the r.h.s.~of \eqref{Gr4_movcone}. On the other hand, let $d$ be a rational point of the interior of $\cMov (X_{\sF})$. There is an integer $m > 0$ and an effective movable divisor $D$ such that $m d = [D]$. If $D$ is nef, then $d \in \Nef(X_{\sF}) \subseteq \Pi$. If $D$ is not nef, we can run the log minimal model program for the klt pair $(X_{\sF}, \varepsilon D)$, $0 < \varepsilon \ll 1$, to find a birational map $f$ such that $f_{\ast} D$ is nef. Note that any birational map between minimal models is decomposed into finitely many flops \cite{Kawamata08}. By the shapes of $\Nef(X_{\sF})$ and $\Nef(X_{\sE})$, the birationa map $f$ must be either $X_{\sF} \dashrightarrow X_{\sF}$ or $X_{\sF} \dashrightarrow X_{\sE}$.

If $f \in \Bir (X_{\sF})$, we claim that $f = \varrho^n \circ \iota_{\sF}$ for some $n \in \bZ$ (up to automorphisms of $X_{\sF}$). Recall that any flopping contraction of a Calabi--Yau manifold is given by a codimension one face of the nef cone. Then the claim follows from the shapes of $\Nef(X_{\sF})$ and $\Nef(X_{\sE})$ and fact that $\theta = \iota_{\sF} \circ \varrho$ and $\varrho \circ \iota_{\sF} = \iota_{\sF} \circ \varrho^{- 1}$. Hence we get that $D \in (\varrho^{n - 1})^{\ast} \Pi$ by $f_{\ast} D \in \Nef (X_{\sF})$ and the fact that $(\rho^{- 1})^{\ast} \Pi = (\iota_{\sF})^{\ast} \Pi$. To treat the case $f \colon X_{\sF} \dashrightarrow X_{\sE}$, we may assume that up to birational automorphisms of $X_{\sF}$ the birational map $f$ is either $ \chi$ or $\iota_{\sE} \circ \chi$. By the definition of $\Pi$, we find that $D \in \Pi$. Therefore $d \in C$ in any case, and hence we get $\cMov (X_{\sF{}}) \subseteq \overline{C}$, which completes the proof.
\end{proof}

%By the chamber structure of $\Pi$, there is an integer $n \in \bZ$ such that $(\varrho^{\ast})^n [D] \in \Pi$, and thus $(\varrho^{\ast})^n d \in \Pi$. Hence we have $\cMov (X_{\sF{}}) \subseteq \overline{C}$, which completes the proof.

%Since both boundary rays of $\overline{C}$ is spanned by eigenvectors of $\varrho^{\ast}$, we see that the closure $\overline{C}$ coincides with the r.h.s.~of \eqref{Gr4_movcone} and thus $\overline{C} \subseteq \cMov (X_{\sF{}})$. 

%\cup  (\varrho^{\ast})^{- 1} \Pi

%Then we see that the closure $\overline{C}$ coincides with the r.h.s.~of \eqref{Gr4_movcone}, which is clear from linear algebra. 

%The same argument in \cite[Lemma 6.7]{Og-CY2} shows that $\cMov (X_{\sF{}}) = \overline{C}$.

\begin{remark}
%From the actions of $\varrho^{\ast}$, we see that$\varrho^{\ast} \Nef(X_{\sF})$ and $\Nef(X_{\sF})$ have non-overlapping interior. Thus $\varrho \notin \Aut (X_{\sF})$. Using the same argument in \cite[Lemma 6.4]{Og-CY2}, we get that the group of birational maps of $X_{\sF}$ is given by $\Bir (X_{\sF}) = \Aut (X_{\sF}) \cdot \left< \varrho, \iota_{\sF} \right>$.
By the above argument in the proof of Theorem \ref{mainGr41}, we see that the group of birational maps of $X_{\sF}$ is given by $\Bir (X_{\sF}) = \Aut (X_{\sF}) \cdot \left< \varrho, \iota_{\sF} \right>$.
\end{remark}
%%%%%%%%%%%%%%%%%%%%
%%%%%%%%%%%%%%%%%%%%

\appendix

\renewcommand{\thesection}{\Alph{section}}
\section{Chern Classes of Virtual Quotient Bundles} \label{chernclsubsec}

For the convenience of the reader, we collect some formulas of Chern and Segre classes that we need (cf.~\cite[Example 3.2.7 (a)]{fulton98}). For bundles $\sA$ and $\sB$, we  write $\sB^{\vee}$ for the dual bundle of $\sB$,
$$
c (\sA - \sB^{\vee}) = c (\sA) / c (\sB^{\vee}) = c (\sA) s (\sB^{\vee}),
$$
and let $c_k (\sA - \sB^{\vee})$ be the $k$th term in this expansion, that is,
$$
c_k (\sA - \sB^{\vee}) = \sum_{i = 0}^k c_i (\sA) s_{k - i} (\sB^{\vee}).
$$
By the definition of Chern and Segre classes, we get $s_1 (\sB^{\vee}) = c_1(\sB)$ and
\begin{align*}
    s_2 (\sB^{\vee}) &= c_1(\sB)^2 - c_2(\sB), \\
    s_3 (\sB^{\vee}) &= c_1(\sB)^3 - 2 c_1(\sB) c_2(\sB) + c_3(\sB), \\
    s_4 (\sB^{\vee}) &=  c_1(\sB)^4 - 3 c_1(\sB)^2 c_2(\sB) + 2 c_1(\sB) c_3(\sB) + c_2(\sB)^2 - c_4(\sB).
\end{align*}

\section{Hodge Numbers}\label{hodgenumsubsec}

The aim of this section is to compute the Hodge numbers of the  smooth Calabi--Yau threefolds $X$ obtained in Sections  \ref{movfansecI} and \ref{movfansecII}. This can be done by using Koszul complexes or the following known result for $\chi_{top} (X) = 2(h^{1, 1}(X) - h^{2, 1} (X))$ (see, for example, \cite[Proposition 2.3]{SSW2016} and \cite[Example 3.8]{NS95}).

\begin{proposition} \label{compEuler}
Let $X \to Y$ be a small resolution of Calabi--Yau threefold $Y$. If $Y$ is smoothable to a smooth Calabi--Yau threefold $\widetilde{Y}$, then $Y$ has only ODPs if and only if
$$
\chi_{top}(X) - \chi_{top}(\widetilde{Y}) = 2 |\Sing (Y)|.
$$
Moreover, the Hodge numbers are given by $h^{1, 1} (X) = \rho (X)$ and
$$
h^{2, 1} (\widetilde{Y}) - h^{2, 1} (X) = |\Sing (Y)| - \rho(X / Y).
$$
\end{proposition}

In our situation, $\tY \in |- K_M|$ is a smooth hypersurface in a smooth Fano fourfold $M$ with $\rho (M) = 1$ and the relative Picard number $\rho(X / Y)$ is $1$. From the Lefschetz hyperplane theorem, we get $h^{1, 1} (\tY) = 1$. Hence to find $\chi_{top} (X)$ (or, equivalently, $h^{2,1}(X)$), it is enough to compute $\chi_{top}(\tY)$. Note the numbers of ODPs are given in Proposition \ref{numbODP}.

\begin{lemma}\label{chitY}
With notation as above, we have
\[
  \chi_{top} (\tY) = \int_M c_1 (T_M)  c_3 (T_M) - c_1 (T_M)^2  c_2 (T_M).
\]
\end{lemma}

\begin{proof}
By the fact that $N_{\tY / M} \cong \cO ( - K_M)$ and
\[
0 \to T_{\tY} \to T_{M} |_{\tY} \to N_{\tY / M} \to 0,
\]
we have $[\tY] = c_1 (T_M) \cap [M]$ in $A_3 (M)$ and
\begin{equation*}
   c_3(T_{\tY}) = (c_3 (T_M) - c_1 (T_M) . c_2 (T_M))|_{\tY}.
\end{equation*}
The lemma follows from the Gauss--Bonnet theorem $\chi_{top} (\tY) = \int_{\tY} c_3(T_{\tY})$.
\end{proof}

According to the above lemma, our problem reduces to computing the Chern classes of the tangent bundles $T_M$. To shorten notation, we use $c_1^{k}\cdot c_{4 - k}$  for  $\int_M c_1 (T_M)^k  c_{4 - k} (T_M)$.

When $M = \bP^4$ or $\Gr (2, 4)$, we find that $(c_1 \cdot c_3,c_1^2\cdot c_2)=(50,250)$ and $(48,224)$ respectively. Therefore we infer that $(\chi_{top} (\tY),h^{2, 1} (\tY)) =( -200,101)$ and $(-176,89)$ respectively.

For the remaining cases, from Lemmas \ref{c2Fano} and \ref{chitY} it follows that
\begin{equation*}\label{chitYdPMu}
    \chi_{top} (\tY) = c_1 \cdot c_3 -
    \begin{cases}
    18 d + 108  & \text{if $M$ is del Pezzo},\\
    4 d + 96  & \text{if $M$ is Mukai}.
    \end{cases}
\end{equation*}
To compute $c_1\cdot c_3$, let us recall the classification of smooth del Pezzo and Mukai fourfolds of Picard number $1$ (see \cite[Theorem 3.3.1 and Section 5.2]{bookAGV} and the references given there).

\begin{notation}
%For a projective variety $G \subseteq \bP^N$, a complete intersection $G \cap \Lambda \subseteq \Lambda$ by a  ($k$-codimensional) linear subspace $\Lambda$ of $\bP^N$ is called a \emph{($k$-codimensional) linear section} of $G \subseteq \bP^N$. Note that $\Lambda = \cap_{i = 1}^k H_i$ is the transverse intersection of $k$ hyperplanes $H_i$ in $\bP^N$. 
We will use the symbol $M_{d_1, d_2, \cdots,d_k}$ to denote a general complete intersection hypersurfaces of indicated degrees in a given polarized  variety.
\end{notation}

We state the classification of smooth del Pezzo fourfolds, classified by Fujita (cf.~\cite{Fujita82} and \cite[Theorem 3.3.1]{bookAGV}).

\begin{theorem}[\cite{Fujita82}]\label{ClassdP4}
Let $(M, H_M)$ be a smooth del Pezzo fourfold of degree $d = H_M^4$. Suppose that $\rho (M) = 1$. Then $1 \leqslant d \leqslant 5$ and $M$ is one of the following:
\begin{enumerate}
    \item If $d=1$, then $M = M_6 \subseteq \bP(1^4, 2, 3)$. %is a hypersurface of degree $6$ in the weighted projective space $\bP(1^4, 2, 3)$.
    \item If $d=2$, then $M = M_4 \subseteq \bP(1^5, 2)$. %is a hypersurface of degree $4$ in the weighted projective space $\bP(1^5, 2)$.
    \item If $d=3$, then $M =M_3 \subseteq \bP^5$. %is a cubic hypersurface.
    \item If $d=4$, then $M = M_{2, 2} \subseteq \bP^6$. %is a complete intersection of two quadrics.
    \item If $d=5$, then $M$ is a $2$-codimensional linear section of the Grassmannian $\Gr (2, 5) \subseteq \bP^9$ in the Pl\"ucker embedding.
\end{enumerate}
\end{theorem}

When the Fano fourfold $(M, H_M)$ is Mukai, there is an integer $g \geqslant 2$, called the genus of $M$, such that $h^0(H_M)=g+3$ and $d = 2 g - 2$ (see \cite[Corollary 2.1.14]{bookAGV}). %In this case, $\varphi_{|H|}:M\rightarrow\bP^{g+2}$ is a morphism of degree 1 or 2 and it is an embedding if $g \geqslant 4$.

\begin{theorem}[\cite{Mukai89}] \label{ClassMu4}
Let $(M, H_M)$ be a smooth Mukai fourfold of genus $g$. Suppose that $\rho (M) = 1$. Then $2 \leqslant g \leqslant 10$.
\begin{enumerate}[(I)]
    \item If $2 \leqslant g \leqslant 5$, $M$ is one of the following:
    \begin{itemize}
    \item[($g = 2$)] $M = M_6 \subseteq \bP(1^5,3)$.
    \item[($g = 3$)] $M = M_4 \subseteq \bP^5$ or $M_{2,4} \subseteq \bP(1^6,2)$.
    \item[($g = 4$)] $M = M_{2,3} \subseteq \bP^6$.
    \item[($g = 5$)] $M = M_{2,2,2} \subseteq \bP^7$.
    \end{itemize}

    \item If $6 \leqslant g \leqslant 10$, then $M$ is a $(n(g) - 4)$-codimensional linear section of
    an $n(g)$-dimensional smooth variety
    \[
    \Sigma^{n(g)}_{2g - 2} \subseteq \bP^{g + n(g) - 2}
    \]
    of degree $2g - 2$, which can be described as follows:
    \begin{itemize}
    \item[($g = 6$)] $\Sigma^{6}_{10} = Q_2 \cap CG \subseteq \bP^{10}$ is a quadric section of a cone $CG\subseteq \bP^{10}$ over the Grassmannian $G = \Gr(2, 5) \subseteq \bP^9$ in the Pl\"ucker embedding.

    \item[($g = 7$)] $\Sigma^{10}_{12} = {\rm OG}_+(5,10) \subseteq \bP^{15}$ is a connected component of the orthogonal Grassmannian ${\rm OG} (5,10)$ in the half-spinor embedding.

    \item[($g = 8$)] $\Sigma^8_{14} = \Gr (2, 6) \subseteq \bP^{14}$ is the Grassmannian $\Gr (2, 6)$ in the Pl\"ucker embedding.

    \item[($g = 9$)] $\Sigma^6_{16} = \mathrm{LG} (3, 6) \subseteq \bP^{13}$ is the Lagrangian Grassmannian $\mathrm{LG} (3, 6)$ in the Pl\"ucker embedding.

    \item[($g = 10$)] $\Sigma^5_{18} \subseteq \bP^{13}$ is the subvariety of $\Gr (5, 7)$ parameterizing isotropic $5$-spaces of a general $4$-form in $\bC^7$ in the Pl\"ucker embedding. %Another description is $\Sigma_{18}^5=G_2/P\subseteq\bP^{13}.$
    \end{itemize}
\end{enumerate}
\end{theorem}

\begin{remark}[$g = 6$]
Let $v \in \bP^{10}$ be the vertex of the cone $CG = C\Gr(2,5)$. One can prove that $M = \Gr (2, 5) \cap Q_2 \cap \Lambda$ (under the projection from $v$) if $v$ is not in the $2$-codimensional linear subspace $\Lambda \subseteq \bP^{10}$. For the case $v \in \Lambda$, the linear space $\Lambda$ is a cone over $\bP^7 \cong L \subseteq \bP^9$. Let $W \coloneqq \Gr (2, 5) \cap L$. Then $M$ is the intersection of the cone $CW$ over $W$ with a quadric $Q_2$ and $M \to W$ is a double cover branched along the (smooth) intersection of $W$ with a quadric.
\end{remark}

For the normal bundle of $\Sigma^{n(g)}_{2g - 2}$ in the Grassmannian for $g = 7,9,10$, we make the following remarks.

\begin{remark}[$g = 7$]
Recall that ${\rm OG}(5,10)\subseteq\Gr(5,10)$ is the zero locus of a global section of the vector bundle $\Sym^2(\cS^\vee)$, and ${\rm OG} (5,10)$ is a disjoint union of two isomorphic connected components ${\rm OG}_{\pm}(5,10)$. Hence the fundemental cycles satisfy
\[
2 [\Sigma_{12}^{10}] = [{\rm OG}_+(5,10)] + [{\rm OG}_-(5,10)] = c_{15} (\Sym^2(\cS^\vee)) \cap [\Gr(5, 10)].
\]
Notice that a hyperplane section of
${\rm OG}_+(5,10)$ via the Pl\"ucker embedding is linearly equivalent to
twice the hyperplane section of the half-spinor embedding ${\rm OG}_+(5,10) \hookrightarrow \bP^{15}$ (see \cite[Proposition 1.7]{Mukai95}).

%It is a component of the set of linear spaces of dimension 4 contained in a smooth eight dimensional quadric hypersurface in $\bP^9$ in its spinor embedding\footnote{The spinor embedding is given by $V^5\mapsto(\det(V),\wedge^2V,\wedge^4V)\in\bP(\wedge^{\rm even}V)$}.
\end{remark}

\begin{remark}[$g = 9, 10$]
The Lagrangian Grassmannian ${\rm LG}(3,6)\subseteq\Gr(3,6)$ is the zero locus of a global section of the vector bundle $\wedge^2(\cS^\vee)$, and $\Sigma_{18}^5\subseteq\Gr(5,7)$ is the zero locus of a global section of the vector bundle $\wedge^4(\cS^\vee)$.

%Another description of $\Sigma_{18}^5$ is a $5$-dimensional homogeneous variety $G_2/P_2$ for a group of type $G_2$.
%(see \cite[Example 5.2]{Mukai03}).
\end{remark}

By Theorems \ref{ClassdP4} and \ref{ClassMu4}, our task now is to compute Chern classes of the tangent bundle of a weighted projective space $\bP(\underline{a}) = \bP (a_0, \cdots,a_m)$ and the Grassmannian $G = \Gr (k, n)$. This follows from the generalized Euler exact sequence (see \cite[Theorem 12.1]{BC94})
\[
0\ra\Omega_{\bP(\underline{a})}\ra\bigoplus_{i=0}^m\cO_{\bP(\underline{a})}(-a_i)\ra\cO_{\bP(\underline{a})}\ra0.
\]
and $T_G \cong \cS^{\vee} \otimes \cQ $, where $\cS$ and $\cQ$ are the universal sub- and quotient bundles.

We are now in position to give tables of the Hodge numbers of the obtained Calabi--Yau $X_{\sF}$ in Sections \ref{movfansecI} and \ref{movfansecII} (see the list in Proposition \ref{cicylist}). Recall that we know that $h^{1, 1}(\tY) = 1$ and $h^{1, 1}(X_{\sF}) = 2$. Therefore the Hodge number $h^{2, 1}(\tY)$ could be computed by the standard tools of intersection theory, Schubert calculus (cf.~\cite[Section 14.7]{fulton98}), the above classification results, and the Hodge number $h^{2, 1}(X_{\sF})$ by Proposition \ref{compEuler} and \ref{numbODP}. All obtained results are summarized in Tables \ref{tab:dP4F}--\ref{tab:P4GrEtnot}.

\begin{table}[H]
\begin{center}
\begin{tabular}{ccccc}
\toprule
$\sF$ & $d$ & $\chi_{top} (X_{\sF})$ & $h^{2,1} (X_{\sF})$ & \# of ODPs\\ \midrule
$\cO(1)^3$ & 2 & -132 & 68 & 12\\
& 3 & -108 & 56 & 18\\
& 4 & -96 & 50 & 24\\
& 5 & -90 & 47 & 30\\
$\cO(1) \oplus \cO(2)$ & 2 & -140 & 72 & 8\\
& 3 & -120 & 62 & 12 \\
& 4 & -112 & 58 & 16 \\
& 5 & -110 & 57 & 20\\
\bottomrule
\end{tabular}
\end{center}
   \caption{M is a smooth dP4 with $\rho (M) = 1$.}
    \label{tab:dP4F}
\end{table}

\begin{table}[H]
\begin{center}
\begin{tabular}{ccc}
\toprule
$d$ & $\chi_{top} (\tY)$ & $h^{2,1} (\tY)$ \\
\midrule
2 & -156 & 79  \\
3 & -144 & 73  \\
4 & -144 & 73  \\
5 & -150 & 76  \\
\bottomrule
\end{tabular}
\end{center}
   \caption{M is a smooth dP4 with $\rho (M) = 1$ and
   $\tY \in |- K_M|$ is smooth.}
    \label{tab:dP4Y}
\end{table}

\begin{table}[H]
\begin{center}
\begin{tabular}{cccccc}
\toprule
$g$ & $\chi_{top} (\tY)$ & $h^{2,1} (\tY)$ & $\chi_{top} (X_{\sF})$ & $h^{2,1} (X_{\sF})$ & \# of ODPs\\
\midrule
2 & -256 & 129 & -252 & 128 & 2\\
3 & -176 & 89 & -168 & 86 & 4\\
4 & -144 & 73 & -132 & 68 & 6\\
5 & -128 & 65 & -112 & 58 & 8\\
6 & -120 & 61 & -100 & 52 & 10\\
7 & -116 & 59 & -92 & 48 & 12\\
8 & -116 & 59 & -88 & 46 & 14\\
9 & -116 & 59 & -84 & 44 & 16\\
10& -120 & 61 & -84 & 44 & 18\\
\bottomrule
\end{tabular}
\end{center}
   \caption{M is a smooth Muaki fourfold with $\rho (M) = 1$,
   $\tY \in |- K_M|$ is smooth and $\sF = \cO (1)^2$.}
    \label{tab:Muk}
\end{table}

\begin{table}[H]
\begin{center}
    \begin{tabular}{cccccc}
        \toprule
        M & rk & $\sF$ & $\chi (X_{\sF})$ & $h^{2, 1} (X_{\sF})$ & \# of ODPs \\ \midrule
        $\bP^4$ & 5 & $\cO (1)^{5}$ & -100 & 52 & 50 \\
        % & 4 & $T_{\bP^4}$ & -102 & 53 & 50 \\
        & 4 & $\cO (1)^{3} \oplus \cO (2)^{\phantom{2}}$ & -108 & 56 & 46 \\
        & 3 & $\cO (1)^{\phantom{2}} \oplus \cO (2)^{2}$ & -112 & 58 & 44 \\
        & 3 & $\cO (1)^{2} \oplus \cO (3)^{\phantom{2}}$ & -132 & 68 & 34 \\
        & 2 & $\cO (1) \oplus \cO (4)$ & -168 & 86 & 16 \\
        & 2 & $\cO (2) \oplus \cO (3)$ & -128 & 66 & 36 \\
        $\Gr (2, 4)$ & 4 & $\cO (1)^4 $ & -96 & 50 & 40 \\
        & 3 & $\cO (1)^2 \oplus \cO (2)$ & -108 & 56 & 34 \\
        % & 3 & $\mathbf{E} (2) \oplus \cO (1)$ & -96 & 50 & 41 \\
        & 2 & $\cO (1) \oplus \cO (3)$ & -140 & 72 & 18 \\
        & 2 & $\cO (2) \oplus \cO (2)$ & -112 & 58 & 32 \\
        \bottomrule
    \end{tabular}
\end{center}
    \caption{$M = \bP^4$ or $\Gr (2, 4)$, $\sE = \cO^{\mathrm{rk}}$.}
    \label{tab:P4GrEt}
\end{table}

\begin{table}[H]
\begin{center}
\begin{tabular}{ccccccc}
\toprule
M & rk & $\sE$ & $\sF$ & $\chi (X_{\sF})$ & $h^{2, 1} (X_{\sF})$ & \# of ODPs \\ \midrule
$\bP^4$ & 3 & $\cO^2 \oplus \cO (1)$ & $\cO (1)^2 \oplus \cO (2)$ & -120 & 62 & 40 \\
 & 2 & $\cO^{\phantom{2}} \oplus \cO (1)$ & $\cO (1)^{\phantom{2}} \oplus \cO (3)$ & -152 & 78 & 24 \\
$\Gr (2, 4)$ & 2 & $\cO^{\phantom{2}} \oplus \cO (1)$ & $\cO (1)^{\phantom{2}} \oplus \cO (2)$ & -128 & 66 & 24\\
\bottomrule
\end{tabular}
\end{center}
    \caption{$M = \bP^4$ or $\Gr (2, 4)$, $\sE \neq \cO^{\mathrm{rk}}$.}
    \label{tab:P4GrEtnot}
\end{table}

%%%%%%%%%%%%%%%%%%%%
%%%%%%%%%%%%%%%%%%%%
%%%%%%%%%%%%%%%%%%%%
%%%%%%%%%%%%%%%%%%%%
\bibliographystyle{alpha}%{amsalpha} 
\bibliography{references}
\end{document}